\documentclass[a4paper,fleqn]{cas-sc}

 \usepackage[square, comma, sort&compress, numbers]{natbib}

\def\tsc#1{\csdef{#1}{\textsc{\lowercase{#1}}\xspace}}
\tsc{WGM}
\tsc{QE}
\begin{document}
	\let\WriteBookmarks\relax
	\def\floatpagepagefraction{1}
	\def\textpagefraction{.001}
	
	% Short title
	\shorttitle{A ST GFDM for solving the transient Stokes/Parabolic interface problem in the moving system}    
	
	% Short author
	\shortauthors{Y.N.Xing,H.B.Zheng}  
	
	% Main title of the paper
	\title [mode = title]{A space-time generalized finite difference method for solving the transient Stokes/Parabolic interface problem in the moving system}  
	
	% Title footnote mark
	% eg: \tnotemark[1]
	%\tnotemark[<tnote number>] 
	
	% Title footnote 1.
	% eg: \tnotetext[1]{Title footnote text}
	%\tnotetext[<tnote number>]{<tnote text>} 
	
	% First author
	%
	% Options: Use if required
	% eg: \author[1,3]{Author Name}[type=editor,
	%       style=chinese,
	%       auid=000,
	%       bioid=1,
	%       prefix=Sir,
	%       orcid=0000-0000-0000-0000,
	%       facebook=<facebook id>,
	%       twitter=<twitter id>,
	%       linkedin=<linkedin id>,
	%       gplus=<gplus id>]
	
	\author[a]{Yanan Xing}
	
	% Corresponding author indication
	%\cormark[]
	
	% Footnote of the first author
	%\fnmark[<footnote mark no>]
	
	% Email id of the first author
	%\ead{xingyanan9529@qq.com}
	
	% URL of the first author
%	\ead[url]{<URL>}
	
	% Credit authorship
	% eg: \credit{Conceptualization of this study, Methodology, Software}
	\credit{Methodology, Writing the initial draft and the final draft}
	
	% Address/affiliation
	\affiliation[a]{
		addressline={School of Mathematical Sciences, East China Normal University}, 
		city={shanghai},
		%          citysep={}, % Uncomment if no comma needed between city and postcode
		postcode={200241}, 
		%   state={},
		country={China}}
	
	\author[b]{Haibiao Zheng}
	
	% Footnote of the second author
	\fnmark[*]
	
	% Email id of the second author
	\ead{hbzheng@math.ecnu.edu.cn}
	
	% URL of the second author
	%\ead[url]{}
	
	% Credit authorship
	\credit{Methodology, Revising the paper}
	
	% Address/affiliation
\affiliation[b]{
	addressline={School of Mathematical Sciences, Key Laboratory of MEA(Ministry of Education), Shanghai Key Laboratory of PMMP, East China Normal University}, 
	city={shanghai},
	%          citysep={}, % Uncomment if no comma needed between city and postcode
	postcode={200241}, 
	%   state={},
	country={China}}
% Corresponding author text
\cortext[a]{Corresponding author}

% Footnote text
%\fntext[*]

% For a title note without a number/mark
%\nonumnote{*}{}

% Here goes the abstract
\begin{abstract}
In this paper, a space-time generalized finite difference method (ST-GFDM) is proposed to solve the transient Stokes/Parabolic moving interface problem which is a type of fluid-structure interaction problem. The ST-GFDM considers the time dimension as the third space dimension, and the 2D time-dependent Stokes/Parabolic moving interface problem can be seen as a 3D interface problem where the interface is formed by the initial interface shape and the moving trajectory. The GFDM has an advantage in dealing with interface problems with complex interface shape and moving interface in the ST domain. More irregular moving direction, more complex interface shape, and the translation and deformation of interface are analyzed to show the advantage of the  ST-GFDM. The interface problem can be transformed into coupled sub-problems and  locally dense nodes is used when the subdomain is too small to satisfy the needs of the numbers of the nodes to improve the performance of the ST-GFDM. Five examples are provided to verify the existence of the good performance of the ST-GFDM for Stokes/Parabolic moving interface problems, including those of the simplicity, accuracy, high efficiency and stability.    
\end{abstract}

% Use if graphical abstract is present
%\begin{graphicalabstract}
%\includegraphics{}
%\end{graphicalabstract}

% Research highlights
%\begin{highlights}
%\item The Stokes/Parabolic moving interface problem under PPE form is proposed to overcome the disadvantage of the pressure term.
%\item The irregular moving direction, complex interface shape and the deformation of complex interface shape are adopted to verify the advantage of the ST-GFDM for solving the moving interface problems.
%\item The domain decomposition method transforms the Stokes/Parabolic moving interface problems into two non-interface problems that a Stokes problem in terms of fluid velocity and pressure and a Parabolic problem in terms of the structure velocity.
%\item The comparison between the original GFDM under PPE form, the ST-GFDM, the DLM/FD FEM [11] is provided to show the accuracy, stability and simplicity of these two proposed methods.
%\end{highlights}
% Keywords
% Each keyword is seperated by \sep
\maketitle

% Main text
\section{Introduction}\label{}
The fluid-structure interaction problem is ubiquitous in nature, it is widely used in aeroelasticity, biomechanics and haemodynamics problems[1-3]. The Stokes/Parabolic moving interface problem is a type of linearized fluid-structure interaction problem[4,5]. It is a classical and important task in the study of the fluid-structure interaction.  The fluid is modelled by Stokes equations in terms of fluid velocity and pressure. As we all know, the Stokes equation can cause the pressure oscillation due to the information deletion of the pressure term. To overcome this difficult problem, researchers have paid much attention to how to formulate well-posed Stokes equations with suitable boundary conditions. The Pressure Poisson Equation (PPE) [6-8] is a popular way.

To deal with the fluid-structure interaction problem has two difficulties: one is these two sub-problems defined in different coordinate system, the interface condition is difficult to handle; another is that this is a moving interface problem, which involves the mesh shape change. In order to overcome these difficulties, many methods coupled with some approaches are proposed to deal with the troublesome mesh reconstruction and coordinate system transformation. Immersed Boundary Method (IBM) and Arbitrary Lagrangian-Eulerian (ALE) approach are two popular methods.

In the past few years, there are many mesh methods have been proposed to analyze the fluid-structure interaction problem. For instance, Chen et al.[9] analyze the Eulerian-Lagrangian flow-vegetation interaction model by using the immersed boundary method and Open FOAM.  Ryzhakov et al.[10] solve the fluid-structure interaction problems involving flows in flexible channels by using a unified arbitrary Lagrangian-Eulerian model. In particular, Sun[11] used the distributed Lagrange multiplier/fictitious domain finite element method to solve the Stokes/Parabolic interface  problem with jump coefficients. Lan et al.[12,13] proposed a monolithic arbitrary Lagrangian-Eulerian finite element and ALE finite element methods to analyze Stokes/Parabolic moving interface problems. Kesler et al.[14] proposed an arbitrary Lagrangian-Eulerian (ALE) finite element method to solve the transient Stokes/Parabolic interface problems.  From these paper, we find that the majority of these methods are the FEM coupled with some approach to deal with this problem on different coordinate system. Meanwhile, it also faces the mesh construction and numerical integration, especially for the changing interface shape.

GFDM, which is based on Taylor series expansions and weighted moving square (MLS) approximation, is truly free from mesh generation and numerical quadrature.  Recently, the GFDM has been extended to study various problems, such as the inverse problem [15] and the interface problem in static and moving systems [16-18], two-phase flows[19], phase transitions[20] and so on. Later, the meshless method was coupled with some other methods to deal with the time-dependent problems, such as the Krylov deferred correction (KDC) method [21,22], the space-time approach [23-28] and the second-order Crank-Nicolson scheme[29] and so on. The main idea of the space-time approach is that it considers the time dimension as the another space dimension in the ST domain and avoids the troublesome time discretization. The ST-GFDM has been successfully used in many problems, especially for the Parabolic PDEs[26], Hyperbolic PDEs [27] and the Zakharov-Kuznetsov-Modified Equal Widthequation[28].In this paper, we adopt the ST-GFDM to solve the Stokes/parabolic moving interface problem.

The ST-GFDM is a novel meshless method and the point collocation on the same Cartesian coordinates, it has advantage in dealing with the interface condition. The moving interface and the complexity of the interface shape only affect the point collocation on the interface, the ST-GFDM has advantage in dealing with the moving interface problems. The ST-GFDM approximates the partial derivatives of unknown functions of different orders by using a linear combination of adjacent node function values, the ST-GFDM has advantage in dealing with the term of time partial derivatives. The ST-GFDM considers the 2D time-dependent problem as a 3D problem and avoids the inconvenience of the time discretization. The Stokes-Parabolic moving interface problem can be transformed into coupled sub-problems and  locally dense nodes is used when the subdomain is too small to satisify the needs of the numbers of the nodes to improve the performance of the ST-GFDM.  In this paper, some Stokes-Parabolic moving interface problems with more irregular moving direction, more complex interface shape ,and the translation and deformation of interface are analyzed to show the advantage of the ST-GFDM.

The rest of the paper is structured as follows: Section 2 introduces the model of Stokes/Parabolic moving interface problems. Section 3 presents the procedure of ST-GFDM,  the handling skill for the Stokes/ Parabolic moving interface problem, and the ST-GFDM for Stokes/Parabolic moving interface problems. In Section 4, five numerical examples are presented to verify the accuracy, high efficiency and stability of the proposed method. Finally, a conclusion is given in Section 5.
\section{The Stokes/Parabolic moving interface problem}
In this paper, we consider the following transient Stokes/Parabolic moving interface problem (From Ref.[11]):
\begin{eqnarray}
	\rho_1 \frac{\partial \mathbf{u}_1}{\partial t}-\nabla \cdot(\beta_1\nabla \mathbf{u}_1)+\nabla p_1&=&\mathbf{f}_1, \quad  in \ \Omega^1_t\times (0,T],\\
	\nabla \cdot \mathbf{u}_1&=&0, \quad \ \ in \ \Omega^1_t\times (0,T],\\
	\rho_2 \frac{\partial \mathbf{u}_2}{\partial t}-\nabla \cdot(\beta_2\nabla \mathbf{u}_2)&=&\mathbf{f}_2, \quad  in \ \Omega^2_t \times (0,T],\\
	\mathbf{u}_1&=&\mathbf{u}_2,  \quad on \ \Gamma_t \times (0,T],\\
	(\beta_1\nabla \mathbf{u}_1-p_1\mathbf{I})n_1+\beta_2\nabla \mathbf{u}_2 n_2&=& \mathbf{\tau},  \quad \ \  on \ \Gamma_t \times (0,T],\\
	\mathbf{u}_1&=&0,  \ \ \quad on\ \partial \Omega_t^1 \textbackslash \Gamma_t \times (0,T],\\
	\mathbf{u}_2&=&0, \ \ \quad on\ \partial \Omega_t^2 \textbackslash \Gamma_t \times (0,T],\\
	\mathbf{u}_1(\mathbf{x},0)&=&\mathbf{u}_1^0,  \quad in\ \Omega_0^1,\\
	\mathbf{u}_2(\mathbf{x},0)&=&\mathbf{u}_2^0,  \quad in\ \Omega_0^2,
\end{eqnarray}
with discontinuous coefficient:
\begin{equation}
	\beta=\begin{cases}
		{\beta_1},&\mbox{in $\Omega_t^1$},\\			
		{\beta_2},&\mbox{in $\Omega_t^2$},
		\end {cases}
		\rho=\begin{cases}
			{\rho_1},&\mbox{in $\Omega_t^1$},\\			
			{\rho_2},&\mbox{in $\Omega_t^2$},
			\end {cases}
		\end{equation}
		where $\beta_1$ and $\rho_1$ can stand for the fluid viscosity and density of the fluid, and $\beta_2$ and $\rho_2$ for the elastic parameter and the density of the structure. The vectors $n_1=(n_{11},n_{12})^T,n_2=(n_{21},n_{22})^T$ denote the unit outward normal vectors of the interface.  $ \Gamma_t $ denotes a moving interface, $\mathbf{u}_1^0=(u_1^0,v_1^0)^T$ and $\mathbf{u}_2^0=(u_2^0,v_2^0)^T$ are the initial conditions.  $\mathbf{u}_1=(u_1,v_1)^T$ is in terms of fluid velocity, $p_1$ is in terms of fluid pressure, $\mathbf{u}_2=(u_2,v_2)^T$ is in terms of the structure velocity. $\mathbf{\tau}=(\tau_1,\tau_2)^T$ is a function related to the $\beta$, $\mathbf{f}_1=(f_{11},f_{12})^T$ and $\mathbf{f}_2=(f_{21},f_{22})^T$ is the source function. From the above Stokes/Parabolic moving interface problem (see Fig.1), we can see that Eq.(2) doesn't have any information about the pressure $p_1$. Therefore, we use the classical pressure poisson equation from [30], which finds divergence on both sides of the equal sign for Eq.(1), then
		\begin{eqnarray} 
			\rho_1 \nabla \cdot\frac{\partial \mathbf{u}_1}{\partial t}-\nabla \cdot\nabla \cdot(\beta_1\nabla\mathbf{u}_1)+\nabla \cdot\nabla p_1&=&\nabla \cdot\mathbf{f}_1, \  \  in \ \Omega_t^1\times(0,T],
		\end{eqnarray}
		we rewrite the above equation
		\begin{eqnarray} 
			\rho_1 \nabla \cdot\frac{\partial \mathbf{u}_1}{\partial t}-\nabla \cdot\nabla \cdot(\beta_1\nabla \mathbf{u}_1)+\Delta p_1&=&\nabla \cdot \mathbf{f}_1,\  \  in \ \Omega_t^1 \times(0,T].
		\end{eqnarray}
		Due to Eq.(2), the above equation can be simplified as follows
		\begin{eqnarray}
			\Delta p_1 &=&\nabla \cdot \mathbf{f}_1,\quad \quad in \ \Omega_t^1 \times(0,T],
		\end{eqnarray}
		Therefore, the Eq.(2) can be exchange into the following part
		\begin{eqnarray}
			\Delta p_1 &=&\nabla \cdot \mathbf{f}_1,\quad \quad in \ \Omega_t^1\times(0,T],	\end{eqnarray}
		\begin{eqnarray}
			\nabla \cdot \mathbf{u}_1&=&0,\quad\quad\quad\quad on \  \partial \Omega_t^1 \times(0,T]. 
		\end{eqnarray}
		In the boundary part of the above equations, we adopt the following scheme
		\begin{eqnarray}
			\nabla \cdot \mathbf{u}_1+p_1-p_1&=&0,\quad\quad\quad\quad on \  \partial \Omega_t^1 \times(0,T], 
		\end{eqnarray}
		then a new boundary condition can be obtained
		\begin{eqnarray}
			\nabla \cdot \mathbf{u}_1+p_1&=&p_1,\quad\quad\quad\quad on \  \partial \Omega_t^1 \times(0,T].
		\end{eqnarray}
		\begin{figure}
			\centering
			\includegraphics[scale=.5]{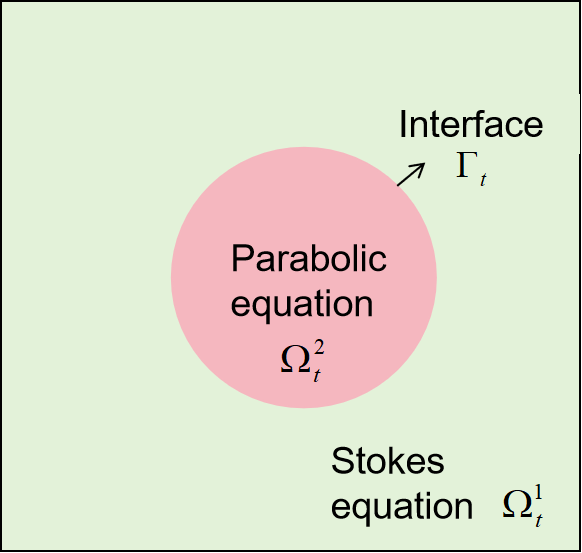}    
			\includegraphics[scale=.5]{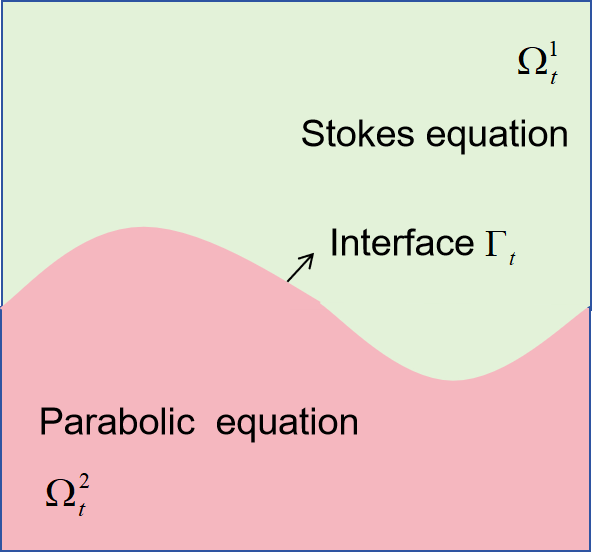}    
			\caption{ The distributions of the Stokes/Parabolic moving interface problem with closed interface(left) and unclosed interface(right).}
		\end{figure}
% Uncomment and use as the case may be
%\begin{theorem} 
%\end{theorem}

% Uncomment and use as the case may be
%\begin{lemma} 
%\end{lemma}

%% The Appendices part is started with the command \appendix;
%% appendix sections are then done as normal sections
%% \appendix

\section{Numerical schemes}
In this section, we briefly describe the numerical scheme of the ST-GFDM. The time-dependent two-dimensional (2D) problem is transformed into a three-dimensional (3D) problem in the $x-y-t$ ST domain. Specifically, for 2D time-dependent interface problems (see Fig. 2), the domain is a 2D domain at all times. We store the domain information at all times, a 3D ST domain can be formed. The interface of the 3D ST domain can be determined by the initial interface shape or the moving direction. The 3D ST domain can be divided by $\Omega^+,\Omega^-$ and $\Gamma$, where $\Omega^+=\Omega_0^+\cup\Omega_1^+\cup...\cup\Omega_i^+$, $\Omega^-=\Omega_0^-\cup\Omega_1^-\cup...\cup\Omega_i^-$, $\Gamma=\Gamma_0\cup\Gamma_1\cup...\cup\Gamma_i,$ $i=0,1,...,n$.
% Figure
\begin{figure}
	\centering
	\includegraphics[scale=.5]{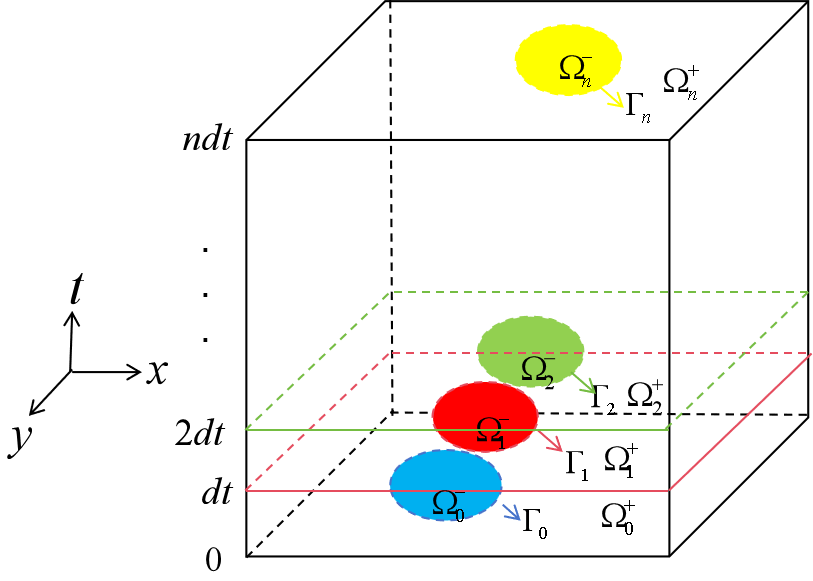}
	\caption{ The process of the 2D spatial domain transformed into 3D ST-domain.}
\end{figure}

\begin{figure}
	\centering
	\includegraphics[scale=0.7]{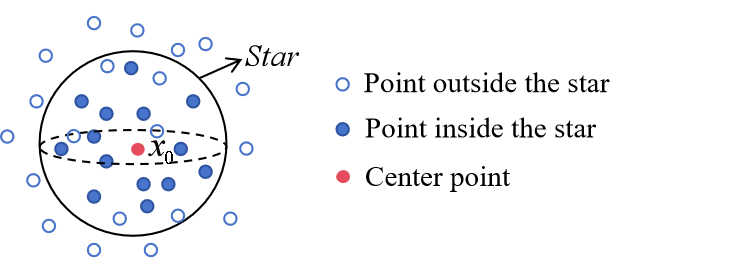}    
	
	\caption{ A star with the center point $x_0$ of the ST-GFDM.}
\end{figure}
\subsection{The Space-Time generalized finite difference method}

The main idea of the ST-GFDM is that the computational domain of the 2D  time-dependent moving interface problem is set in an $x-y-t$ ST domain. For a given interior node $(x_0,y_0,t_0)$, we select $m$ nearest nodes $(x_k,y_k,t_k)(k=1,2,...,m)$ around the central node $(x_0,y_0,t_0)$ and make these $m+1$nodes from a star, as shown in Fig.3. Let $u_0=u(x_0,y_0,t_0)$ is the function value at the centre node $(x_0,y_0,t_0)$ and $u_k=u(x_k,y_k,t_k) (k=1,2,...,m)$ are the function values of the rest nodes $(x_k,y_k,t_k) (k=1,2,...,m)$ inside the star. To exspand these function values $u_k=u(x_k,y_k,t_k)$ around the central node
$(x_0,y_0,t_0)$ by using Taylor serise
\begin{eqnarray}
	u_k=u_0+h_k\frac{\partial u_0}{\partial x}+l_k\frac{\partial u_0}{\partial y}+r_k\frac{\partial u_0}{\partial t}+\frac{1}{2}{h_k}^2\frac{\partial^{2}u_0}{\partial x^{2}}+\frac{1}{2}{l_k}^2\frac{\partial^{2}u_0}{\partial y^{2}}+\frac{1}{2}{r_k}^2\frac{\partial^{2}u_0}{\partial t^{2}}
	+{h_k}{l_k}\frac{\partial^{2}  u_0}{\partial x \partial y}+{h_k}{r_k}\frac{\partial^{2}u_0}{\partial x \partial t}+{l_k}{r_k}\frac{\partial^{2}u_0}{\partial y \partial t}				
\end{eqnarray}
$$+O(\rho^3).$$
To truncate the Taylor series after the second-order derivatives, a weighted residual functional $B(u)$ can be defined as
\begin{equation}
	\begin{aligned}
		B(u)=&\left(\sum_{k=1}^m [(u_0-u_k+h_k\frac{\partial u_0}{\partial x}+l_k\frac{\partial u_0}{\partial y}+r_k\frac{\partial u_0}{\partial t}+\frac{1}{2}{h_k}^2\frac{\partial^{2}u_0}{\partial x^{2}}+\frac{1}{2}{l_k}^2\frac{\partial^{2}u_0}{\partial y^{2}}+\frac{1}{2}{r_k}^2\frac{\partial^{2}u_0}{\partial t^{2}}+{h_k}{l_k}\frac{\partial  u_0}{\partial x \partial y}\right.\\
		&\left.{{+{h_k}{r_k}\frac{\partial^{2}u_0}{\partial x \partial t}+{l_k}{r_k}\frac{\partial^{2}u_0}{\partial y \partial t})\omega_k]^2}}\right),
	\end{aligned}
\end{equation}
where $\ h_k$, $\ l_k$ and $\ r_k$ are the distances between the nodes$(x_k,y_k,t_k)$and the central node $(x_0,y_0,t_0)$ along $\ x$ axis , $\ y$ axis and $\ t$ axis, respectively,
$	h_k=x_k-x_0, l_k=y_k-y_0, r_k=t_k-t_0,  k=1,2,...,m.$
$\omega_k (k=1,2,...,m)$denote the weighting coefficients at the nodes $(x_k,y_k,t_k)$. We adopted the following weighting function:
\begin{eqnarray}
	\omega_k=
	\begin{cases}
		1-6(\frac{d_k}{d_m})^2+8(\frac{d_k}{d_m})^3-3(\frac{d_k}{d_m})^4&d_k\leq d_m,\\
		0&d_k\geq d_m.\\
	\end{cases}
\end{eqnarray}
Where $d_k=\sqrt{(x_k-x_0)^2+(y_k-y_0)^2+(t_k-t_0)^2}$ is the distance between nodes $(x_k,y_k,t_k)$ and $(x_0,y_0,t_0)$, and $d_m$ denotes the maximum value of $d_k(k=1,2,...,m)$. The further away from the centre point $(x_0, y_0,t_0)$, the smaller the value of the weight function. The weighting function is used to indicate that the approximation $u_k$ is more important when the nodes $(x_k,y_k,t_k)(k = 1, 2, …, m)$ are more closer to the central node $(x_0,y_0,t_0)$ in the star.
To minimize the above residual function $B(u)$ with respect to the partial derivatives
$ D_u=(\frac{\partial u_0}{\partial x}, \frac{\partial u_0}{\partial y}, \frac{\partial u_0}{\partial t}, \frac{\partial^{2}u_0}{\partial x^{2}},\frac{\partial^{2}u_0}{\partial y^{2}}, \frac{\partial^{2}u_0}{\partial t^{2}},\frac{\partial^{2}u_0}{\partial x \partial y },\frac{\partial^{2}u_0}{\partial x \partial t },\frac{\partial^{2}u_0}{\partial y \partial t })^T,$
that is, let$	\frac{\partial B(u)}{\partial D_u}=0,$
then the following system of linear equations can be obtained
\begin{equation}
	AD_u=b. 
\end{equation}
where $A$ is a symmetric matrix given by
\begin{equation}
	A=\sum_{k=1}^m \left( 
	\begin{array}{lllllllll}
		{{h_k^2}} &{h_kl_k} &{h_kr_k} &{\frac{h_k^3}{2}} &{\frac{h_kl_k^2}{2}} &{\frac{h_kr_k^2}{2}} &{h_k^2l_k} &{h_k^2r_k} &{h_kl_kr_k} \\
		{h_kl_k} &{l_k^2}&{l_kr_k} &{\frac{h_k^2l_k}{2}} &{\frac{l_k^3}{2}} &{\frac{l_kr_k^2}{2}} &{h_kl_k^2} &{h_kl_kr_k} &{l_k^2r_k} \\
		{h_kr_k} &{l_kr_k} &{r_k^2} &{\frac{h_k^2r_k}{2}} &{\frac{l_k^2r_k}{2}} &{\frac{r_k^3}{2}} &{h_kl_kr_k} &{h_kr_k^2} &{l_k^2r_k}
		\\
		{\frac{h_k^3}{2}} &{\frac{h_k^2l_k}{2}} &{\frac{h_k^2r_k}{2}}  &{\frac{h_k^4}{4}} &{\frac{h_k^2l_k}{4}} &{\frac{h_k^2r_k^2}{4}} &{\frac{h_k^3l_k}{2}} &{\frac{h_k^2r_k}{2}}
		&{\frac{h_k^2l_k}{2}}
		\\
		{\frac{h_kl_k^2}{2}} &{\frac{l_k^3}{2}} &{\frac{l_k^2r_k}{2}} &{\frac{h_k^2l_k^2}{4}}  &{\frac{l_k^4}{4}} &{\frac{l_k^2r_k^2}{4}} &{\frac{h_kl_k^3}{2}} &{\frac{h_kl_k^2r_k}{2}}&{\frac{l_k^3r_k}{2}}
		\\
		{\frac{h_kr_k^2}{2}} &{\frac{l_kr_k^2}{2}} &{\frac{r_k^3}{2}} &{\frac{h_k^2r_k^2}{4}} &{\frac{l_k^2r_k^2}{4}} &{\frac{r_k^4}{4}} &{\frac{h_kl_kr_k^2}{2}} &{\frac{h_kr_k^3}{2}} &{\frac{l_kr_k^3}{2}}
		\\
		{h_k^2l_k} &{h_kl_k^2} &{h_kl_kr_k} &{\frac{h_k^3l_k}{2}} &{\frac{h_kl_k^3}{2}} &{\frac{h_kl_kr_k^2}{2}} &{h_k^2l_k^2} &{h_k^2l_kr_k} &{h_kl_k^2r_k}\\
		{h_k^2r_k} &{h_kl_kr_k} &{h_kr_k^2} &{\frac{h_k^2r_k}{2}} &{\frac{h_kl_k^2r_k}{2}} &{\frac{h_kr_k^3}{2}} &{h_k^2l_kr_k} &{h_k^2r_k^2} &{h_kl_kr_k^2}
		\\
		{h_kl_kr_k} &{l_k^2r_k} &{l_k^2r_k} &{\frac{h_k^2l_k}{2}} &{\frac{l_k^3r_k}{2}} &{\frac{l_kr_k^3}{2}} &{h_kl_k^2r_k} &{h_kl_kr_k^2} &{l_k^2r_k^2}
	\end{array}
	\right) \omega_k^2.
\end{equation}
\begin{equation}
	A=\sum_{k=1}^m\omega_k \left( 
	\begin{array}{lllllllll}
		{{h_k}}\\
		{l_k}\\
		{r_k}\\
		{\frac{h_k^2}{2}}\\
		{\frac{l_k^2}{2}}\\
		{\frac{r_k^2}{2}}\\
		{h_kl_k}\\
		{h_kr_k}\\
		{l_kr_k}
	\end{array}
	\right)
	\left( 
	\begin{array}{lllllllll}
		{h_k},
		{l_k},
		{r_k},
		{\frac{h_k^2}{2}},
		{\frac{l_k^2}{2}},
		{\frac{r_k^2}{2}},
		{h_kl_k},
		{h_kr_k},
		{l_kr_k}
	\end{array}
	\right) \omega_k=\sum_{k=1}^mS_k^TS_k.
\end{equation}
Here $S_k=\left( 
\begin{array}{lllllllll}
	{{h_k}},
	{l_k},
	{r_k},
	{\frac{h_k^2}{2}},
	{\frac{l_k^2}{2}},
	{\frac{r_k^2}{2}},
	{h_kl_k},
	{h_kr_k},
	{l_kr_k}
\end{array}
\right) \omega_k,$
and
\begin{equation}
	b=BU=(-\sum_{k=1}^m \omega_kS_{k},\omega_1S_{1},\omega_2S_{2},\cdots,\omega_mS_{m})_{9\times(m+1)}(u_0,u_1,\cdots,u_m)^T,
\end{equation}
in which $U=(u_0,u_1,\cdots,u_m)^T$ are function values of all nodes inside the star. According to Eq. (21) and Eq. (24), the partial derivative vector $D_u$ can be expressed as follows:
\begin{equation}
	D_u=A^{-1}b=A^{-1}(BU)=(A^{-1}B)U=EU,
\end{equation}
where
\begin{equation}
	E=A^{-1}B.
\end{equation}
To repeat the same procedure for each node inside the ST domain, we can obtain the algebraic equations for other interior and boundary points. Enforcing the satisfaction of the governing equations at the interior nodes and the boundary conditions at the boundary nodes, the final system of linear algebraic equations will be yielded. Once this algebraic equation system is solved, the values of all unknown variables will be obtained.
\subsection{The handling skills for the Stokes/ Parabolic moving interface problem}
In this section, the original transient Stokes/Parabolic moving interface problem (Eqs.(1),(3)-(9),(14),(17)) is divided into two two non-interface subproblems:
\begin{eqnarray}
	\rho_1 \frac{\partial \mathbf{u}_1}{\partial t}-\nabla \cdot(\beta_1\nabla \mathbf{u}_1)+\nabla p_1&=&\mathbf{f}_1, \ \quad \quad \quad  in \ \Omega_t^1\times(0,T],\\
	\Delta p_1 &=&\nabla \cdot f_1,\quad \quad in \ \Omega_t^1\times(0,T],\\
	u_1(\mathbf{x},0)&=&\mathbf{u}_1^0, \ \ \quad \quad \quad in\ \Omega_0^1,\\
	\nabla \cdot \mathbf{u}_1+p_1&=&p_1,\ \  \ \  \quad\quad\quad on \  \partial \Omega_t^1\times(0,T], \\
	\mathbf{u}_1&=&0,\ \  \ \  \quad\quad \quad on\  \partial \Omega_t^1 \textbackslash \Gamma_t\times(0,T], \\
	\mathbf{u}_1&=&\mathbf{u}_2, \  \ \ \quad\quad \quad on \ \Gamma_t\times(0,T],
\end{eqnarray}
and 
\begin{eqnarray}
	\rho_2 \frac{\partial \mathbf{u}_2}{\partial t}-\nabla \cdot(\beta_2\nabla \mathbf{u}_2)&=&\mathbf{f}_2,  \ \quad \quad  \quad   in \ \Omega_t^2\times(0,T],\\
	u_2(\mathbf{x},0)&=&\mathbf{u}_2^0,\ \quad \quad  \quad in\ \Omega_0^2,\\
	\mathbf{u}_2&=&0, \ \ \quad \quad \quad on\  \partial \Omega_t^2 \textbackslash \Gamma_t\times(0,T],\\
	(\beta_1\nabla\mathbf{u}_1-p_1\mathbf{I})n_1+\beta_2\nabla \mathbf{u}_2n_2&=& \mathbf{\tau}, \  \quad  \quad \quad \  on \ \Gamma_t\times(0,T].
\end{eqnarray}

\begin{figure}
	\centering
	\includegraphics[scale=.4]{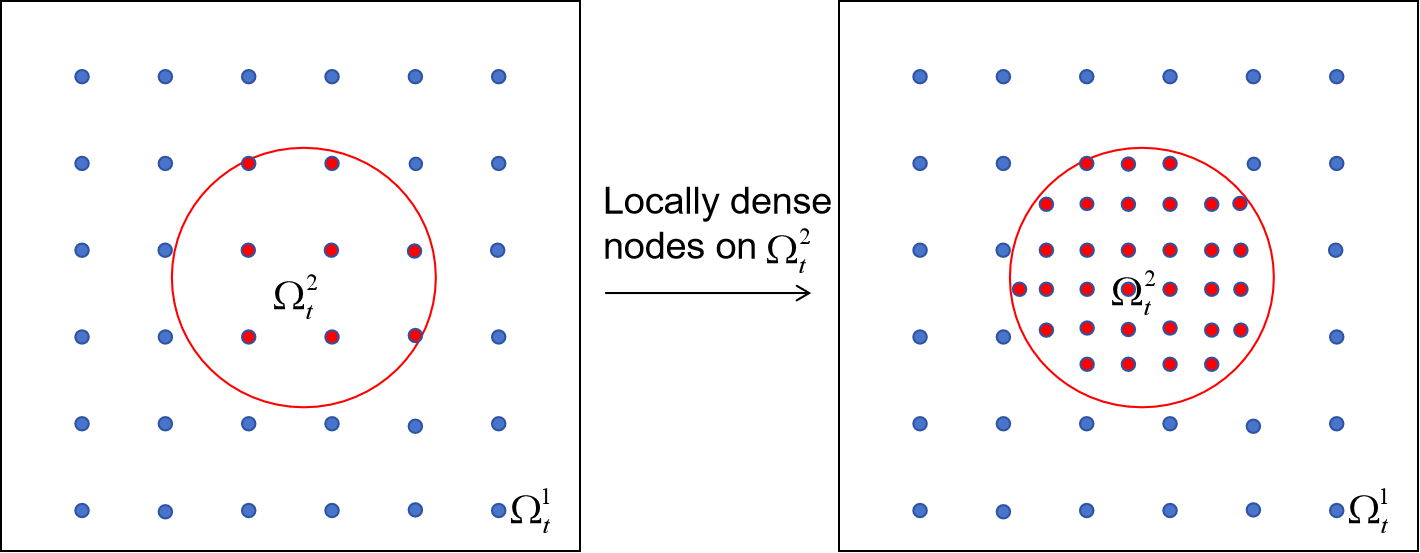}    
	\caption{ The method of Locally dense nodes.}
\end{figure}
Note that the interface can be seen as the boundary of the subproblem and we can use the ST-GFDM twice to solving this Stokes/Parabolic moving interface problem. In addition, the computational subdomains $\Omega_t^1$ and  $\Omega_t^2$ may not have the same shape. Sometimes, the number of nodes inside the subdomain can't meet the computational needs when the subdomain is small, but the number of nodes inside another subdomain that can meet the computational needs. To dense nodes may increase the computational cost.  Therefore, a method of locally dense nodes(see Fig.4) is proposed to overcome the difficulty which we face in Example 1 and Example 2 in this paper. In which, the locally dense nodes are those nodes which use the locally dense nodes method to add many more points in some small subdomain that the number of the nodes can’t satisfy the needs of the computational simulation.

\subsection{The Space-Time generalized finite difference method for the Stokes/Parabolic moving interface problem}
The Stokes/Parabolic moving interface problem can be expressed as the following scheme,
\begin{eqnarray}
	\rho_1 \frac{\partial u_1}{\partial t}-\beta_1 \frac{\partial^2 u_1}{\partial x^2}-\beta_1\frac{\partial^2 u_1}{\partial y^2}+\frac{\partial p_1}{\partial x}&=&f_{11}, \ \quad in \ \Omega_t^1\times(0,T],\\
	\rho_1 \frac{\partial v_1}{\partial t}-\beta_1 \frac{\partial^2 v_1}{\partial x^2}-\beta_1\frac{\partial^2 v_1}{\partial y^2}+\frac{\partial p_1}{\partial y}&=&f_{12},   \  \quad in \  \Omega_t\times(0,T],\\
	\frac{\partial^2 p_1}{\partial x^2}+\frac{\partial^2 p_1}{\partial y^2} &=&\frac{\partial f_{11}}{\partial x}+\frac{\partial f_{12}}{\partial y}, \quad in \ \Omega_t^1\times(0,T],\\
	u_1(\mathbf{x},0)&=&u_1^0, \quad \quad in\  \Omega_t^1, \\
	v_1(\mathbf{x},0)&=&v_1^0,  \quad \quad in\ \Omega_t^1, \\
	\frac{\partial u_1}{\partial x}+\frac{\partial v_1}{\partial y}+p_1&=&p_1,  \quad\quad on \  \partial \Omega_t^1 \times(0,T], \\
	u_1&=&0, \quad\quad  on\ \partial\Omega_0^1  \textbackslash \Gamma_t\times(0,T],\\
	v_1&=&0,  \quad\quad   on \partial\Omega_0^1  \textbackslash \Gamma_t\times(0,T],\\
	u_1&=&u_2,   \quad \quad on \ \Gamma_t\times(0,T],\\
	v_1&=&v_2, \quad \quad on \ \Gamma_t\times(0,T].
\end{eqnarray}
and 
\begin{eqnarray}
	\rho_2 \frac{\partial u_2}{\partial t}-\beta_2 \frac{\partial^2 u_2}{\partial x^2}-\beta_2\frac{\partial^2 u_2}{\partial y^2}&=&f_{21}, \quad \quad    in \ \Omega_t^2\times(0,T],\\
	\rho_2 \frac{\partial v_2}{\partial t}-\beta_2 \frac{\partial^2 v_2}{\partial x^2}-\beta_2\frac{\partial^2 v_2}{\partial y^2}&=&f_{22},   \quad \quad    in \ \Omega_t^2\times(0,T],\\
	u_2(\mathbf{x},0)&=&u_2^0,\ \quad \quad in\ \Omega_0^2(t),\\
	v_2(\mathbf{x},0)&=&v_2^0, \ \quad \quad   in\ \Omega_0^2(t),\\
	u_2&=&0, \ \ \quad \quad on\ \partial\Omega_0^2  \textbackslash \Gamma_t \times(0,T],\\
	v_2&=&0, \  \ \quad \quad on\  \partial\Omega_0^2   \textbackslash \Gamma_t\times(0,T],\\
	(\beta_1\frac{\partial u_1}{\partial x}+p_1)n_{11}+\beta_1\frac{\partial v_1}{\partial x}n_{12}+\beta_2\frac{\partial u_2}{\partial x}n_{21}+\beta_2\frac{\partial v_2}{\partial x}n_{22}&=& \tau_1,  \quad\quad on \ \Gamma_t\times(0,T],\\
	(\beta_1\frac{\partial v_1}{\partial y}+p_1)n_{21}+\beta_1\frac{\partial u_1}{\partial y}n_{11}+\beta_2\frac{\partial u_2}{\partial y}n_{12}+\beta_2\frac{\partial v_2}{\partial y}n_{22}&=& \tau_2, \quad \quad on \ \Gamma_t\times(0,T].
\end{eqnarray}
To enforce interior nodes satisfy the governing equations in $\Omega_t^1$ and $\Omega_t^2$ , boundary nodes satisfy the boundary condition on $\partial \Omega_t=\partial \Omega_t^1\cup \partial\Omega_t^2$, interface nodes satisfy two interface conditions at $\Gamma_t$ , $N_{1,inp}+N_{2,inp}+N_{1,bp}+N_{2,bp}+N_{1,\Gamma_t}+N_{2,\Gamma_t}$  linear algebraic equations can be obtained, and the final linear system can be rewritten in matrix form:
\begin{equation}
	KX=F,
\end{equation}
where
\begin{eqnarray}
	K=\left(
	\begin{array}{lll}
		{M}&{Y}&{C_1} \\
		{L}&{Q}&{C_2}\\
		{G_1}&{G_2}&{P}
	\end{array}
	\right),
	X=\left(
	\begin{array}{ll}
		U\\
		V\\
		P_1
	\end{array}
	\right),
	F=\left (
	\begin{array}{ll}
		D\\
		H\\
		F_p
	\end{array}
	\right),
\end{eqnarray}
here
\begin{eqnarray}
	M=\left(
	\begin{array}{l|l}
		{M_{11}}&{M_{12}} \\
		\hline
		{M_{21}}&{M_{22}}
	\end{array}
	\right),
	Y=\left(
	\begin{array}{ll}
		{0}&{0} \\
		{Y_{21}}&{Y_{22}}
	\end{array}
	\right),
	C_1=\left(
	\begin{array}{l}
		{C_{11}} \\
		{0}
	\end{array}
	\right),
\end{eqnarray}
\begin{eqnarray}
	L=\left(
	\begin{array}{ll}
		{0}&{0} \\
		{L_{21}}&{L_{22}}
	\end{array}
	\right),
	Q=\left(
	\begin{array}{l|l}
		{Q_{11}}&{Q_{12}} \\
		\hline
		{Q_{21}}&{Q_{22}}
	\end{array}
	\right),
	C_2=\left(
	\begin{array}{l}
		{C_{21}} \\
		{0}
	\end{array}
	\right),
\end{eqnarray}
\begin{eqnarray}
	G_1=\left(
	\begin{array}{ll}
		{G_{11}}&{0}
	\end{array}
	\right),
	G_2=\left(
	\begin{array}{ll}
		{G_{21}}&{0}
	\end{array}
	\right),
\end{eqnarray}
where
\begin{eqnarray}
	M_{11}=\left(
	\begin{array}{lll}
		{E_{N_{1,bp}\times{N_{1,bp}}}} &{0}&{0}\\
		{0} &{E_{N_{1,\Gamma_t}\times{N_{1,\Gamma_t}}}} &{0}\\
		{0} &{0} &{M^1_{N_{1,inp}\times N_{1,inp}}}\\
	\end{array}
	\right),
	M_{12}=\left(
	\begin{array}{lll}
		{0} &{0}&{0}\\
		{0} &{-E_{N_{2,\Gamma_t}\times{N_{2,\Gamma_t}}}} &{0}\\
		{0} &{0} &{0}\\
	\end{array}
	\right),
\end{eqnarray}

\begin{eqnarray}
	M_{21}=\left(
	\begin{array}{lll}
		{0} &{0}&{0}\\
		{0} &{I^1_{N_{1,\Gamma_t}\times{N_{1,\Gamma_t}}}} &{0}\\
		{0} &{0} &{0}\\
	\end{array}
	\right),
	M_{22}=\left(
	\begin{array}{lll}
		{E_{N_{2,bp}\times{N_{2,bp}}}} &{0}&{0}\\
		{0} &{I^2_{N_{2,\Gamma_t}\times{N_{2,\Gamma_t}}}} &{0}\\
		{0} &{0} &{M^2_{N_{2,inp}\times N_{2,inp}}}\\
	\end{array}
	\right),
\end{eqnarray}
\begin{eqnarray}
	Y_{21}=\left(
	\begin{array}{lll}
		{0} &{0}&{0}\\
		{0} &{I^3_{N^+_{1,\Gamma_t}\times{N_{1,\Gamma_t}}}} &{0}\\
		{0} &{0} &{0}\\
	\end{array}
	\right),
	Y_{22}=\left(
	\begin{array}{lll}
		{0} &{0}&{0}\\
		{0} &{I^4_{N_{2,\Gamma_t}\times{N_{2,\Gamma_t}}}} &{0}\\
		{0} &{0} &{0}\\
	\end{array}
	\right),
	C_{11}=\left(
	\begin{array}{lll}
		{0} &{0}&{0}\\
		{0} &{I^5_{N^+_{1,\Gamma_t}\times{N_{1,\Gamma_t}}}} &{0}\\
		{0} &{0} &{P^1_{N_{1,inp}\times N_{1,inp}}}\\
	\end{array}
	\right),
\end{eqnarray}
\begin{eqnarray}
	L_{21}=\left(
	\begin{array}{lll}
		{0} &{0}&{0}\\
		{0} &{I^6_{N^+_{1,\Gamma_t}\times{N_{1,\Gamma_t}}}} &{0}\\
		{0} &{0} &{0}\\
	\end{array}
	\right),
	L_{22}=\left(
	\begin{array}{lll}
		{0} &{0}&{0}\\
		{0} &{I^7_{N_{2,\Gamma_t}\times{N_{2,\Gamma_t}}}} &{0}\\
		{0} &{0} &{0}\\
	\end{array}
	\right),
\end{eqnarray}
\begin{eqnarray}
	Q_{11}=\left(
	\begin{array}{lll}
		{E_{N_{1,bp}\times{N_{1,bp}}}} &{0}&{0}\\
		{0} &{E_{N_{1,\Gamma_t}\times{N_{1,\Gamma_t}}}} &{0}\\
		{0} &{0} &{Q^1_{N_{1,inp}\times N_{1,inp}}}\\
	\end{array}
	\right),
	Q_{12}=M_{12},
\end{eqnarray}
\begin{eqnarray}
	Q_{21}=\left(
	\begin{array}{lll}
		{0} &{0}&{0}\\
		{0} &{I^8_{N_{1,\Gamma_t}\times{N_{1,\Gamma_t}}}} &{0}\\
		{0} &{0} &{0}\\
	\end{array}
	\right),
	Q_{22}=\left(
	\begin{array}{lll}
		{E_{N_{2,bp}\times{N_{2,bp}}}} &{0}&{0}\\
		{0} &{I^9_{N_{2,\Gamma_t}\times{N_{2,\Gamma_t}}}} &{0}\\
		{0} &{0} &{M^2_{N_{2,inp}\times N_{2,inp}}}\\
	\end{array}
	\right),
\end{eqnarray}
\begin{eqnarray}
	C_{21}=\left(
	\begin{array}{lll}
		{0} &{0}&{0}\\
		{0} &{I^{10}_{N_{1,\Gamma_t}\times{N_{1,\Gamma_t}}}} &{0}\\
		{0} &{0} &{P^2_{N_{1,inp}\times N_{1,inp}}}\\
	\end{array}
	\right),
\end{eqnarray}
\begin{eqnarray}
	G_{11}=\left(
	\begin{array}{lll}
		{J^1_{N_{1,bp}\times{N_{1,bp}}}} &{0}&{0}\\
		{0} &{0} &{0}\\
		{0} &{0} &{0}\\
	\end{array}
	\right)
	G_{21}=\left(
	\begin{array}{lll}
		{J^2_{N_{1,bp}\times{N_{1,bp}}}} &{0}&{0}\\
		{0} &{0} &{0}\\
		{0} &{0} &{0}\\
	\end{array}
	\right),
\end{eqnarray}
\begin{eqnarray}
	P=\left(
	\begin{array}{lll}
		{P^3_{N_{1,bp}\times N_{1,bp}}} &{0}&{0}\\
		{0} &{P^4_{N_{1,bp}\times N_{1,bp}}} &{0}\\
		{0} &{0} &{P_{N_{1,inp}\times N_{1,inp}}}\\
	\end{array}
	\right),
\end{eqnarray}
in which
\begin{eqnarray}
	M^1_{N_{1,inp}\times N_{1,inp}}(i,j,t)=\rho_1 E(3,:,t)-\beta_1E(4,:,t)-\beta_1E(5,:,t)+E(1,:,t),
\end{eqnarray}
\begin{eqnarray}
	M^2_{N_{2,inp}\times N_{2,inp}}(i,j,t)=\rho_2 E(3,:,t)-\beta_2E(4,:,t)-\beta_2E(5,:,t),
\end{eqnarray}
\begin{eqnarray}
	Q^1_{N_{1,inp}\times N_{1,inp}}(i,j,t)=\rho_1 E(3,:,t)-\beta_1E(4,:,t)-\beta_1E(5,:,t)+E(2,:,t),
\end{eqnarray}
\begin{eqnarray}
	p^{1}_{N_{1,inp}\times N_{1,inp}}(i,j,t)=E(1,:,t),
\end{eqnarray}
\begin{eqnarray}
	p^{2}_{N_{1,inp}\times N_{1,inp}}(i,j,t)=E(2,:,t),
\end{eqnarray}
\begin{eqnarray}
	p^{3}_{N_{1,bp}\times N_{1,bp}}(i,j,t)=I,
\end{eqnarray}
\begin{eqnarray}
	p^{4}_{N_{1,bp}\times N_{1,bp}}(i,j,t)=I,
\end{eqnarray}
\begin{eqnarray}
	p_{N_{1,inp}\times N_{1,inp}}(i,j,t)=E(4,:,t)+E(5,:,t),
\end{eqnarray}
\begin{eqnarray}
	J^1_{N_{1,bp}\times N_{1,bp}}(i,j,t)=E(1,:,t),
\end{eqnarray}
\begin{eqnarray}
	J^2_{N_{1,bp}\times N_{1,bp}}(i,j,t)=E(2,:,t),
\end{eqnarray}
\begin{eqnarray}
	I^1_{N_{1,\Gamma_t}\times N_{1,\Gamma_t}}(i,j,t)=\beta_1 E(1,:,t)n_{11},
\end{eqnarray}
\begin{eqnarray}
	I^2_{N_{2,\Gamma_t}\times N_{2,\Gamma_t}}(i,j,t)=\beta_1 E(1,:,t)n_{12},
\end{eqnarray}
\begin{eqnarray}
	I^3_{N_{1,\Gamma_t}\times N_{1,\Gamma_t}}(i,j,t)=\beta_2 E(1,:,t)n_{21},
\end{eqnarray}
\begin{eqnarray}
	I^4_{N_{2,\Gamma_t}\times N_{2,\Gamma_t}}(i,j,t)=\beta_2 E(1,:,t)n_{22},
\end{eqnarray}
\begin{eqnarray}
	I^5_{N_{1,\Gamma_t}\times N_{1,\Gamma_t}}(i,j,t)=n_{11},
\end{eqnarray}
\begin{eqnarray}
	I^6_{N_{1,\Gamma_t}\times N_{1,\Gamma_t}}(i,j,t)=\beta_1 E(2,:,t)n_{21},
\end{eqnarray}
\begin{eqnarray}
	I^7_{N_{2,\Gamma_t}\times N_{2,\Gamma_t}}(i,j,t)=\beta_1 E(2,:,t)n_{11},
\end{eqnarray}
\begin{eqnarray}
	I^8_{N_{1,\Gamma_t}\times N_{1,\Gamma_t}}(i,j,t)=\beta_2 E(2,:,t)n_{12},
\end{eqnarray}
\begin{eqnarray}
	I^9_{N_{2,\Gamma_t}\times N_{2,\Gamma_t}}(i,j,t)=\beta_2 E(2,:,t)n_{22},
\end{eqnarray}
\begin{eqnarray}
	I^{10}_{N_{1,\Gamma_t}\times N_{1,\Gamma_t}}(i,j,t)=n_{21},
\end{eqnarray}
and
\begin{eqnarray}
	U=\left(
	\begin{array}{ll}
		{U_1} \\
		\hline
		{U_2}
	\end{array}
	\right),
	V=\left(
	\begin{array}{ll}
		{V_1} \\
		\hline
		{V_2}
	\end{array}
	\right)
	P_1=[P_1],
\end{eqnarray}
\begin{eqnarray}
	D=\left(
	\begin{array}{ll}
		{D_1} \\
		\hline
		{D_2}
	\end{array}
	\right),
	H=\left(
	\begin{array}{ll}
		{H_1} \\
		\hline
		{H_2}
	\end{array}
	\right),
	F_p=[F_p],
\end{eqnarray}
here
\begin{eqnarray}
	U_1=(u_{1,1},\cdots, u_{1,N_{1,bp}},u_{1,N_{1,bp}+1},\cdots, u_{1,N_{1,bp}+N_{1,\Gamma_t}},u_{1,N^+_{bp}+N_{1,\Gamma_t}+1},\cdots, u_{1,N_{1,bp}+N_{1,\Gamma_t}+N_{1,inp}})^T,
\end{eqnarray}
\begin{eqnarray}
	U_2=(u_{2,1},\cdots, u_{2,N_{2,bp}},u_{2,N_{2,bp}+1},\cdots, u_{2,N_{2,bp}+N_{2,\Gamma_t}},u_{2,N^+_{bp}+N_{2,\Gamma_t}+1},\cdots, u_{2,N_{2,bp}+N_{2,\Gamma_t}+N_{2,inp}})^T,
\end{eqnarray}
\begin{eqnarray}
	V_1=(v_{1,1},\cdots, v_{1,N_{1,bp}},v_{1,N_{1,bp}+1},\cdots, v_{1,N_{1,bp}+N_{1,\Gamma_t}},v_{1,N^+_{bp}+N_{1,\Gamma_t}+1},\cdots, v_{1,N_{1,bp}+N_{1,\Gamma_t}+N_{1,inp}})^T,
\end{eqnarray}
\begin{eqnarray}
	V_2=(v_{2,1},\cdots, v_{2,N_{2,bp}},v_{2,N_{2,bp}+1},\cdots, v_{2,N_{2,bp}+N_{2,\Gamma_t}},v_{2,N^+_{bp}+N_{2,\Gamma_t}+1},\cdots, v_{2,N_{2,bp}+N_{2,\Gamma_t}+N_{2,inp}})^T,
\end{eqnarray}
\begin{eqnarray}
	P_1=(p_{1,1},\cdots, p_{1,N_{1,bp}},p_{1,N_{1,bp}+1},\cdots, p_{1,N_{1,bp}+N_{1,\Gamma_t}},p_{1,N^+_{bp}+N_{1,\Gamma_t}+1},\cdots, p_{1,N_{1,bp}+N_{1,\Gamma_t}+N_{1,inp}})^T,
\end{eqnarray}
\begin{eqnarray}
	D_1=(0_{1},\cdots, 0_{N_{1,bp}},0_{N_{1,bp}+1},\cdots, 0_{N_{1,bp}+N_{1,\Gamma_t}},{f}_{11,N_{1,bp}+N_{1,\Gamma_t}+1},\cdots, {f}_{11,N_{1,bp}+N_{1,\Gamma_t}+N_{1,inp}})^T,
\end{eqnarray}
\begin{eqnarray}
	D_2=(0_{1},\cdots, 0_{N_{1,bp}},\tau_{1,N_{2,bp}+1},\cdots, \tau_{1,N_{2,bp}+N_{2,\Gamma_t}},{f}_{12,N_{2,bp}+N_{2,\Gamma_t}+1},\cdots, {f}_{12,N_{2,bp}+N_{2,\Gamma_t}+N_{2,inp}})^T,
\end{eqnarray}
\begin{eqnarray}
	H_1=(0_{1},\cdots, 0_{N_{1,bp}},0_{N_{1,bp}+1},\cdots, 0_{N_{1,bp}+N_{1,\Gamma_t}},{f}_{21,N_{1,bp}+N_{1,\Gamma_t}+1},\cdots, {f}_{21,N_{1,bp}+N_{1,\Gamma_t}+N_{1,inp}})^T,
\end{eqnarray}
\begin{eqnarray}
	H_2=(0_{1},\cdots, 0_{N_{1,bp}},\tau_{2,N_{2,bp}+1},\cdots, \tau_{2,N_{2,bp}+N_{2,\Gamma_t}},{f}_{22,N_{2,bp}+N_{2,\Gamma_t}+1},\cdots, {f}_{22,N_{2,bp}+N_{2,\Gamma_t}+N_{2,inp}})^T,
\end{eqnarray}
\begin{eqnarray}
	F_p=(p_{1,1},\cdots, p_{1,N_{1,bp}},p_{1,N_{1,bp}+1},\cdots, p_{1,N_{1,bp}+N_{1,\Gamma_t}},{f}_{3,N_{1,bp}+N_{1,\Gamma_t}+1},\cdots, {f}_{3,N_{1,bp}+N_{1,\Gamma_t}+N_{1,inp}})^T,
\end{eqnarray}
in which 
\begin{eqnarray}
	f_3=\frac{\partial f_{11}}{\partial x}+\frac{\partial f_{12}}{\partial y}.
\end{eqnarray}
Here I is the unit matrix. The matrices $M$ and $Q$ are formed according to the similarly idea of domain decomposition, $M_{11}$ and $Q_{11},M_{22}$ and $Q_{22}$are formed by the governing equations (Eq.(1), Eq.(3)), the Dirichlet boundary conditions (Eq.(6), Eq.(7)) and the interface condition (Eq.(4)) contain the information in $\Omega_t^1$ and $\Omega_t^2$, respectively. The interface parts of $M_{21}, M_{22}, Q_{21}, Q_{22}, Y, L$ are all created by the interface condition (Eq.(5)). Obviously, it is simple to use the above matrices to describe the interface conditions. $C_1$ and $C_2$ obtain the pressure information in Eq.(1). $P$ is formed according to the governing equation (Eq.(14)). $G_1$ and $G_2$ are related with the boundary condition (Eq.(17)), we can obtain different matrix if we use different form which is used to enrich the information about the pressure $p_1$.We can see that the matrix K is sparse.
\section{Numerical examples}
In this section, five examples are provided to show the simplicity, accuracy and the stability of ST-GFDM for the Stokes/Parabolic moving interface problem with a linear moving circle interface in Example 1, an irregular moving circle interface in Example 2, an irregular moving octagon interface in Example 3, a deformation of a three-petalled flower interface in Example 4 and an eight-petalled flower interface in Example 5.

For simplicity, we defined the error norms as follows:
\begin{eqnarray}
	L_\infty=max|u_i-u(x_i)|(i=1,2,\cdots, N_{T}),
\end{eqnarray}

\begin{eqnarray}
	L_2=[\sum_{i=1}^{N_{T}}{\frac{(u_i-u(x_i))^2}{N_{T}}}]^{\frac{1}{2}}(i=1,2,\cdots, N_{T}),
\end{eqnarray}
\begin{eqnarray}
	H^1=[\sum_{i=1}^{N_{T}}{\frac{(\nabla u_i-\nabla u(x_i))^2}{N_{T}}}]^{\frac{1}{2}}(i=1,2,\cdots, N_{T}).
\end{eqnarray}
and the relative errors are defined as follows
\begin{eqnarray}
	L_{\infty,relative}=\frac{max|u_i-u(x_i)|}{max|u(x_i)|}(i=1,2,\cdots, N_{T}),
\end{eqnarray}
\begin{eqnarray}
	L_{2,relative}=\frac{[\sum_{i=1}^{N_{T}}{\frac{(u_i-u(x_i))^2}{N_{T}}}]^{\frac{1}{2}}}{[\sum_{i=1}^{N_{T}}{\frac{u(x_i)^2}{N_{T}}}]^{\frac{1}{2}}}(i=1,2,\cdots, N_{T}),
\end{eqnarray}
\begin{eqnarray}
	H^1_{relative}=\frac{[\sum_{i=1}^{N_{T}}{\frac{(\nabla u_i-\nabla u(x_i))^2}{N_{T}}}]^{\frac{1}{2}}}{[\sum_{i=1}^{N_{T}}{\frac{(\nabla u(x_i))^2}{N_{T}}}]^{\frac{1}{2}}}(i=1,2,\cdots, N_{T}).
\end{eqnarray}
$u_i$ and $u(x_i)$ are the numerical and exact solution at point $x_i$, respectively. Let the domain $\Omega$=[0,1] $\times$ [0,1]$ \times$ [0,1], $\Omega_1$ is outside the interface $\Gamma_t$ and $\Omega^2=\Omega/(\Omega^1\cup\Gamma_t)$. $N_{T}$ is the number of all scattered nodes in $\Omega^1_t$, $\Omega^2_t$, $\partial\Omega_t$ and $\Gamma_t$. Namely, $N_{T}=N_{1,inp}+N_{2,inp}+N_{bp}+N_{1,\Gamma}+N_{2,\Gamma}$.  $u$ and $v$ are the x-component and the y-component of $\mathbf{u}=(u,v)$, respectively. We adopt $m=60$for the ST-GFDM and we assume $dt=0.1,t=0.5$ for all figures and tables unless they are stated otherwise.

\subsection{Example 1: The Stokes/Parabolic interface problem with a linear moving circle interface.}
In this example, we consider a Stokes/Parabolic interface problem with a linear moving circle interface (From Ref.[11]) (see Fig.6(left)) $\Gamma_t:\varphi_1=(x-0.3-0.1t)^2+(y-0.3-0.1t)^2-0.01=0.$
The coefficient $\beta_1=100,\beta_2=1,\rho_1=100,\rho_2=1.$
The exact solution is
\begin{eqnarray}
	u&=&(u_1,u_2)=(y-0.3-0.1t)((x-0.3-0.1t)^2+(y-0.3-0.1t)^2-0.01)t/\beta,\\
	v&=&(v_1,v_2)=-(x-0.3-0.1t)((x-0.3-0.1t)^2+(y-0.3-0.1t)^2-0.01)t/\beta, \\
	p&=&0.1(x^3-y^3)((x-0.3-0.1t)^2+(y-0.3-0.1t)^2-0.01).
\end{eqnarray}
\begin{figure}
	\centering
	\includegraphics[scale=.4]{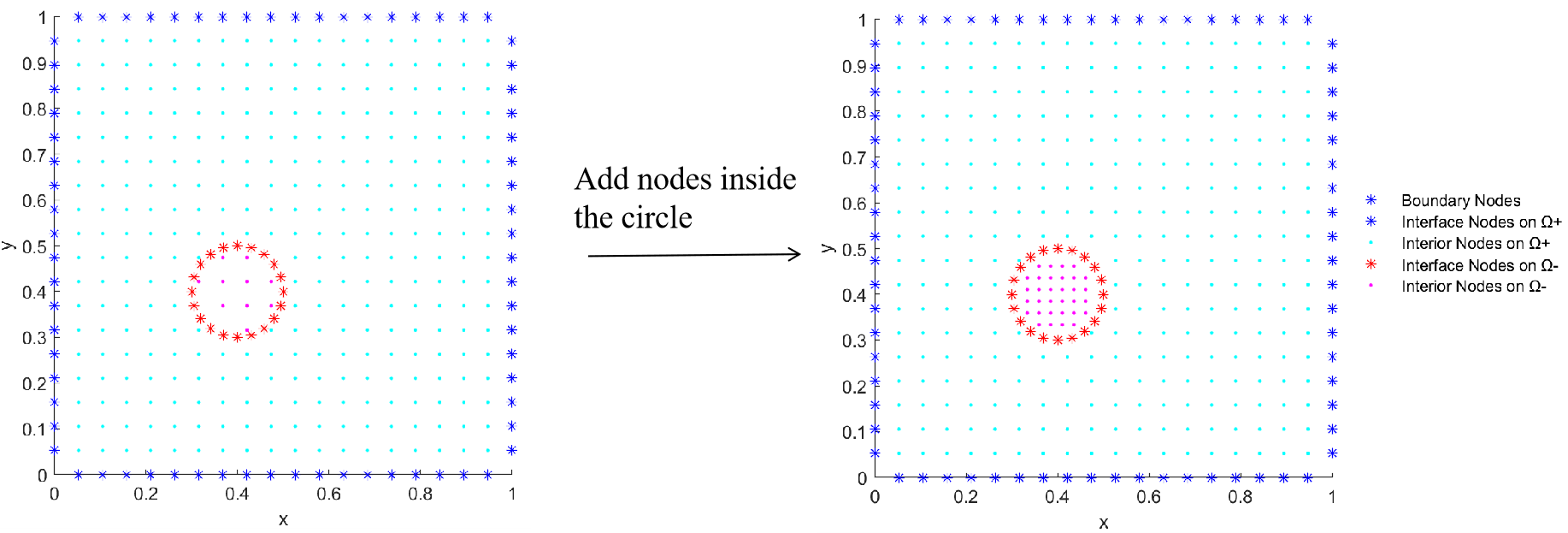}     
	\caption{ The method of locally dense nodes for Example 1.}
\end{figure}
\begin{figure}
	\centering
	\includegraphics[scale=.4]{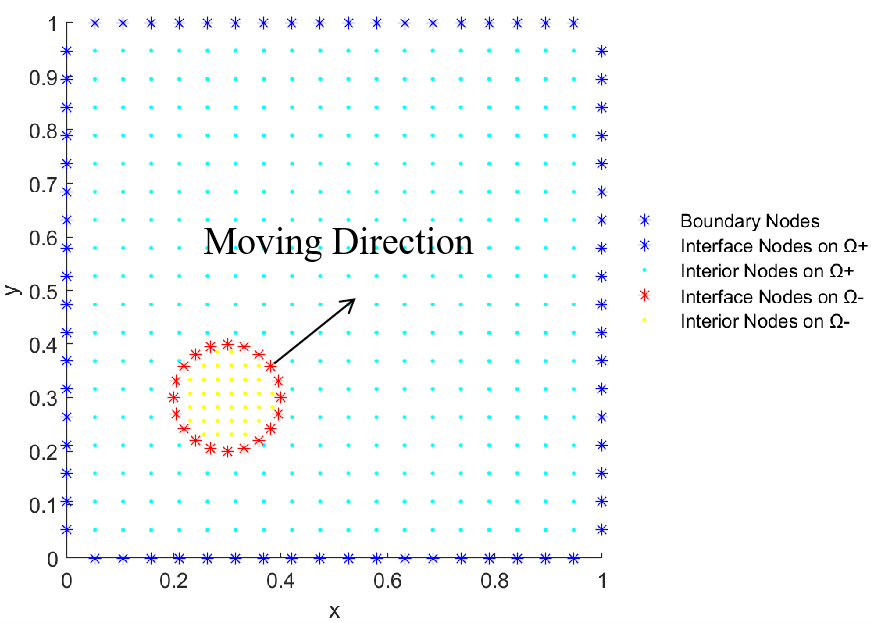}    
	\includegraphics[scale=.4]{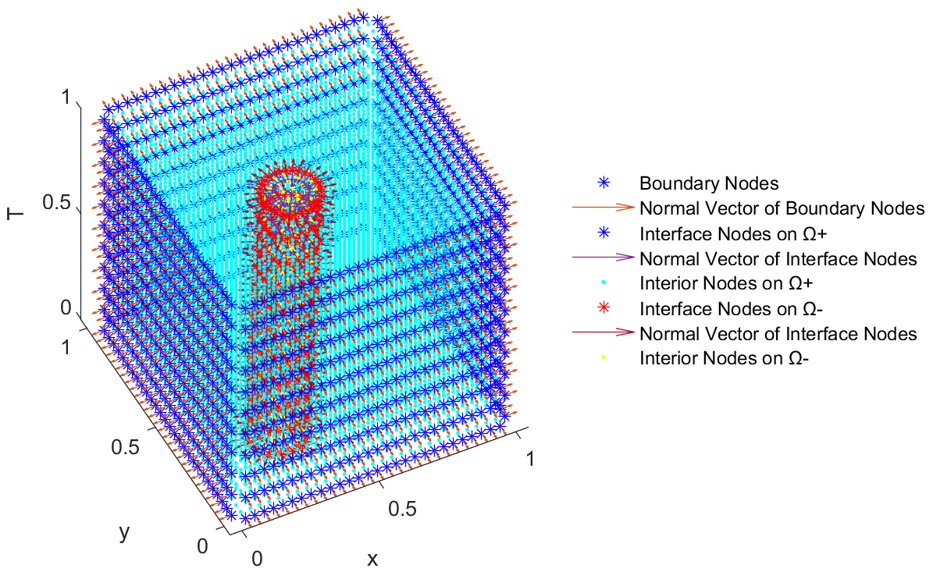}    
	\caption{ The point collocation at t=0 (left) and the point collocation for all time (right) for Example 1.}
\end{figure}
% Figure

\begin{figure}
	\centering
	\includegraphics[scale=.4]{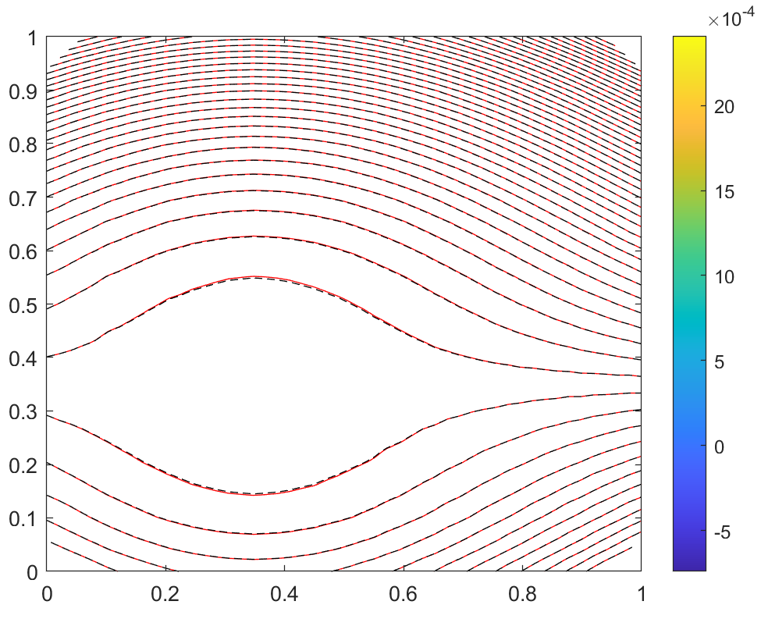}    
	\includegraphics[scale=.4]{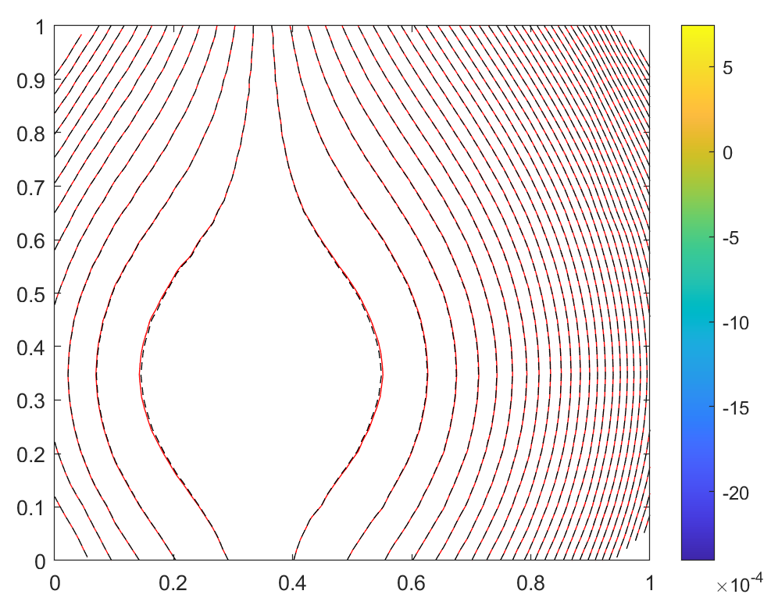}    
	\includegraphics[scale=.4]{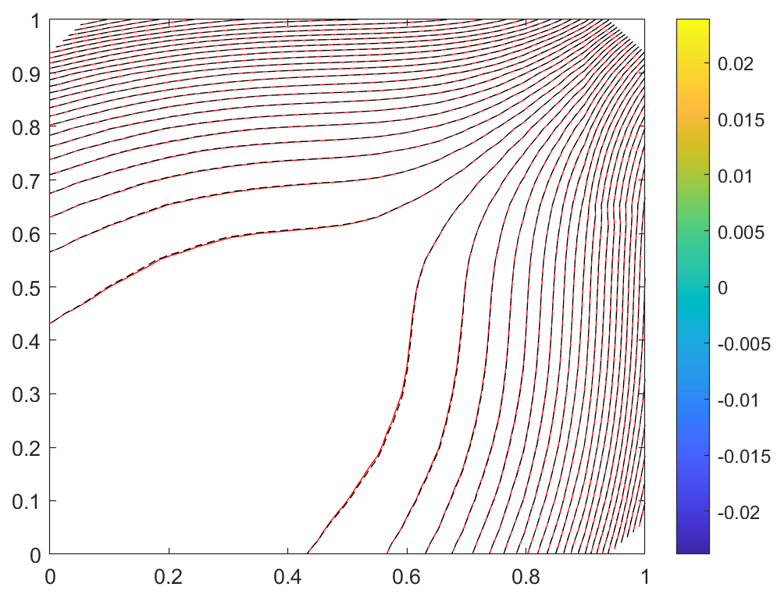} 
	\includegraphics[scale=.4]{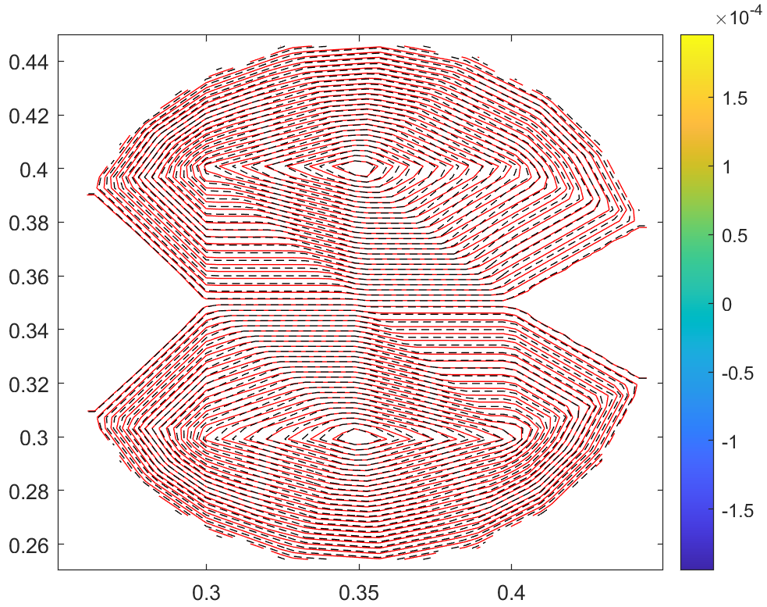}    
	\includegraphics[scale=.4]{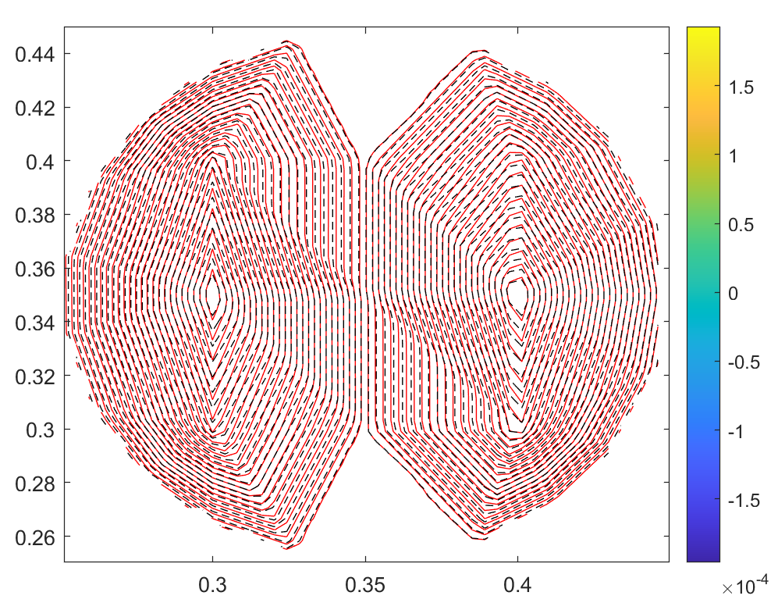}    
	\caption{The contour of numerical solution (black line) and exact solution (red lines) when $N_{T_1}=5213$ for $u_1$(left up),$v_1$(left middle),$p$(right up) and $u_2$(left down), $v_2$(right down) for Example 1.}
\end{figure}
\begin{figure}
	\centering
	\includegraphics[scale=.4]{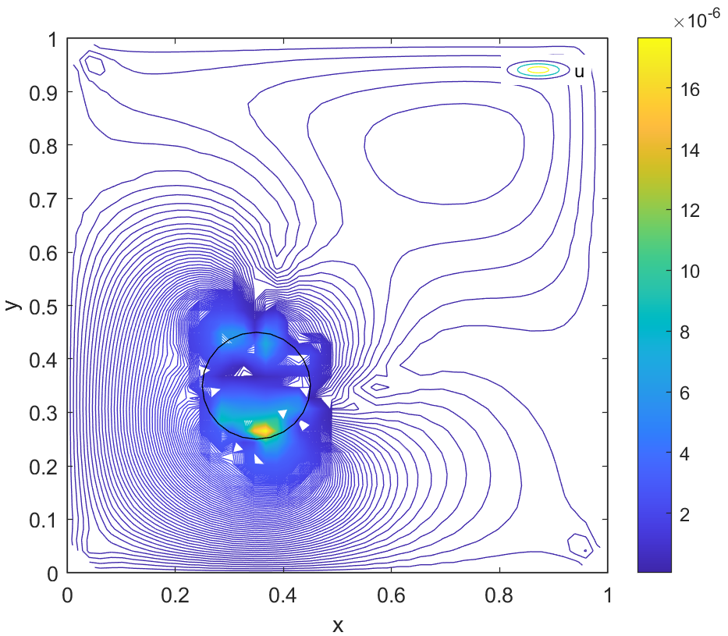}
	\includegraphics[scale=.4]{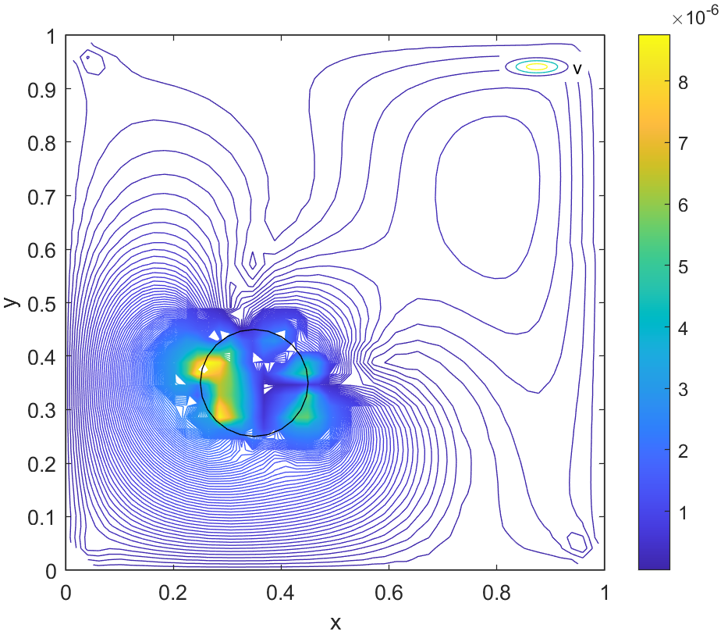}       
	\includegraphics[scale=.4]{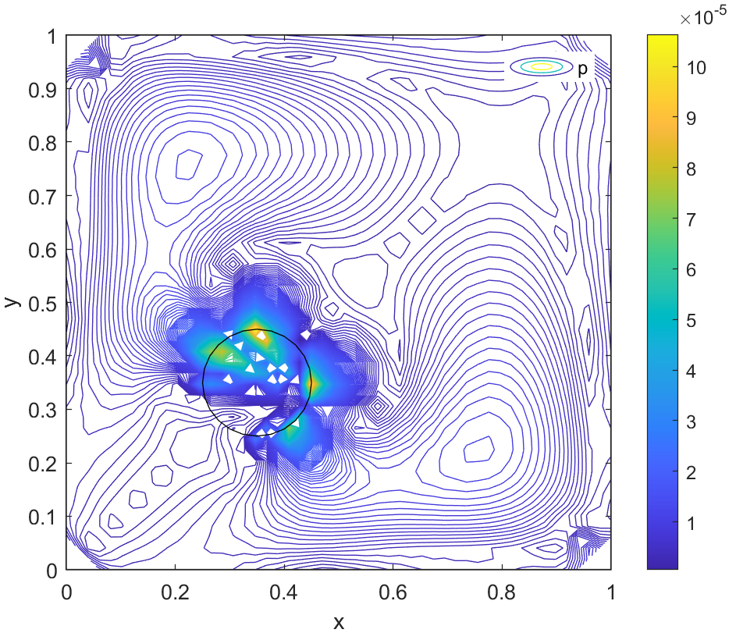}              
	\caption{ The contour of the absolute error when $N_{T_1}=13347$ for Example 1.}
\end{figure}
% Figure
\begin{center}
	\begin{table*}	
		\scriptsize
		\caption{ $L_{\infty}$,$L_2$ and $H^1$ errors of ST-GFDM for Example 1}
		\begin{tabular}{ccccccccccc}
			\hline
			\multirow{1}{*}{$N_{T}$} & \multicolumn{3}{c}{u} & \multicolumn{3}{c}{v}& \multicolumn{3}{c}{p} & \multicolumn{1}{c}{$Time(s)$}\\  
			\hline		
			& $L_{\infty}$      &  $L_2 $   &   $H^1 $
			
			& $L_{\infty}$      &  $L_2$   &   $H^1$
			
			& $L_{\infty}$      &  $L_2$   &   $H^1$ & $ $\\
			\hline

			$10230$ & $1.98\times10^{-5}$ & $8.15\times10^{-6}$ & $6.29\times10^{-4}$ &  $1.13\times10^{-5}$ & $4.40\times10^{-6}$ & $6.79\times10^{-4}$&
			$1.32\times10^{-4}$ & $ 1.87\times10^{-5}$ & $8.30\times10^{-4}$& $19.9$\\
			
			$13347$ & $ 1.86\times10^{-5}$ & $5.73\times10^{-6}$ & $ 4.82\times10^{-4}$ &  $1.15\times10^{-5}$ & $4.43\times10^{-6}$ & $4.75\times10^{-4}$&
			$1.14\times10^{-4}$ & $ 1.29\times10^{-5}$ & $5.58\times10^{-4}$ &$34.4$\\
			\hline
		\end{tabular}
	\end{table*}
\end{center} 
\begin{center}
	\begin{table*}	
		\scriptsize
		\caption{ $L_{\infty, relative}$,$L_{2,relative}$ and $H^1_{relative}$ errors of ST-GFDM for Example 1}
		\begin{tabular}{ccccccccccc}
			\hline
			\multirow{1}{*}{$N_{T}$} & \multicolumn{3}{c}{u} & \multicolumn{3}{c}{v}& \multicolumn{3}{c}{p} & \multicolumn{1}{c}{$Time(s)$}\\  
			\hline		
			& $L_{\infty, relative}$  &$L_{2,relative}$   &   $H^1_{relative}$
			
			& $L_{\infty, relative}$  &$L_{2,relative}$   &   $H^1_{relative}$
			
			& $L_{\infty, relative}$  &$L_{2,relative}$   &   $H^1_{relative}$ & $ $\\
			\hline

			$10230$ & $1.03\times10^{-1}$ & $9.69\times10^{-2}$ & $1.08\times10^{-1}$ &  $7.02\times10^{-2}$ & $7.18\times10^{-2}$ & $1.14\times10^{-1}$&
			$5.24\times10^{-3}$ & $2.36\times10^{-3}$ & $1.42\times10^{-2}$& $21.2$\\
			
			$13347$ & $9.81\times10^{-2}$ & $6.34\times10^{-2}$ & $8.75\times10^{-2}$ &  $6.04\times10^{-2}$ & $6.15\times10^{-2}$ & $8.63\times10^{-2}$&
			$4.49\times10^{-3}$ & $1.66\times10^{-3}$ & $9.67\times10^{-3}$ &$34.7$\\
			\hline
		\end{tabular}
	\end{table*}
\end{center} 
\begin{center}
	\begin{table*}	
		\scriptsize
		\caption{ The comparison between the ST-GFDM and the DLM/FD FEM[11] when $\frac{\beta_2}{\beta_1}=100, \frac{\rho_2}{\rho_1}=100, t=1$ for Example 1}
		\begin{tabular}{ccccccccccc}
			\hline
			\multirow{1}{*}{$ $}&\multicolumn{1}{c}{$N_x$} & \multicolumn{2}{c}{$u_1$} & \multicolumn{2}{c}{$v_1$}& \multicolumn{1}{c}{$p$} & \multicolumn{2}{c}{$u_2$} & \multicolumn{2}{c}{$v_2$}\\
			\hline
			
			$ $&$ $& $H^1$      &  $L_2$  
			& $H^1$      &  $L_2$  
			&  $L_2$  
			& $H^1$      &  $L_2$  
			
			& $H^1$      &  $L_2$ \\ 
			\hline		
			DLM/FD& $16$ & $5.56\times10^{-4}$ & $3.00\times10^{-5}$ & $5.11\times10^{-4}$ &  $2.96\times10^{-5}$ & $1.20\times10^{-3}$ & $1.70\times10^{-6}$&
			$2.03\times10^{-6}$ & $2.39\times10^{-6}$ & $2.32\times10^{-6}$\\
			
			FEM[11]&$32$& $3.43\times10^{-4}$ & $1.65\times10^{-5}$ & $3.50\times10^{-4}$ &  $1.66\times10^{-5}$ & $6.38\times10^{-4}$ & $1.37\times10^{-6}$&
			$1.08\times10^{-6}$ & $ 1.29\times10^{-6}$ & $1.12\times10^{-6}$\\
			
			$ $& $Order$ & $0.70$ & $0.86$ & $0.55$ &  $0.83$ & $0.91$ & $0.31$&
			$0.91$ & $0.89$ & $1.05$\\
			\hline		
			ST GFDM&$16$ & $6.36\times10^{-5}$ & $ 5.68\times10^{-6}$ & $ 6.20\times10^{-5}$ &  $5.66\times10^{-6}$ & $5.70\times10^{-5}$ & $ 1.77\times10^{-6}$&
			$1.27\times10^{-5}$ & $ 1.85\times10^{-6}$ & $1.27\times10^{-5}$\\
			
			$ $ &$32$  & $ 2.00\times10^{-5}$ & $5.26\times10^{-7}$ & $1.99\times10^{-5}$ &  $5.16\times10^{-7}$ & $7.89\times10^{-6}$ & $  4.67\times10^{-7}$&	$ 1.93\times10^{-7}$ & $4.97\times10^{-7}$ & $9.18\times10^{-8}$\\
			
			$ $& $Order$ & $1.67$ & $3.43$ & $ 1.64$ &  $3.46$ & $1.92$ & $5.80$&
			$2.72$ & $1.90$ & $7.11$\\
			\hline
			ST GFDM&$16$ & $ 6.02\times10^{-5}$ & $  4.87\times10^{-6}$ & $5.58\times10^{-5}$ &  $4.79\times10^{-6}$ & $4.00\times10^{-5}$ & $ 1.03\times10^{-6}$&
			$1.10\times10^{-5}$ & $1.02\times10^{-6}$ & $1.10\times10^{-5}$\\
			
			Locally  &$32$  & $7.01\times10^{-7}$ & $5.17\times10^{-7}$ & $ 1.99\times10^{-5}$ &  $ 5.12\times10^{-7}$ & $7.94\times10^{-6}$ & $  4.16\times10^{-8}$&	$2.90\times10^{-8}$ & $3.55\times10^{-7}$ & $2.12\times10^{-8}$\\
			
			dense nodes& $Order$ & $6.42$ & $3.23$ & $1.49$ &  $3.23$ & $2.33$ & $4.63$&
			$8.57$ & $1.52$ & $5.70$\\
			\hline			
		\end{tabular}
	\end{table*}
\end{center} 

\begin{table*}	
	\scriptsize
	\caption{ $L_{\infty}$,$L_2$ and $H^1$ errors of ST-GFDM with different $m$  when $N_T=13347$ for Example 1}
	\begin{tabular}{cccccccccc}
		\hline
		\multirow{1}{*}{m} & \multicolumn{3}{c}{u} & \multicolumn{3}{c}{v}& \multicolumn{3}{c}{p} \\
		\hline
		
		& $L_{\infty}$      &  $L_2 $   &   $H^1$
		
		& $L_{\infty}$      &  $L_2$   &   $H^1$
		
		& $L_{\infty}$      &  $L_2$   &   $H^1$\\
		\hline
		$55$ & $2.41\times10^{-4}$ & $8.28\times10^{-5}$ & $3.96\times10^{-3}$ &  $2.41\times10^{-4}$ & $8.84\times10^{-5}$ & $4.01\times10^{-3}$&
		$3.34\times10^{-3}$ & $3.47\times10^{-4}$ & $4.72\times10^{-3}$ \\

		$57$ & $ 6.86\times10^{-4}$ & $1.68\times10^{-4}$ & $ 7.31\times10^{-3}$ &  $  6.85\times10^{-4}$ & $1.77\times10^{-4}$ & $  7.36\times10^{-3}$&
		$7.26\times10^{-4}$ & $6.53\times10^{-6}$ & $8.71\times10^{-3}$ \\
		
		$59$ & $5.33\times10^{-4}$ & $1.36\times10^{-4}$ & $ 5.57\times10^{-3}$ &  $  9.39\times10^{-4}$ & $1.87\times10^{-4}$ & $7.23\times10^{-3}$&
		$6.39\times10^{-3}$ & $5.76\times10^{-4}$ & $9.07\times10^{-3}$ \\
		
		$61$ & $1.83\times10^{-5}$ & $6.58\times10^{-6}$ & $ 5.00\times10^{-4}$ &  $  1.33\times10^{-5}$ & $5.44\times10^{-6}$ & $4.91\times10^{-4}$&
		$1.96\times10^{-4}$ & $1.92\times10^{-5}$ & $6.04\times10^{-4}$ \\
		
		$63$ & $2.02\times10^{-5}$ & $6.31\times10^{-6}$ & $  4.83\times10^{-4}$ &  $ 2.13\times10^{-5}$ & $6.13\times10^{-6}$ & $ 4.90\times10^{-4}$&
		$9.70\times10^{-5}$ & $1.22\times10^{-5}$ & $5.54\times10^{-4}$ \\
		
		\hline
	\end{tabular}
\end{table*}
% Figure
\begin{figure}
	\centering
	\includegraphics[scale=.4]{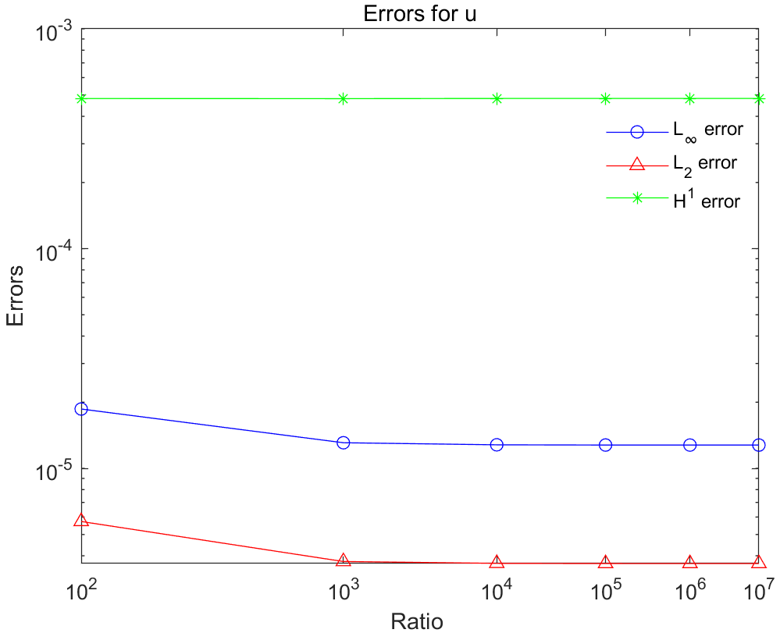}    
	\includegraphics[scale=.4]{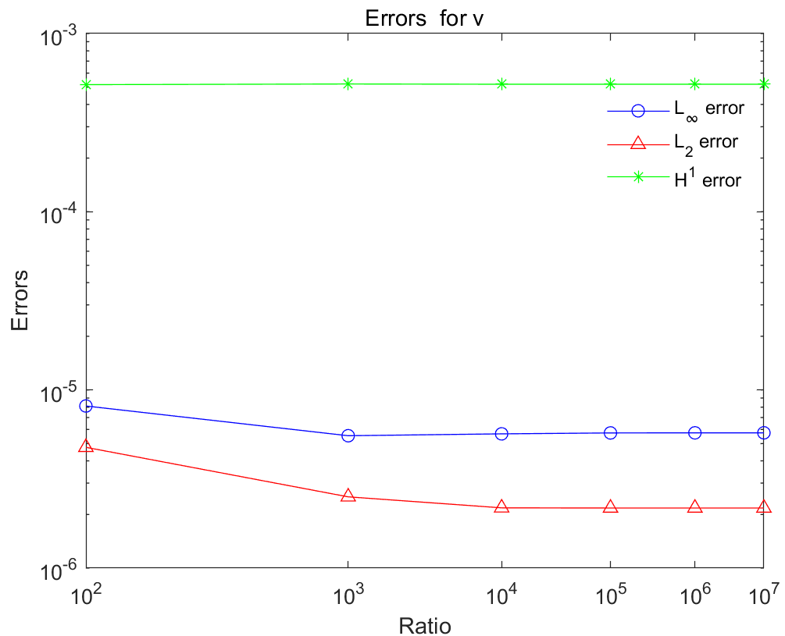}    
	\includegraphics[scale=.4]{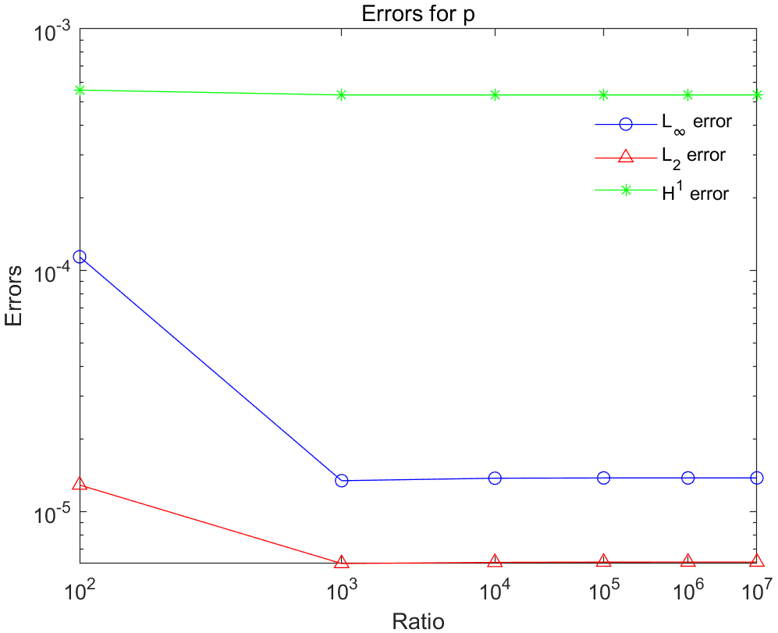}    
	\caption{ $L_{\infty}$,$L_2$ and $H^1$ errors of ST-GFDM with different jump $Ratio=\frac{\beta_1}{\beta_2}=\frac{\rho_1}{\rho_2}=(10^2,10^3,10^4,10^5,10^6,10^7)$ when $N_{T}=13347$ for Example 1.}
\end{figure}
% Figure
\begin{figure}
	\centering
	\includegraphics[scale=.4]{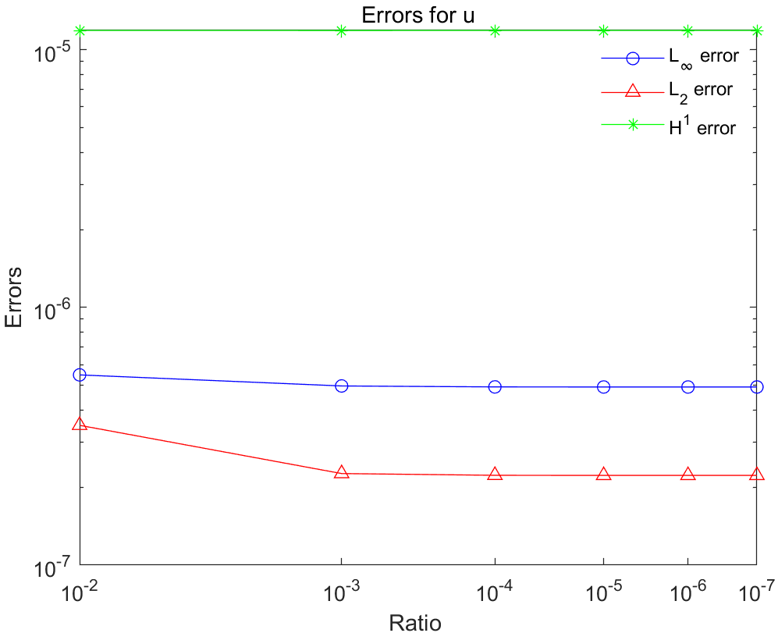}    
	\includegraphics[scale=.4]{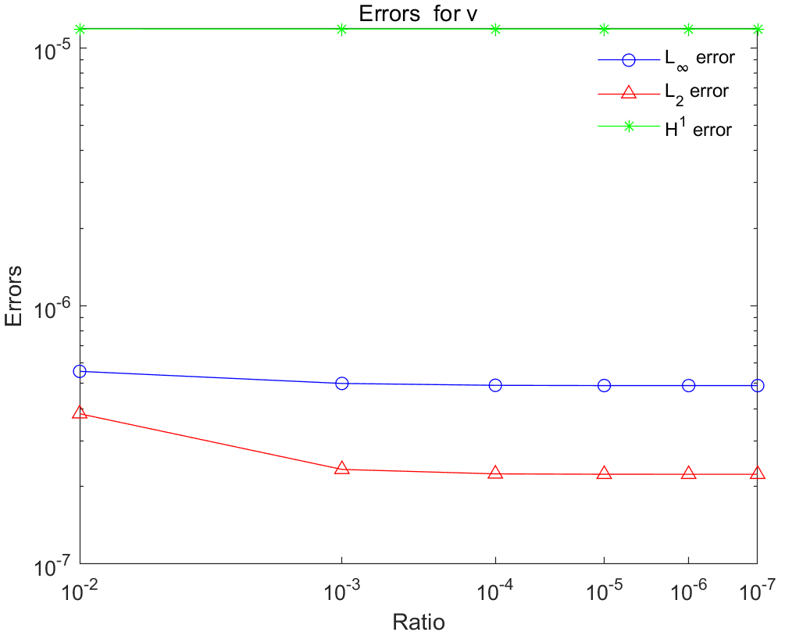}    
	\includegraphics[scale=.4]{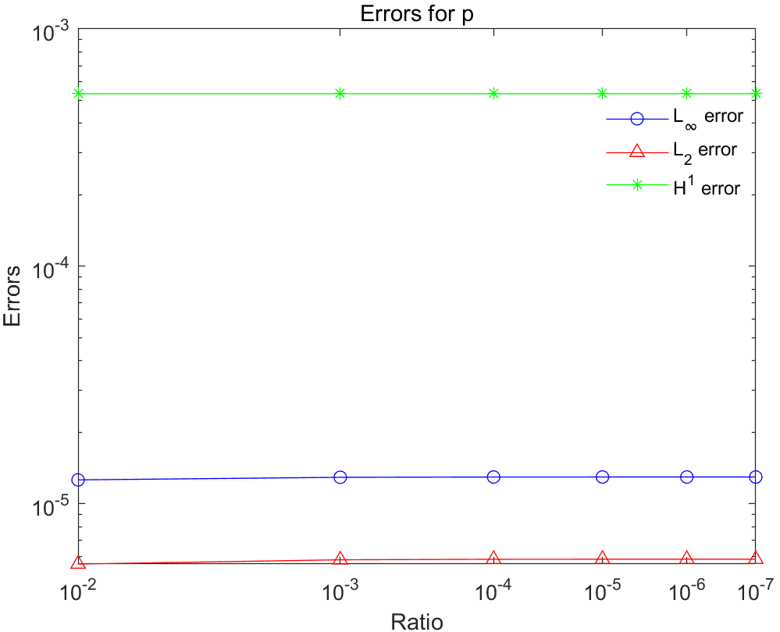}    
	\caption{ $L_{\infty}$,$L_2$ and $H^1$ errors of ST-GFDM with different jump $Ratio=\frac{\beta_1}{\beta_2}=\frac{\rho_1}{\rho_2}=(10^{-2},10^{-3},10^{-4},10^{-5},10^{-6},10^{-7})$ when $N_{T}=13347$ for Example 1.}
\end{figure}
\begin{figure}
	\centering
	\includegraphics[scale=.4]{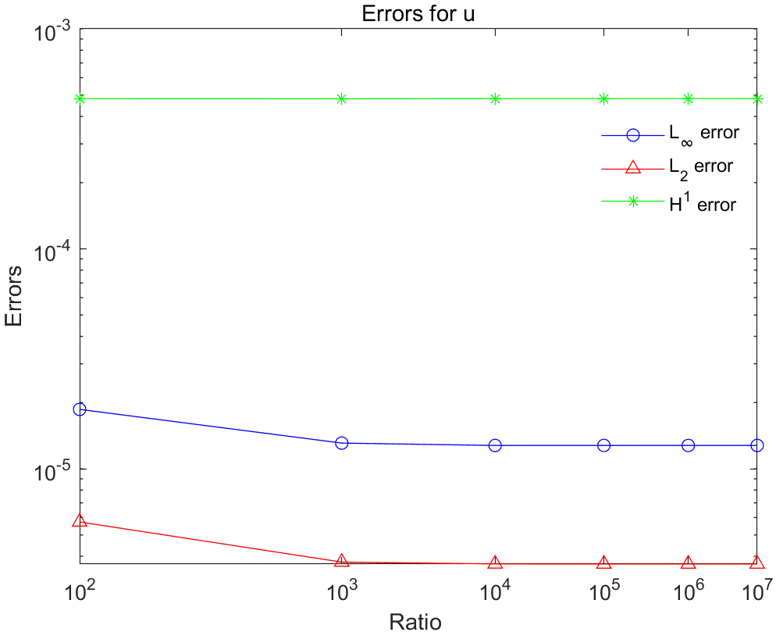}    
	\includegraphics[scale=.4]{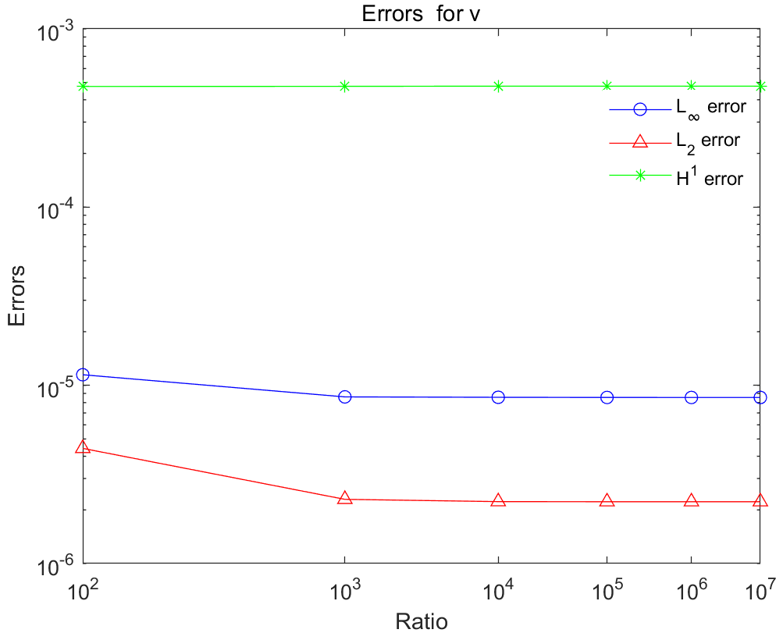}    
	\includegraphics[scale=.4]{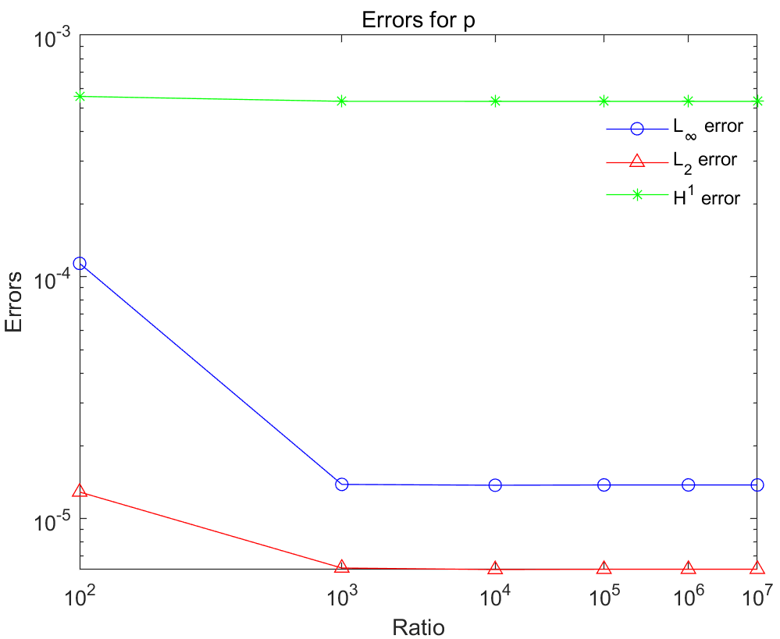}    
	\caption{$L_{\infty}$,$L_2$ and $H^1$ errors of ST-GFDM with different jump $Ratio=\frac{\beta_1}{\beta_2}=(10^2,10^3,10^4,10^5,10^6,10^7)$ when $\frac{\rho_1}{\rho_2}=100$, $N_{T}=13347$ for Example 1.}
\end{figure}

Fig.5 presents the method of locally dense nodes to show the point collocation in the subdomain $\Omega_t^2$, which is a small circle like a hole. Fig.6 presents the point collocation at $t=0$ and the point collocation for all times, we can see that the initial interface shape and moving direction of this example. Fig.7 pressents the contour of the numerical solution and the exact solution, we can see that the numerical solution is coincide with the exact solution. In Table. 1, $L_{\infty}$,$L_2$ and $H^1$ errors of ST-GFDM and computational time are provided to show the accuracy and high efficiency of the ST-GFDM. The $L_{\infty, relative}$,$L_{2,relative}$ and $H^1_{relative}$ errors of ST-GFDM are provided in Table.2. We can see that the relative error at acceptable range and our results are accurate. The comparison between the ST-GFDM and the DLM/FD FEM[11] when $\frac{\beta_2}{\beta_1}=100, \frac{\rho_2}{\rho_1}=100$ is provided to show the accuracy and stability of the ST-GFDM than the DLM/FD FEM[11] in Table 3. From Table 3, we can see that the numerical results of the ST-GFDM are more accurate and stable and can maintain a high order convergence. The results show that the ST-GFDM has an advantage in all errors. Furthermore, the more accurate results can be obtained by using the ST GFDM coupled with the locally dense nodes. The contour of the absolute error is shown in Fig.8. We can see that the large errors appear on the interface. In order to test the stability of ST-GFDM for Stokes/Parabolic moving interface problems, the number 'm' and the jump ratio are changed. From Table 4, we can see that the numerical errors are still stable when the number of 'm' is changed and the more accurate results can be obtained in the range of $m=59$ to $m=63$. According to our experience, the best results can be obtained when $m=60$.  In Fig.9, Fig.10 and Fig.11, we can see that the $L_{\infty}$, $L_2$ and $H^1$ errors are stable when the jump $Ratio=\frac{\beta_1}{\beta_2}=\frac{\rho_1}{\rho_2}$ and even in $Ratio=\frac{\beta_1}{\beta_2},\frac{\rho_1}{\rho_2}=100$. Therefore, the accuracy, stability and high efficiency of the ST-GFDM for the Stokes/Parabolic interface problem with linear moving circle interface are verified.

\subsection{Example 2: The Stokes/Parabolic interface problem with an irregular moving circle interface.}
In this example, we consider a Stokes/Parabolic interface problem with an irregular moving circle interface (From Ref.[11])(see Fig.12(left)) $\Gamma_t:\varphi_2=(x-0.3-0.1(t+sin(5t)))^2+(y-0.3-0.1(t+t^3))^2-0.01=0.$
The coefficient $\beta_1=100,\beta_2=1,\rho_1=100,\rho_2=1.$
The exact solution are
\begin{eqnarray}
	u&=&(y-0.3-0.1(t+t^3))((x-0.3-0.1(t+sin(5t)))^2+(y-0.3-0.1(t+t^3))^2-0.01)t/\beta,\\
	v&=&-(x-0.3-0.1(t+sin(5t)))((x-0.3-0.1(t+sin(5t)))^2+(y-0.3-0.1(t+t^3))^2-0.01)t/\beta, \\
	p&=&0.1(x^3-y^3)(x-0.3-0.1(t+sin(5t)))^2+(y-0.3-0.1(t+t^3))^2-0.01)t.
\end{eqnarray}
\begin{figure}
	\centering
	\includegraphics[scale=.4]{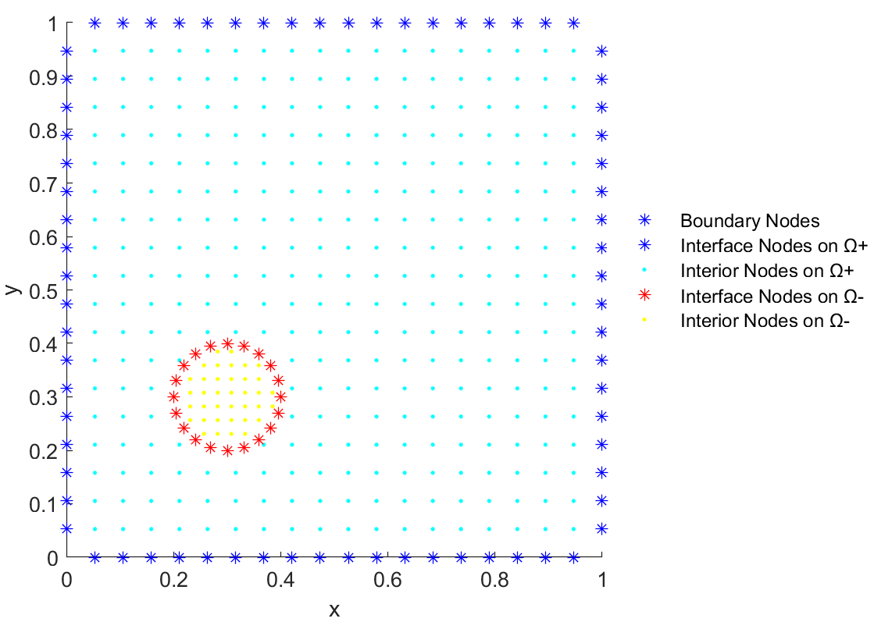}    
	\includegraphics[scale=.4]{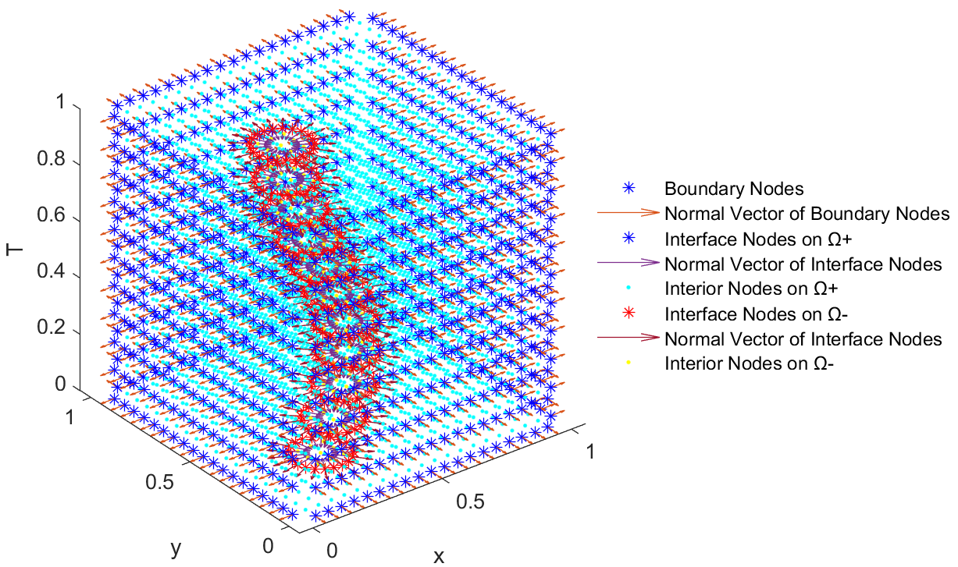}    
	\caption{ The point collocation at t=0 (left) and the point collocation for all time (right) for Example 2.}
\end{figure}
% Figure
\begin{figure}
	\centering
	\includegraphics[scale=.4]{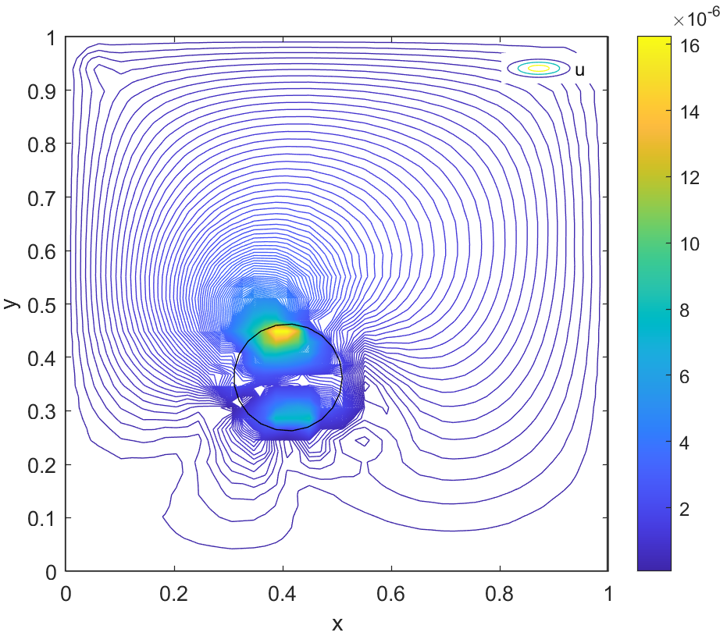}
	\includegraphics[scale=.4]{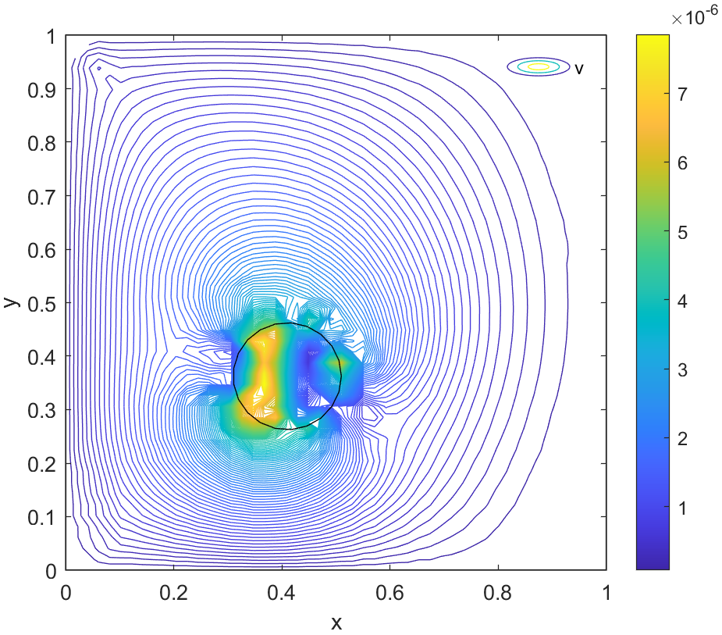}       
	\includegraphics[scale=.4]{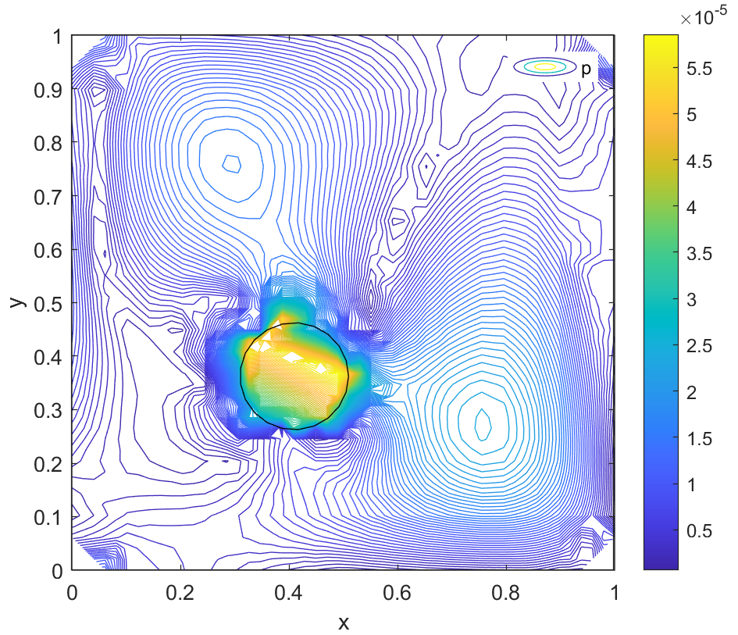}              
	\caption{ The contour of the absolute error when $N_{T}= 10460$ for Example 2.}
\end{figure}
\begin{center}
	\begin{table*}	
		\scriptsize
		\caption{ $L_{\infty}$,$L_{2}$ and $H^1$ errors of ST-GFDM when $N_T=10460$ for Example 2}
		\begin{tabular}{cccccccccc}
			\hline
			\multicolumn{3}{c}{u} & \multicolumn{3}{c}{v}& \multicolumn{3}{c}{p} \\  
			\hline		
			$L_{\infty}$  &$L_{2}$   &   $H^1$
			
			&  $L_{\infty}$ &$L_{2}$   &   $H^1$
			
			&  $L_{\infty}$  &$L_{2}$   &   $H^1$ \\
			\hline
			$1.89\times10^{-5}$ & $6.27\times10^{-6}$ & $4.62\times10^{-4}$ &  $8.12\times10^{-6}$ & $4.76\times10^{-6}$ & $5.17\times10^{-4}$&
			$6.36\times10^{-5}$ & $1.26\times10^{-5}$ &$7.12\times10^{-4}$\\
			\hline		
			$L_{\infty, relative}$  &$L_{2,relative}$   &   $H^1_{relative}$
			
			& $L_{\infty, relative}$  &$L_{2,relative}$   &   $H^1_{relative}$
			
			& $L_{\infty, relative}$  &$L_{2,relative}$   &   $H^1_{relative}$ \\
			\hline
			$1.00\times10^{-1}$ & $6.55\times10^{-2}$ & $8.38\times10^{-2}$ &  $4.76\times10^{-2}$ & $6.50\times10^{-2}$ & $9.64\times10^{-2}$&
			$2.42\times10^{-3}$ & $1.75\times10^{-3}$ &$1.34\times10^{-2}$\\
			\hline
		\end{tabular}
	\end{table*}
\end{center} 
\begin{table*}	
	\scriptsize
	\caption{$L_{\infty}$,$L_2$ and $H^1$ errors of ST-GFDM with different $m$ when $N_T=10460$ for Example 2}
	\begin{tabular}{cccccccccc}
		\hline
		\multirow{1}{*}{m} & \multicolumn{3}{c}{u} & \multicolumn{3}{c}{v}& \multicolumn{3}{c}{p} \\
		\hline
		
		& $L_{\infty}$      &  $L_2$   &   $H^1$
		
		& $L_{\infty}$      &  $L_2$   &   $H^1$
		
		& $L_{\infty}$      &  $L_2$   &   $H^1$\\
		\hline
		$55$ & $6.82\times10^{-5}$ & $3.14\times10^{-5}$ & $8.66\times10^{-4}$ &  $5.64\times10^{-5}$ & $2.47\times10^{-5}$ & $ 8.61\times10^{-4}$&
		$7.19\times10^{-4}$ & $8.34\times10^{-5}$ & $ 1.97\times10^{-3}$ \\

		$57$ & $9.12\times10^{-5}$ & $2.64\times10^{-5}$ & $ 1.11\times10^{-3}$ &  $ 4.93\times10^{-5}$ & $1.99\times10^{-5}$ & $7.29\times10^{-4}$&
		$9.74\times10^{-4}$ & $1.01\times10^{-4}$ & $1.82\times10^{-3}$ \\
		
		$59$ & $1.09\times10^{-4}$ & $3.03\times10^{-5}$ & $ 1.21\times10^{-3}$ &  $ 4.44\times10^{-5}$ & $1.87\times10^{-5}$ & $7.87\times10^{-4}$&
		$1.14\times10^{-3}$ & $1.11\times10^{-4}$ & $1.70\times10^{-3}$ \\
		
		$61$ & $4.85\times10^{-5}$ & $1.31\times10^{-5}$ & $ 6.33\times10^{-4}$ &  $   5.71\times10^{-5}$ & $1.88\times10^{-5}$ & $6.44\times10^{-4}$&
		$3.90\times10^{-4}$ & $4.30\times10^{-5}$ & $8.77\times10^{-4}$ \\

		$63$ & $1.26\times10^{-4}$ & $4.23\times10^{-5}$ & $ 1.43\times10^{-3}$ &  $ 8.80\times10^{-5}$ & $3.09\times10^{-5}$ & $9.80\times10^{-4}$&
		$1.30\times10^{-3}$ & $1.28\times10^{-4}$ & $1.78\times10^{-3}$ \\
		\hline
	\end{tabular}
\end{table*}
% Figure
\begin{figure}
	\centering
	\includegraphics[scale=.4]{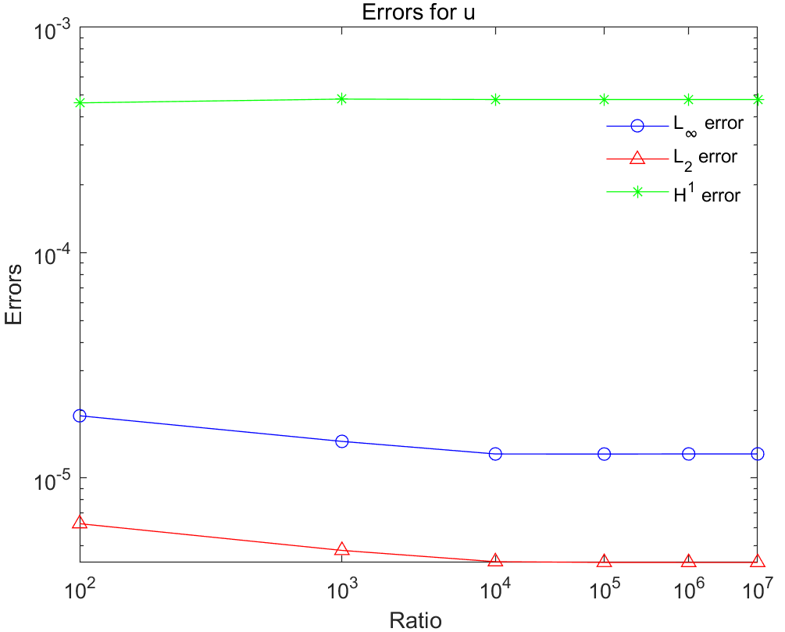}    
	\includegraphics[scale=.4]{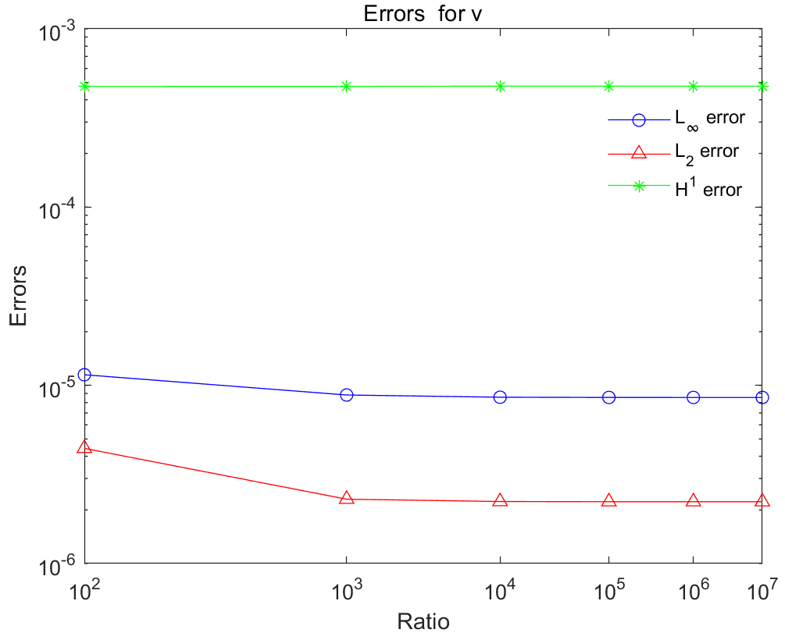}    
	\includegraphics[scale=.4]{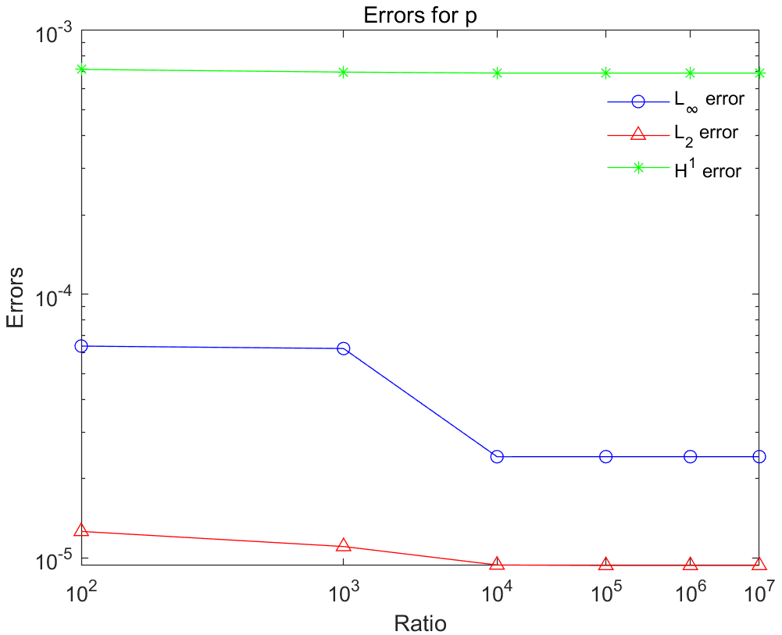}    
	\caption{$L_{\infty}$,$L_2$ and $H^1$ errors of ST-GFDM with different jump $Ratio=\frac{\beta_1}{\beta_2}=\frac{\rho_1}{\rho_2}=(10^2,10^3,10^4,10^5,10^6,10^7)$ when $N_{T}=10460$ for Example 2.}
\end{figure}
\begin{figure}
	\centering
	\includegraphics[scale=.4]{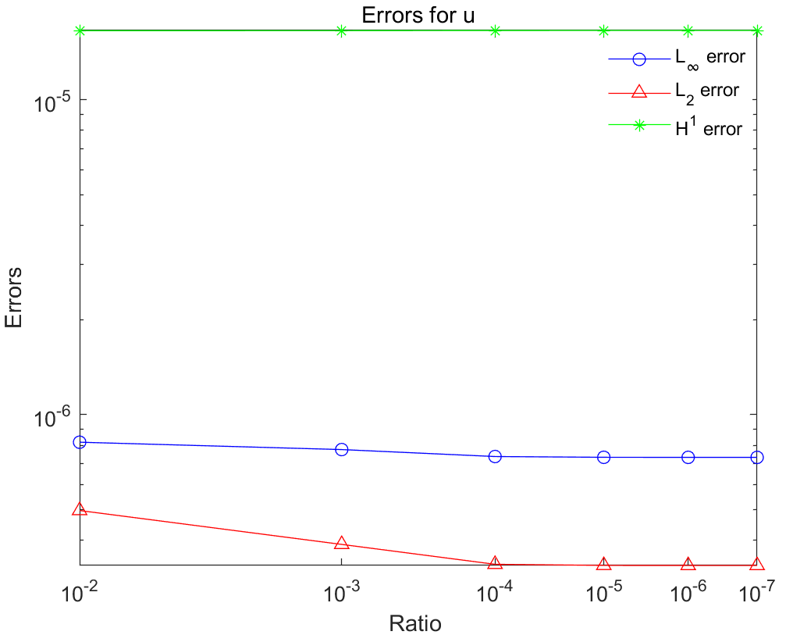}    
	\includegraphics[scale=.4]{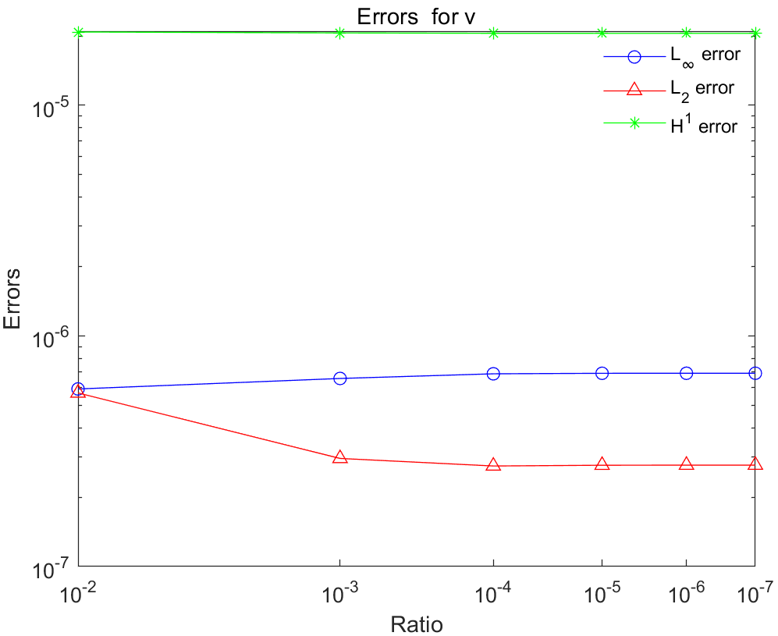}    
	\includegraphics[scale=.4]{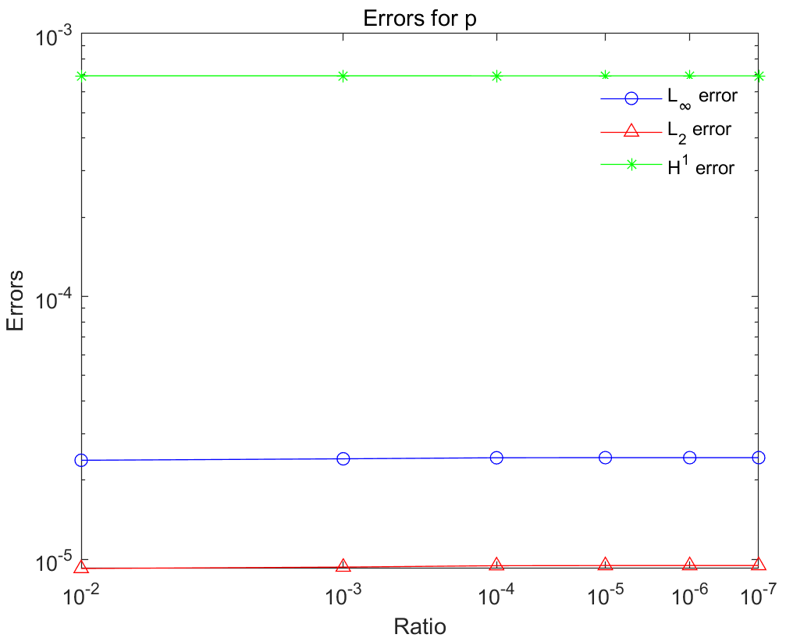}    
	\caption{$L_{\infty}$,$L_2$ and $H^1$ errors of ST-GFDM with different jump $Ratio=\frac{\beta_1}{\beta_2}=\frac{\rho_1}{\rho_2}=(10^{-2},10^{-3},10^{-4},10^{-5},10^{-6},10^{-7})$ when $N_{T}=10460$ for Example 2.}
\end{figure}
\begin{figure}
	\centering
	\includegraphics[scale=.4]{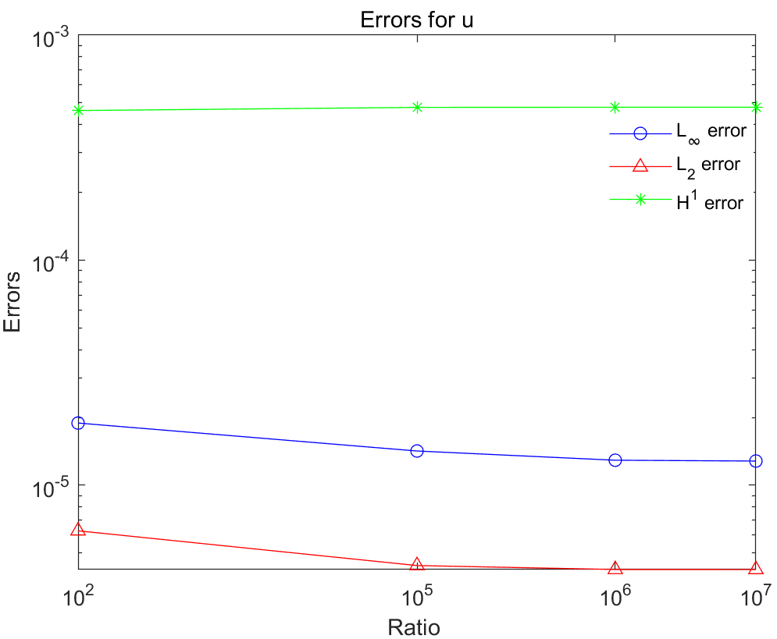}    
	\includegraphics[scale=.4]{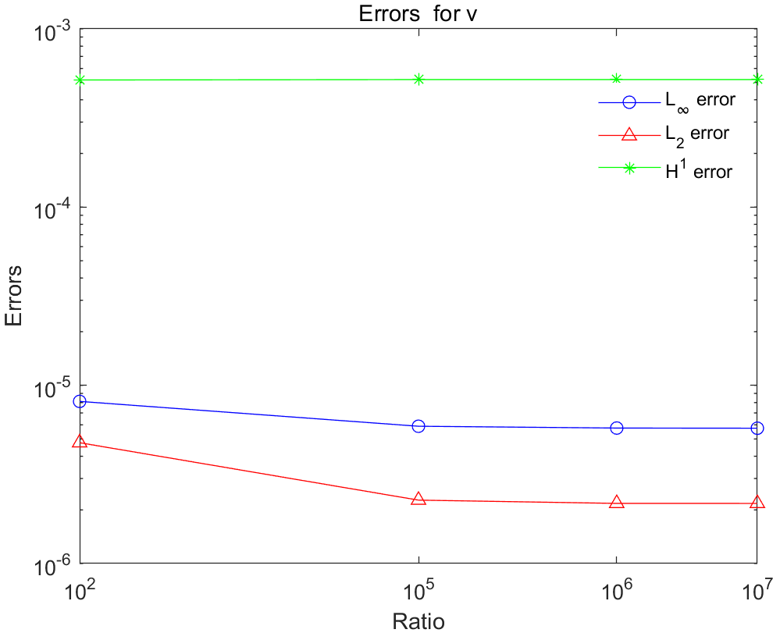}    
	\includegraphics[scale=.4]{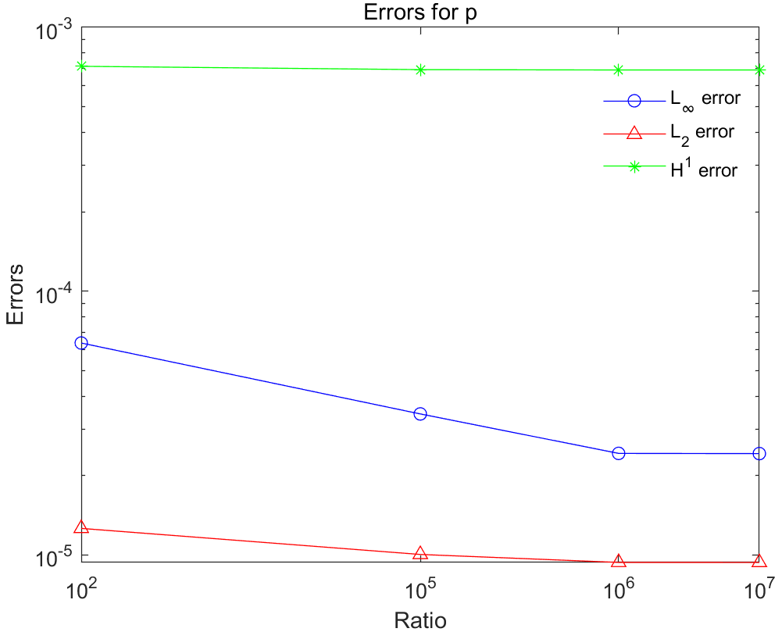}    
	\caption{$L_{\infty}$,$L_2$ and $H^1$ errors of ST-GFDM with different jump $Ratio=\frac{\beta_1}{\beta_2}=(10^2,10^3,10^4,10^5,10^6,10^7)$ when $\frac{\rho_1}{\rho_2}=100$, $N_{T}=10460$ for Example 2.}
\end{figure}
Fig.12 presents the point collocation at $t=0$ and the point collocation for all times, we can see the initial interface shape and the 3D interface shape of this example. In Table. 5, $L_{\infty}$,$L_{2}$ and $H^1$ errors of ST-GFDM when $N_T=10460$ are provided to show the efficiency of the ST-GFDM for the Stokes/Parabolic moving interface problem with an irregular moving circle interface. We can see that the numerical results are accurate. The contour of the absolute error is shown in Fig.13. We can see the relative large errors appear on the interface. In order to test the stability of ST-GFDM for Stokes/Parabolic moving  interface problems, the number 'm' and the jump ratio are changed. From Table 6, we can see that the numerical errors are still stable when the number of 'm' is changed. In Fig. 14, Fig. 15 and Fig. 16, we can see that the $L_{\infty}$,$L_2$ and $H^1$ errors are stable when the jump $Ratio=\frac{\beta_1}{\beta_2}=\frac{\rho_1}{\rho_2}$ and even in $Ratio=\frac{\beta_1}{\beta_2},\frac{\rho_1}{\rho_2}=100$. Therefore, the accuracy, stability and high efficiency of the ST-GFDM for the Stokes/Parabolic interface problem with an irregular moving circle interface are verified. Furthermore, this example shows the advantage of the ST-GFDM in dealing with the irregular moving directive. 
\subsection{Example 3: The Stokes/Parabolic interface problem with an irregular moving octagon interface.}
In this example, we consider the above problem with  an irregular moving octagon interface (see Fig.18(left)) $\Gamma_t:$\\
$\varphi_3=(x-0.3-0.1(t+sin(5t)))^2+(y-0.3-0.1(t+t^3))^2-(0.5(3/(10n^2))(1+2n+n^2-(n+1)cos(n\theta)))=0,$\\
$n=8,\theta\in [0,2\pi].$
\begin{figure}
	\centering
	\includegraphics[scale=.3]{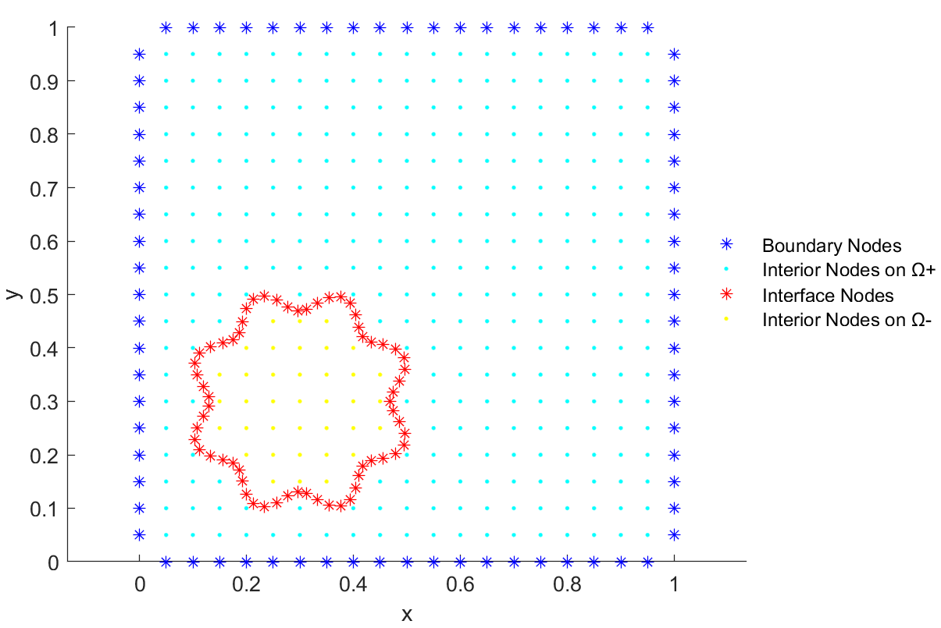}    
	\includegraphics[scale=.3]{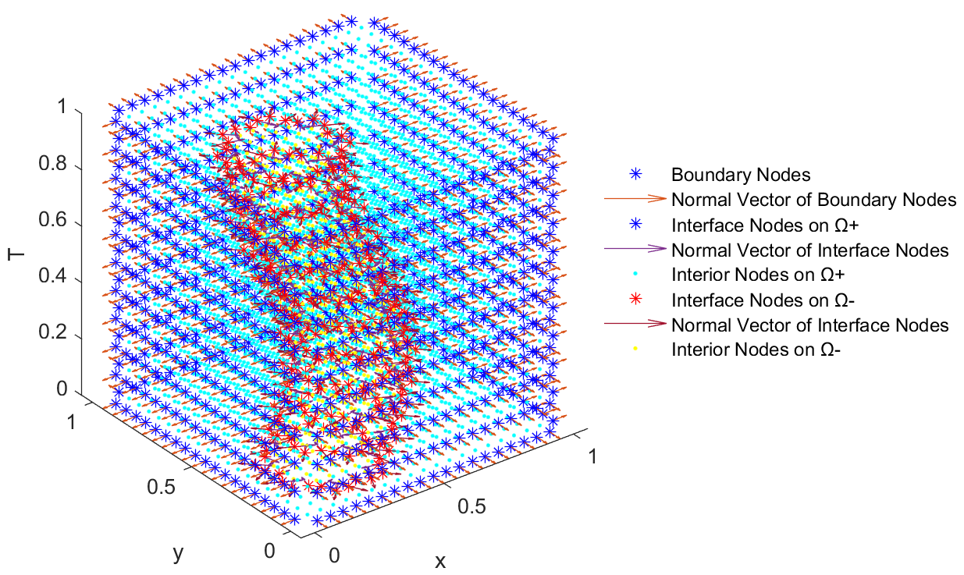}    
	\caption{ The point collocation at t=0 (left) and the point collocation for all time (right) for Example 3.}
\end{figure}% Figure
\begin{figure}
	\centering
	\includegraphics[scale=.4]{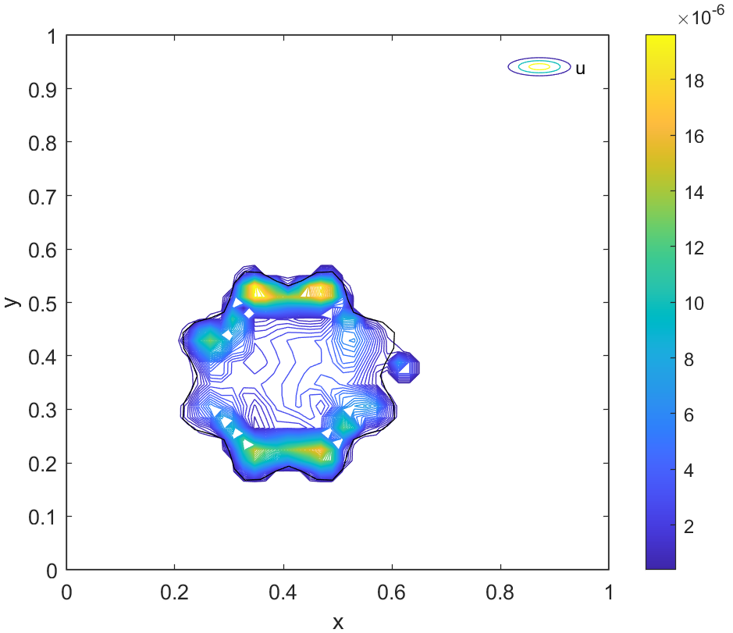}
	\includegraphics[scale=.4]{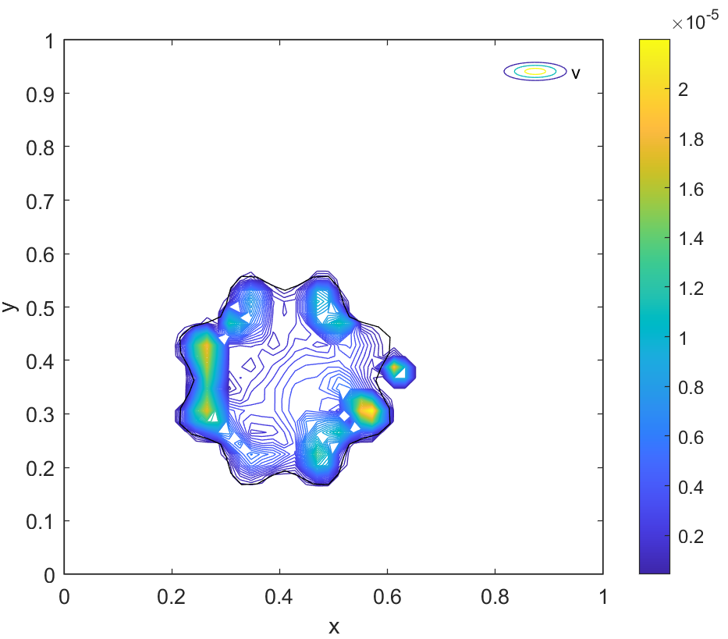}       
	\includegraphics[scale=.4]{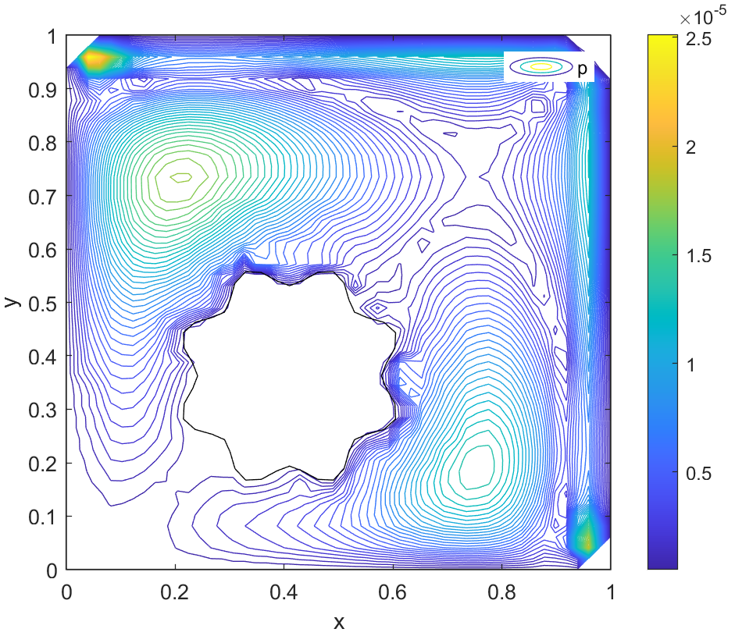}              
	\caption{ The contour of the absolute error when $N_{T}=7254$ for Example 3.}
\end{figure}
\begin{table*}	
	\scriptsize
	\caption{ $L_{\infty}$,$L_2$ and $H^1$ errors of ST-GFDM for Example 3}
	\begin{tabular}{cccccccccccc}
		\hline
		\multirow{1}{*}{$ $}&\multirow{1}{*}{$N_{T}$} & \multicolumn{3}{c}{u} & \multicolumn{3}{c}{v}& \multicolumn{3}{c}{p} \\   
		\hline
		
		&$ $ & $L_{\infty}$      &  $L_2$   &   $H^1$
		
		& $L_{\infty}$      &  $L_2$   &   $H^1$
		
		& $L_{\infty}$      &  $L_2$   &   $H^1$ \\
		\hline
		ST-GFDM &$3405$ & $1.21\times10^{-4}$ & $2.91\times10^{-5}$ & $ 2.47\times10^{-3}$ &  $5.04\times10^{-5}$ & $2.20\times10^{-5}$ & $ 2.43\times10^{-3}$&
		$ 1.22\times10^{-4}$ & $2.00\times10^{-5}$ & $1.18\times10^{-3}$  \\
		
		$ $&$7254$ & $ 3.50\times10^{-5}$ & $1.17\times10^{-5}$ & $1.58\times10^{-3}$ &  $4.86\times10^{-5}$ & $9.90\times10^{-6}$ & $1.57\times10^{-3}$&
		$1.81\times10^{-4}$ & $ 2.06\times10^{-5}$ & $8.10\times10^{-4}$  \\
		\hline
	\end{tabular}
\end{table*}

\begin{center}
	\begin{table*}	
		\scriptsize
		\caption{ $L_{\infty, relative}$,$L_{2,relative}$ and $H^1_{relative}$ errors of ST-GFDM for Example 3}
		\begin{tabular}{cccccccccccc}
			\hline
			\multirow{1}{*}{$N_{T}$} & \multicolumn{3}{c}{u} & \multicolumn{3}{c}{v}& \multicolumn{3}{c}{p} & \multicolumn{1}{c}{$Time(s)$}\\  
			\hline		
			$ $ & $L_{\infty, relative}$  &$L_{2,relative}$   &   $H^1_{relative}$
			
			& $L_{\infty, relative}$  &$L_{2,relative}$   &   $H^1_{relative}$
			
			& $L_{\infty, relative}$  &$L_{2,relative}$   &   $H^1_{relative}$ & $ $\\
			\hline

			$3405$ & $3.44\times10^{-2}$ & $1.91\times10^{-2}$ & $8.20\times10^{-2}$ &  $1.59\times10^{-2}$ & $1.60\times10^{-2}$ & $8.52\times10^{-2}$&
			$1.62\times10^{-3}$ & $1.28\times10^{-3}$ &$2.12\times10^{-2}$& $7.16$\\
			
			$7254$ & $ 1.07\times10^{-2}$ & $9.31\times10^{-3}$ & $6.00\times10^{-2}$ & $1.21\times10^{-2}$ &$6.87\times10^{-3}$ & $5.93\times10^{-2}$&
			$1.10\times10^{-3}$ & $9.51\times10^{-4}$ & $1.08\times10^{-2}$ &$12.8$\\
			\hline
		\end{tabular}
	\end{table*}
\end{center} 

\begin{table*}	
	\scriptsize
	\caption{$L_{\infty}$,$L_2$ and $H^1$ errors of ST-GFDM with different $m$ when $N_T=7254$ for Example 3}
	\begin{tabular}{cccccccccc}
		\hline
		\multirow{1}{*}{m} & \multicolumn{3}{c}{u} & \multicolumn{3}{c}{v}& \multicolumn{3}{c}{p} \\
		\hline
		
		& $L_{\infty}$      &  $L_2$   &   $H^1$
		
		& $L_{\infty}$      &  $L_2$   &   $H^1$
		
		& $L_{\infty}$      &  $L_2$   &   $H^1$\\
		\hline
		$55$ & $3.91\times10^{-4}$ & $1.35\times10^{-4}$ & $ 3.85\times10^{-3}$ &  $3.52\times10^{-4}$ & $1.30\times10^{-4}$ & $4.41\times10^{-3}$&
		$6.92\times10^{-3}$ & $8.01\times10^{-4}$ & $2.50\times10^{-2}$ \\

		$57$ & $3.27\times10^{-5}$ & $1.10\times10^{-5}$ & $ 1.56\times10^{-3}$ &  $ 3.96\times10^{-5}$ & $9.60\times10^{-6}$ & $1.54\times10^{-3}$&
		$2.39\times10^{-4}$ & $ 2.98\times10^{-5}$ & $1.05\times10^{-3}$ \\
		
		$59$ & $3.31\times10^{-5}$ & $1.23\times10^{-5}$ & $ 1.58\times10^{-3}$ &  $  4.25\times10^{-5}$ & $1.20\times10^{-5}$ & $1.57\times10^{-3}$&
		$3.82\times10^{-4}$ & $4.04\times10^{-5}$ & $1.33\times10^{-3}$ \\
		
		$61$ & $ 9.30\times10^{-5}$ & $3.90\times10^{-5}$ & $ 1.85\times10^{-3}$ &  $ 7.50\times10^{-5}$ & $ 3.03\times10^{-5}$ & $1.81\times10^{-3}$&
		$1.20\times10^{-3}$ & $1.52\times10^{-4}$ & $4.11\times10^{-3}$ \\
		
		$63$ & $3.35\times10^{-5}$ & $1.39\times10^{-5}$ & $ 1.64\times10^{-3}$ &  $ 5.07\times10^{-5}$ & $1.38\times10^{-5}$ & $1.64\times10^{-3}$&
		$3.33\times10^{-4}$ & $4.00\times10^{-5}$ & $1.20\times10^{-3}$ \\
		
		\hline
	\end{tabular}
\end{table*}
\begin{figure}
	\centering
	\includegraphics[scale=.4]{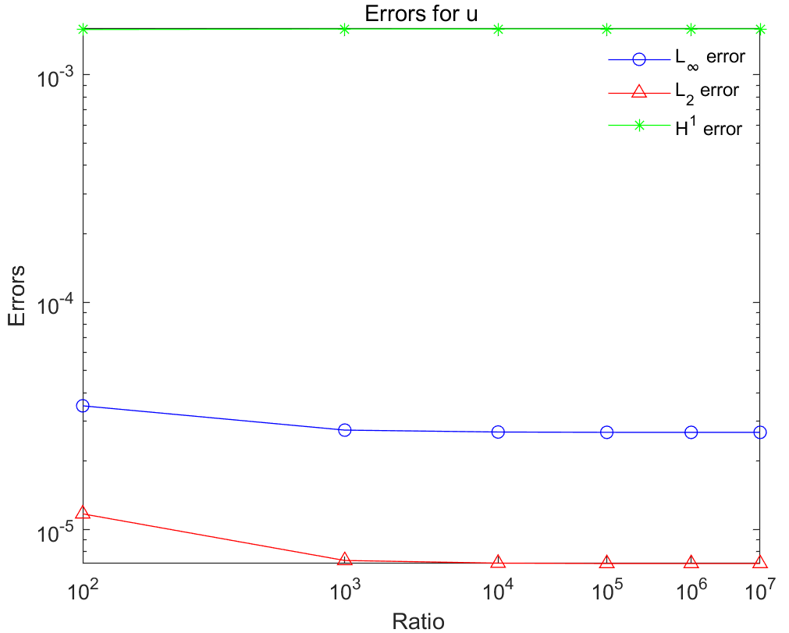}    
	\includegraphics[scale=.4]{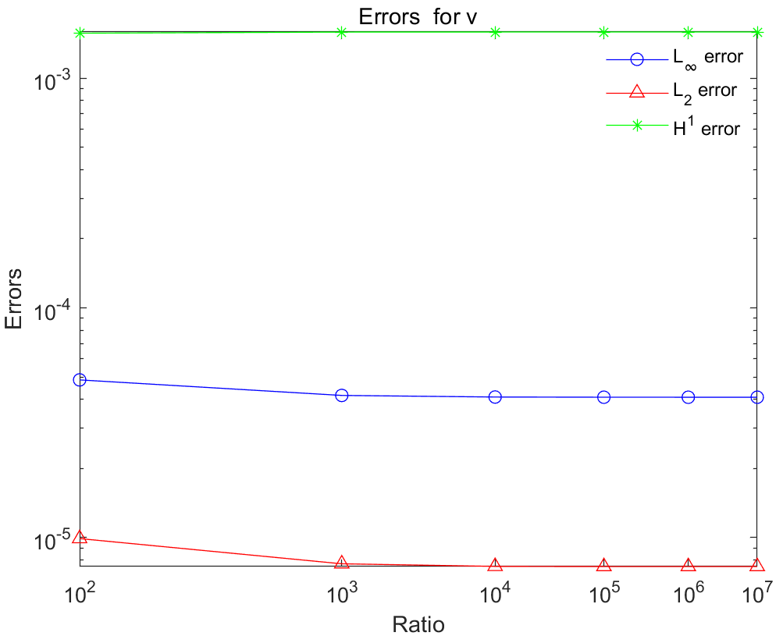}   
	\includegraphics[scale=.4]{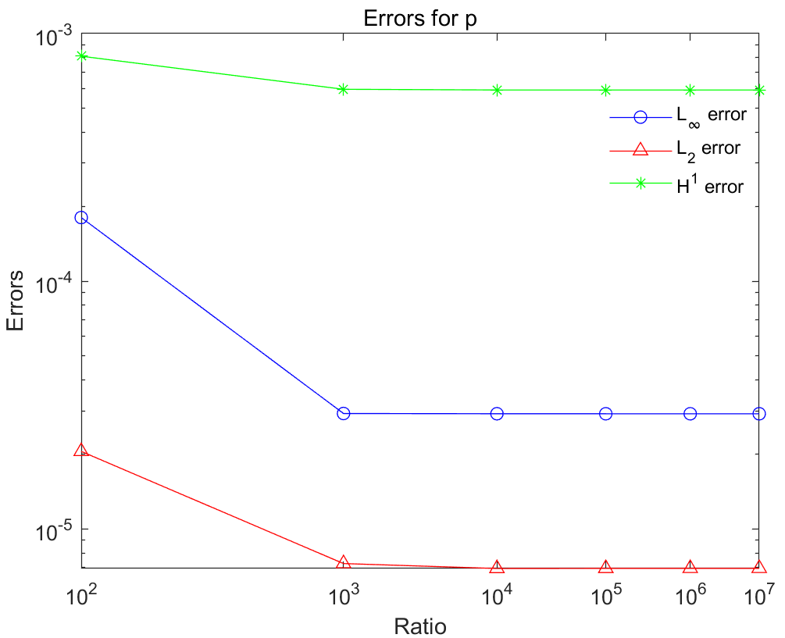}    
	\caption{$L_{\infty}$, $L_2$ and $H^1$ errors of ST-GFDM with different jump $Ratio=\frac{\beta_1}{\beta_2}=\frac{\rho_1}{\rho_2}=(10^2,10^3,10^4,10^5,10^6,10^7)$ when $N_{T}=7254$ for Example 3.}
\end{figure}
\begin{figure}
	\centering
	\includegraphics[scale=.4]{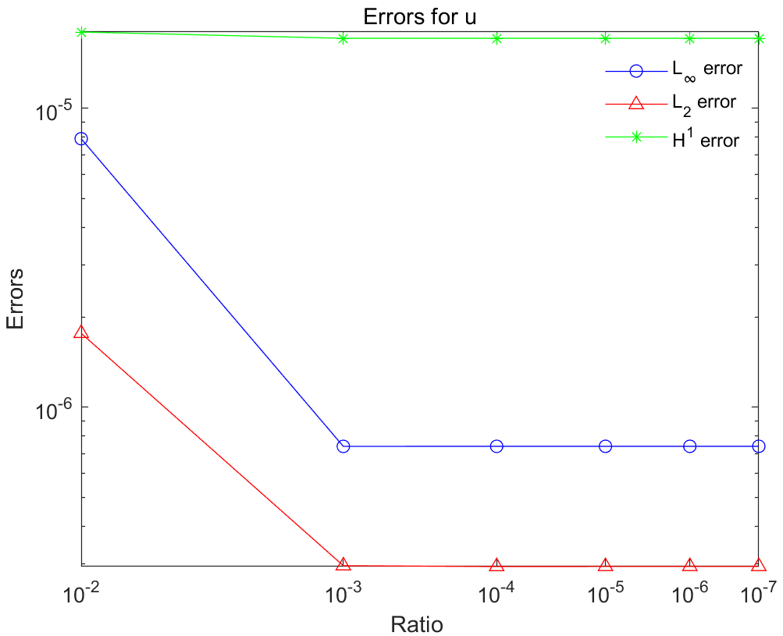}    
	\includegraphics[scale=.4]{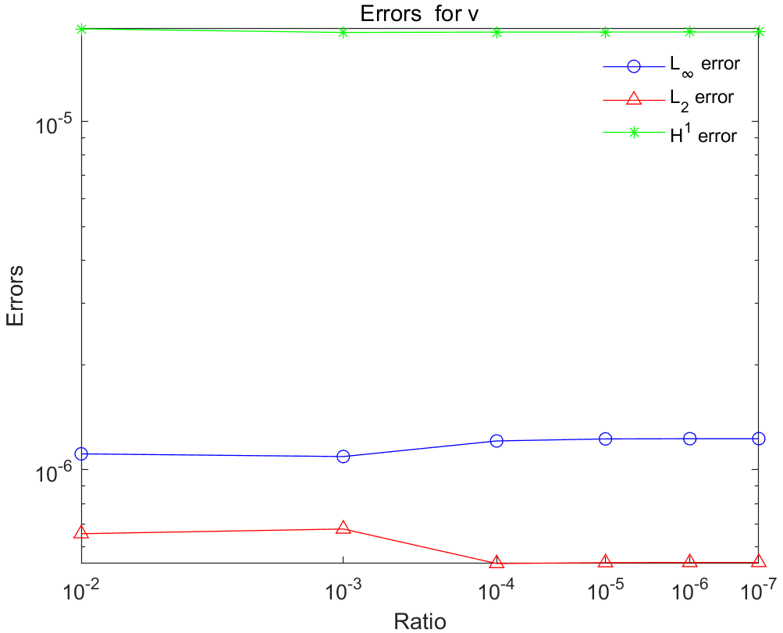}   
	\includegraphics[scale=.4]{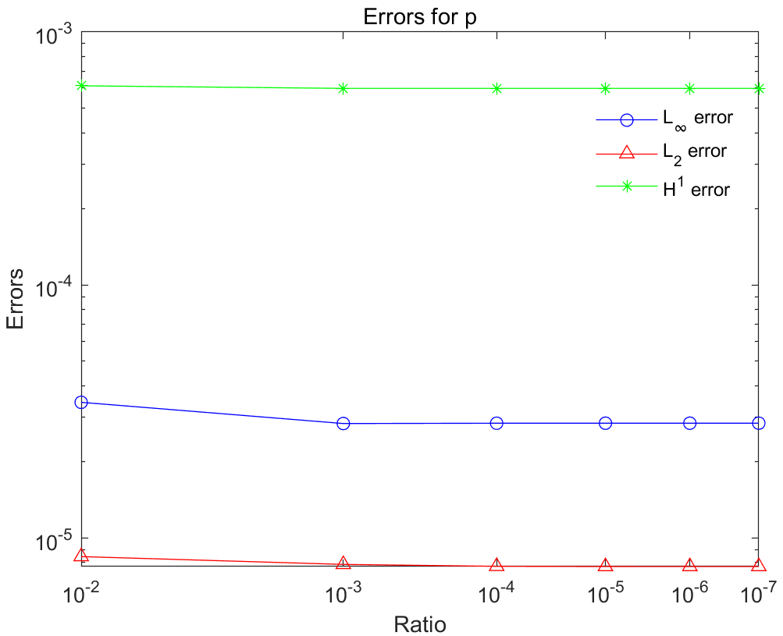}    
	\caption{$L_{\infty}$,$L_2$ and $H^1$ errors of ST-GFDM with different jump $Ratio=\frac{\beta_1}{\beta_2}=\frac{\rho_1}{\rho_2}=(10^{-2},10^{-3},10^{-4},10^{-5},10^{-6},10^{-7})$ when $N_{T}=7254$ for Example 3.}
\end{figure}
\begin{figure}
	\centering
	\includegraphics[scale=.4]{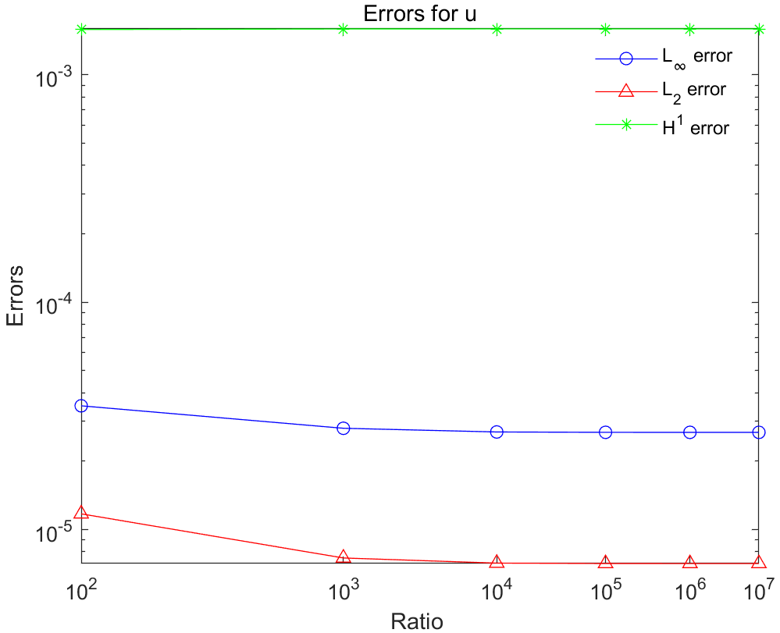}    
	\includegraphics[scale=.4]{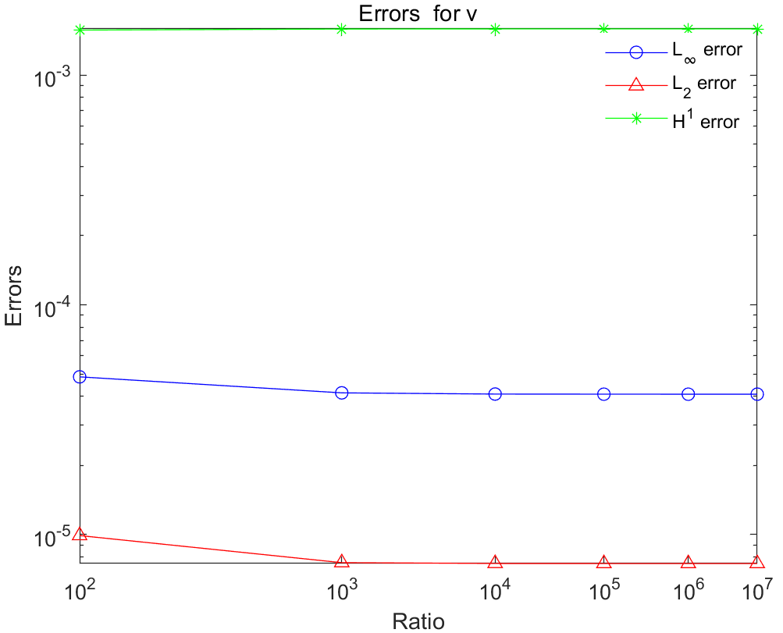}    
	\includegraphics[scale=.4]{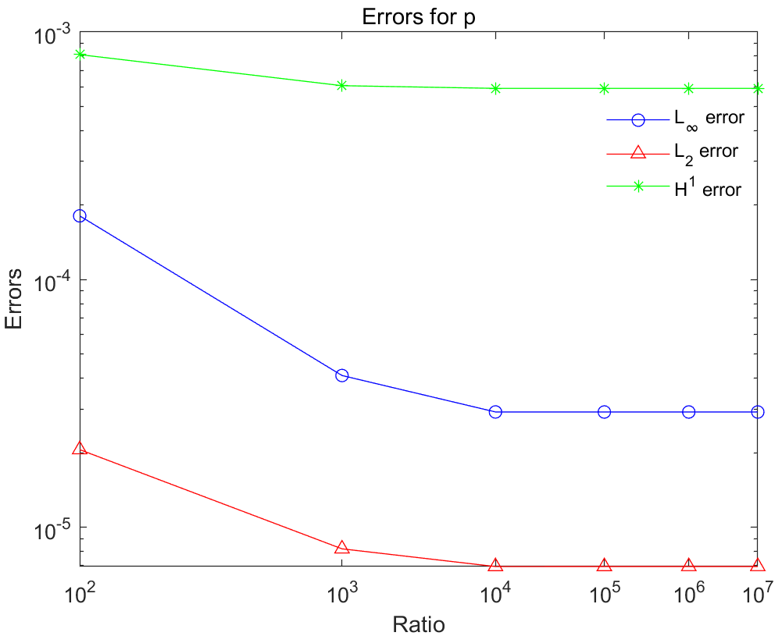}    
	\caption{$L_{\infty}$,$L_2$ and $H^1$ errors of ST-GFDM with different jump $Ratio=\frac{\beta_1}{\beta_2}=(10^2,10^3,10^4,10^5,10^6,10^7)$ when $\frac{\rho_1}{\rho_2}=100$, $N_{T}=7254$ for Example 3.}
\end{figure}

Fig.17 presents the point collocation at $t=0$ and the point collocation for all times, we can see the initial interface shape, moving direction and the final 3D interface shape of this example. In Table. 7, the $L_{\infty, relative}$,$L_{2,relative}$ and $H^1_{relative}$ errors of ST-GFDM  are provided in Table.8. We can see that our results are accurate and high efficiency. The contour of the absolute error is shown in Fig.18. We can see that the relatively large error appears on the interface for displacement. In order to test the stability of ST-GFDM for Stokes/Parabolic moving interface problems, the number 'm' and the jump ratio are changed. From Table 9, we can see that the numerical errors are still stable when the number of 'm' is changed. In Fig. 19, Fig. 20 and Fig. 21, we can see that the $L_{\infty}$,$L_2$ and $H^1$ errors are stable when the jump $Ratio=\frac{\beta_1}{\beta_2}=\frac{\rho_1}{\rho_2}$ and even in $Ratio=\frac{\beta_1}{\beta_2},\frac{\rho_1}{\rho_2}=100$. Therefore, the accuracy, stability and high efficiency of the ST-GFDM for the Stokes/Parabolic interface problem with irregular moving octagon interface are verified.This example shows the advantage of the ST-GFDM in dealing with the complex moving interface. 
\subsection{Example 4: The Stokes/Parabolic interface problem with a deformation of three-petalled flower interface.}
In this example, we consider the Stokes/Parabolic interface problem with a deformation of three-petalled flower interface (From Ref.[12] and Ref. [31]), the initial interface (see Fig.22(left)) is $\Gamma_t^4:r=0.4(0.8+0.2sin(3\theta)), 0\leq \theta \leq 2\pi.$ The interface will relax to its equilibrium, a circle with radius $r_0=0.2.$ The exact solution are
\begin{eqnarray}
	u&=&(u_1,u_2)=(y-wt)sin((x-wt)^2+(y-wt)^2-0.0625)sin(t)/\beta,\\
	v&=&(v_1,v_2)=-(x-wt)sin((x-wt)^2+(y-wt)^2-0.0625)sin(t)/\beta, \\
	p&=&0.1(x^3-y^3)((x-wt)^2+(y-wt)^2-0.0625).
\end{eqnarray}
In this example, we take $dt=\frac{1}{N_x}=\frac{1}{20}.$ The coefficient $\beta_1=1000,\beta_2=\rho_1=\rho_2=1, w=0.1.$ The vector$\mathbf{u}=(u,v)$ and the velocity field $u=\sqrt{(u^2+v^2)}$ is the contour of the vector $\mathbf{u}.$
\begin{figure}
	\centering
	\includegraphics[scale=.5]{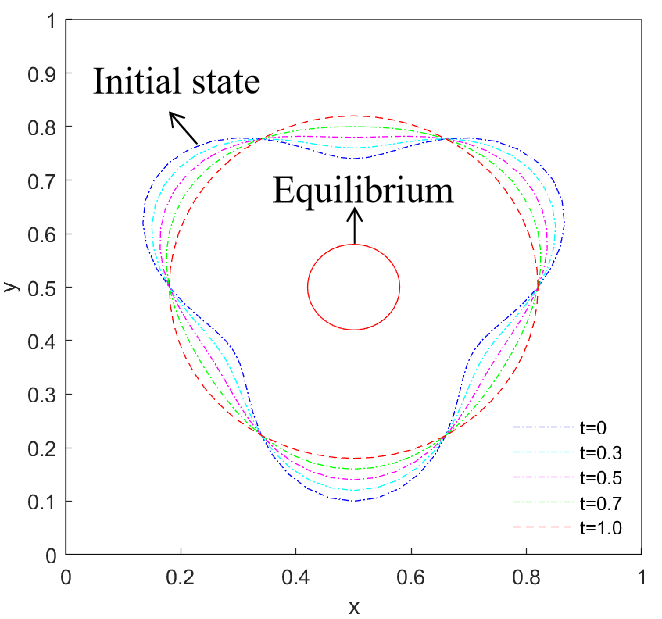}    
	\includegraphics[scale=.5]{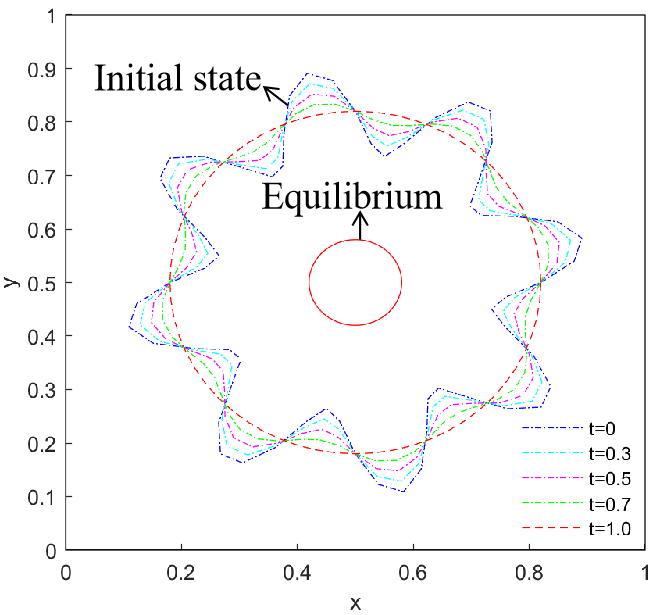}    
	\caption{ The interface collocation at different times in a square domain (Three-petaled flower interface(left) and Eight-petaled flower interface(right))(From Ref.[31]).}
\end{figure}
\begin{table*}	
	\scriptsize
	\caption{ $L_{\infty}$,$L_2$ and $H^1$ errors of ST-GFDM for Example 4}
	\begin{tabular}{ccccccccccc}
		\hline
		\multirow{1}{*}{$t$} & \multicolumn{3}{c}{u} & \multicolumn{3}{c}{v}& \multicolumn{3}{c}{p} & \multicolumn{1}{c}{$Time(s)$}\\
		\hline
		
		$ $& $L_{\infty}$      &  $L_2$   &   $H^1$
		
		& $L_{\infty}$      &  $L_2$   &   $H^1$
		
		& $L_{\infty}$      &  $L_2$   &   $H^1$ & $ $\\
		\hline
		
		$t=0.1$& $4.06\times10^{-5}$ & $ 1.85\times10^{-5}$ & $3.67\times10^{-4}$ &  $4.15\times10^{-5}$ & $1.95\times10^{-5}$ & $3.76\times10^{-4}$&
		$4.16\times10^{-4}$ & $1.31\times10^{-4}$ & $2.35\times10^{-3}$ & $ 18.2$ \\
		
		$t=0.5$& $1.26\times10^{-4}$ & $5.96\times10^{-5}$ & $1.06\times10^{-3}$ &  $1.30\times10^{-4}$ & $6.08\times10^{-5}$ & $1.11\times10^{-3}$&
		$1.87\times10^{-3}$ & $6.07\times10^{-4}$ & $1.14\times10^{-2}$ & $ 19.5$ \\
		
		$t=1.0$& $3.51\times10^{-4}$ & $1.65\times10^{-4}$ & $2.35\times10^{-3}$ &  $3.51\times10^{-4}$ & $1.65\times10^{-4}$ & $2.35\times10^{-3}$&
		$3.14\times10^{-3}$ & $1.10\times10^{-3}$ & $2.05\times10^{-2}$ & $ 17.5$ \\
		\hline
	\end{tabular}
\end{table*}

\begin{center}
	\begin{table*}	
		\scriptsize
		\caption{ $L_{\infty, relative}$,$L_{2,relative}$ and $H^1_{relative}$ errors of ST-GFDM for Example 4}
		\begin{tabular}{ccccccccccc}
			\hline
			\multirow{1}{*}{$t$} & \multicolumn{3}{c}{u} & \multicolumn{3}{c}{v}& \multicolumn{3}{c}{p} & \multicolumn{1}{c}{$Time(s)$}\\  
			\hline		
			& $L_{\infty, relative}$  &$L_{2,relative}$   &   $H^1_{relative}$
			
			& $L_{\infty, relative}$  &$L_{2,relative}$   &   $H^1_{relative}$
			
			& $L_{\infty, relative}$  &$L_{2,relative}$   &   $H^1_{relative}$ & $ $\\
			\hline

			$t=0.1$ & $7.35\times10^{-2}$ & $2.27\times10^{-2}$ & $1.31\times10^{-1}$ &  $7.35\times10^{-2}$ & $2.29\times10^{-2}$ & $1.24\times10^{-1}$&
			$4.11\times10^{-2}$ & $2.97\times10^{-2}$ & $8.05\times10^{-2}$& $14.1 $\\
			
			$t=0.5$ & $3.00\times10^{-2}$ & $9.19\times10^{-3}$ & $8.78\times10^{-2}$ &  $1.82\times10^{-2}$ & $8.70\times10^{-3}$ & $6.48\times10^{-2}$&
			$4.11\times10^{-2}$ & $3.08\times10^{-2}$ & $8.54\times10^{-2}$ &$15.4$\\
			
			$t=1.0$ & $5.26\times10^{-3}$ & $3.16\times10^{-3}$ & $1.20\times10^{-2}$ &  $5.26\times10^{-3}$ & $3.16\times10^{-3}$ & $1.20\times10^{-2}$&
			$3.96\times10^{-2}$ & $3.18\times10^{-2}$ & $8.61\times10^{-2}$ &$14.1$\\
			\hline
		\end{tabular}
	\end{table*}
\end{center} 

\begin{figure}
	\centering
	\includegraphics[scale=.4]{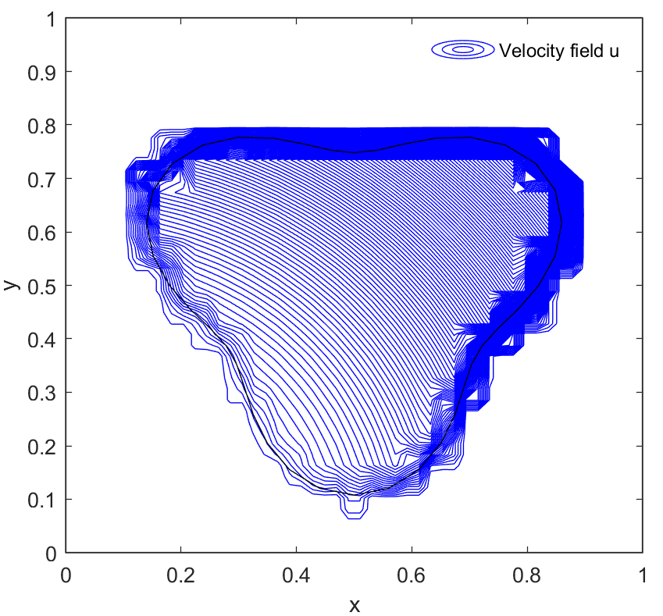}
	\includegraphics[scale=.4]{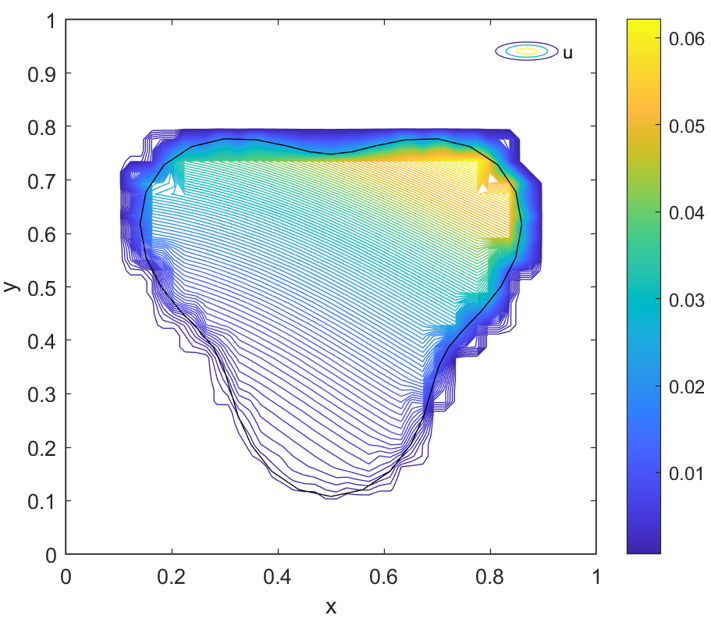}       
	\includegraphics[scale=.4]{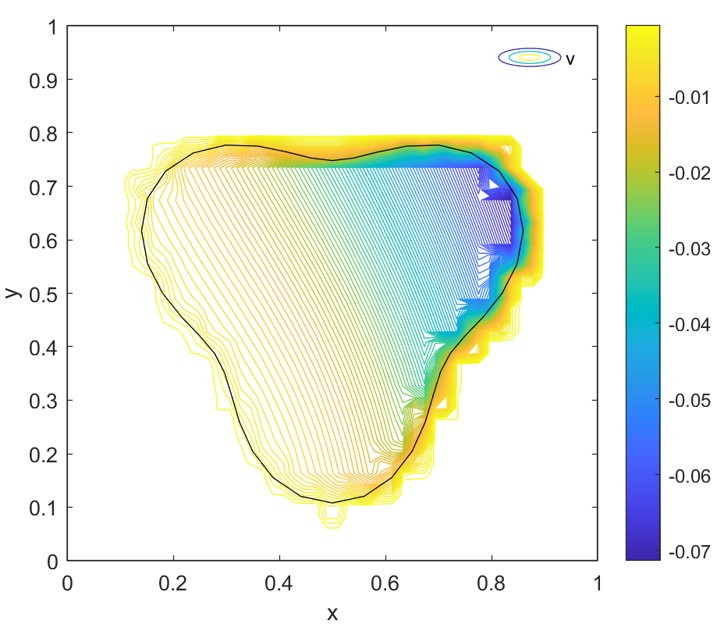}              
	\caption{ The contour of the numerical solution when $t=0.1, N_{T}=9996$ for Example 4.}
\end{figure}
\begin{figure}
	\centering
	\includegraphics[scale=.4]{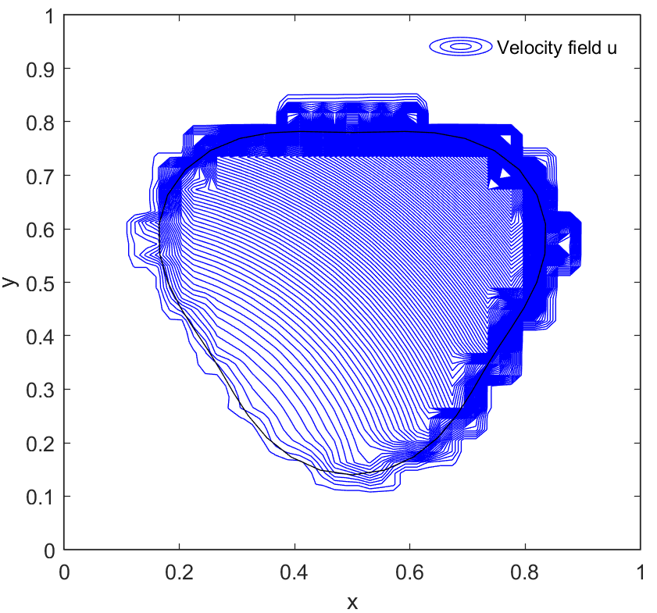}
	\includegraphics[scale=.4]{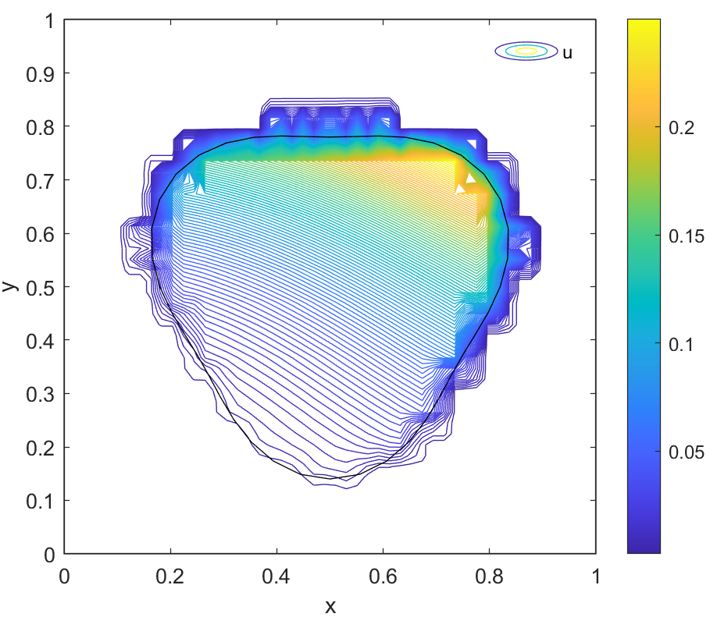}       
	\includegraphics[scale=.4]{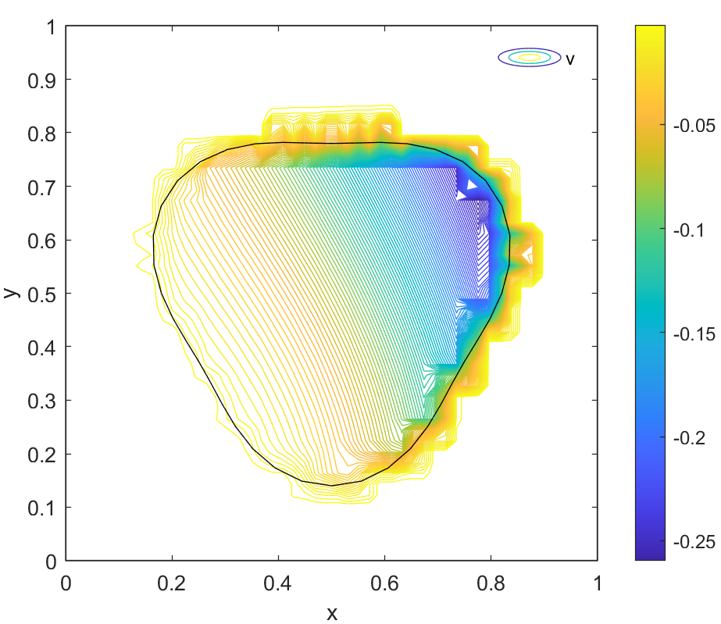}              
	\caption{ The contour of the numerical solution when $t=0.5,N_{T}=9156$ for Example 4.}
\end{figure}
\begin{figure}
	\centering
	\includegraphics[scale=.4]{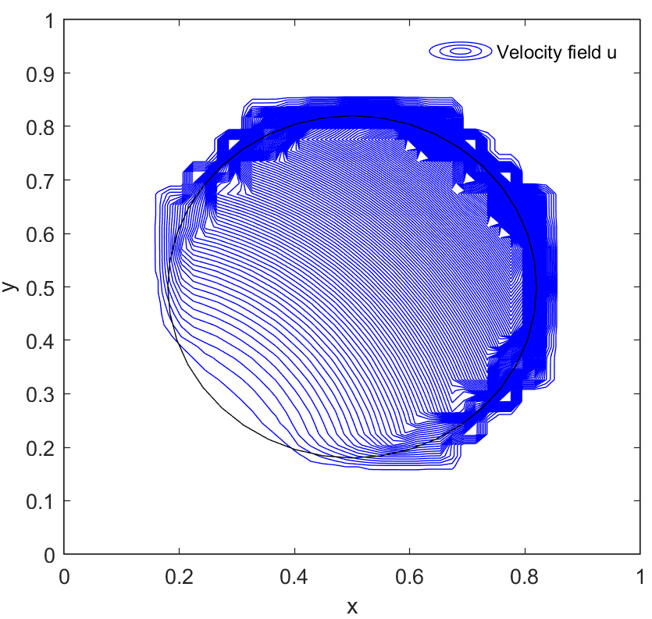}
	\includegraphics[scale=.4]{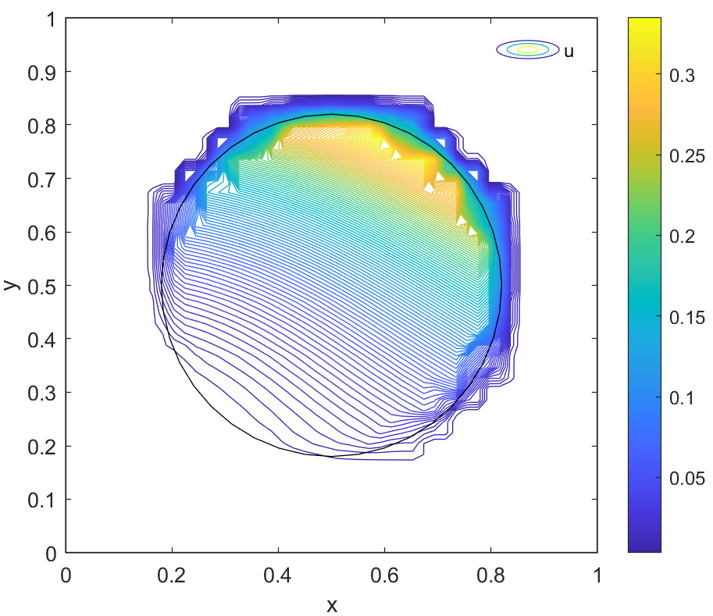}       
	\includegraphics[scale=.4]{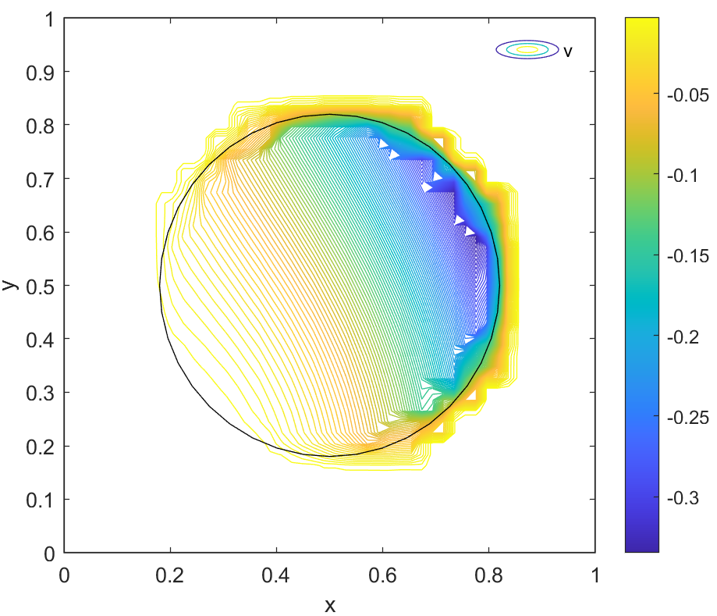}              
	\caption{ The contour of the numerical solution when $t=1, N_{T}=9828$ for Example 4.}
\end{figure}

\begin{figure}
	\centering
	\includegraphics[scale=.35]{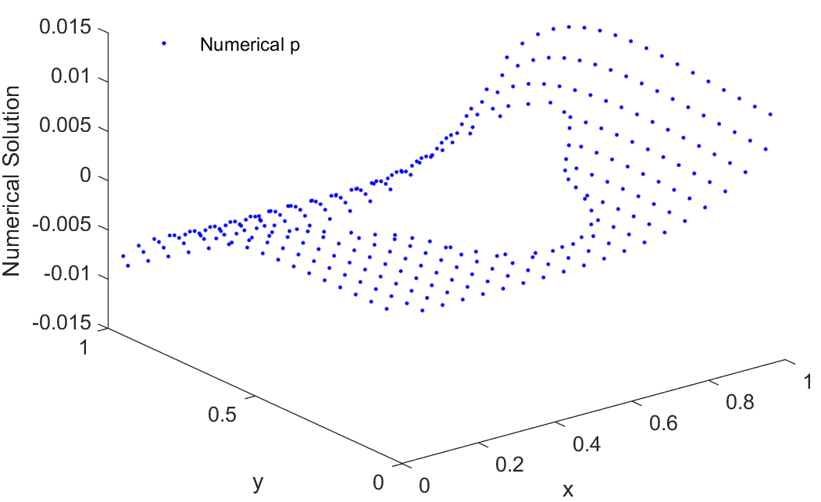}
	\includegraphics[scale=.35]{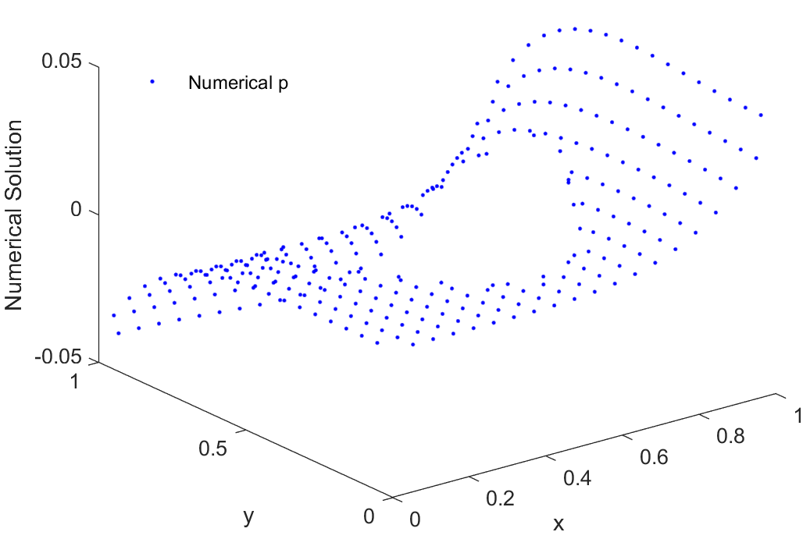}       
	\includegraphics[scale=.35]{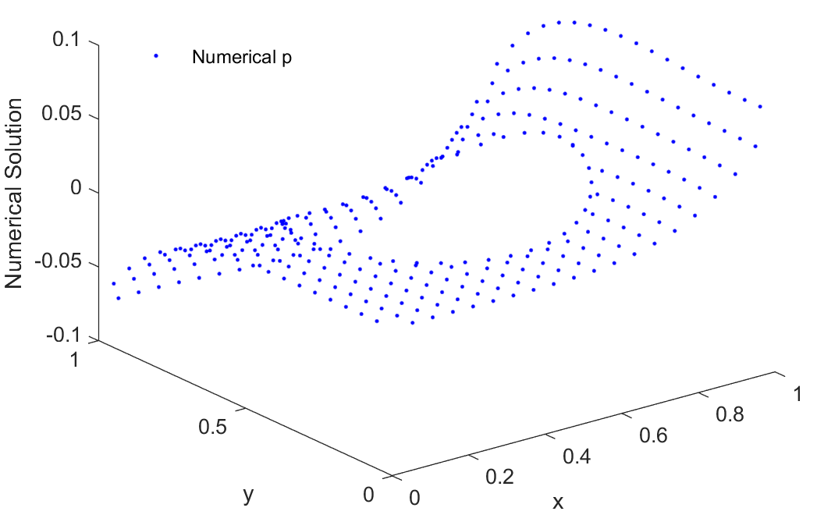}              
	\caption{ The numerical solution for p when $t=0.1$(left)$N_{T}=9996$, $t=0.5$(middle)$N_{T}=9156$ and $t=1$(right), $N_{T}=9828$ for Example 4.}
\end{figure}
Fig. 22 (left) shows the interface initial shape. In Table.10, $L_{\infty}$, $L_2$ and $H^1$ errors of ST-GFDM when $t=0.1, N_{T}=9996$, $t=0.5, N_{T}=9156$ and $t=1.0, N_{T}=9828$ are provided to show the accuracy and the stability of the ST-GFDM for solving the Stokes/Parabolic interface problem with the interface deformation. We can see that the $H^1$errors are all accurate and stable at different times. The $L_{\infty, relative}$,$L_{2,relative}$ and $H^1_{relative}$ errors of ST-GFDM are provided in Table.11.The contours of the numerical solution when $t=0.1$,$t=0.5$ and $t=1$ are presented in Fig. 23, Fig. 24 and Fig. 25, respectively. From these figures, we can see that the numerical solution changes with the interface deformation and we can capture the interface shape at any time. The numerical solution of the pressure is shown in Fig.26. We can see the distribution of the numerical solution. Therefore, the efficiency of the ST GFDM for solving the Stokes/Parabolic interface problem with a deformation of the three-petalled flower interface is verified. 
\subsection{Example 5: The Stokes/Parabolic interface problem with a deformation of eight-petalled flower interface.}
In this example, we consider the above Stokes/Parabolic interface problem with a deformation of eight-petalled flower interface (From Ref.[31]), the initial interface (see Fig.22(right)) is $\Gamma_t^5:r=0.4(0.8+0.2sin(8\theta)), 0\leq \theta \leq 2\pi.$ The interface will relax to its equilibrium, a circle with radius $r_0=0.2.$ In this example, we take$\beta_1=10000,\beta_2=\rho_1=\rho_2=1,$ $dt=\frac{1}{N_x}=\frac{1}{20}.$
\begin{table*}	
	\scriptsize
	\caption{ $L_{\infty}$,$L_2$ and $H^1$ errors of ST-GFDM for Example 5}
	\begin{tabular}{ccccccccccc}
		\hline
		\multirow{1}{*}{$ $} & \multicolumn{3}{c}{u} & \multicolumn{3}{c}{v}& \multicolumn{3}{c}{p} & \multicolumn{1}{c}{$Time(s)$}\\
		\hline
		
		& $L_{\infty}$      &  $L_2$   &   $H^1$
		
		& $L_{\infty}$      &  $L_2$   &   $H^1$
		
		& $L_{\infty}$      &  $L_2$   &   $H^1$ & $ $\\
		\hline
		
		$t=0.1$& $3.84\times10^{-5}$ & $1.54\times10^{-5}$ & $3.81\times10^{-4}$ &  $3.82\times10^{-5}$ & $1.52\times10^{-5}$ & $3.78\times10^{-4}$&
		$3.16\times10^{-4}$ & $9.96\times10^{-5}$ & $2.05\times10^{-3}$ & $  26.9$ \\
		
		$t=0.5$& $1.21\times10^{-4}$ & $5.27\times10^{-5}$ & $1.05\times10^{-3}$ &  $1.22\times10^{-4}$ & $5.25\times10^{-5}$ & $1.07\times10^{-3}$&
		$1.63\times10^{-3}$ & $5.34\times10^{-4}$ & $ 1.01\times10^{-2}$ & $ 21.1$ \\
		
		$t=1.0$& $3.45\times10^{-4}$ & $1.51\times10^{-4}$ & $2.18\times10^{-3}$ &  $3.45\times10^{-4}$ & $1.51\times10^{-4}$ & $ 2.18\times10^{-3}$&
		$3.14\times10^{-3}$ & $ 1.06\times10^{-3}$ & $2.08\times10^{-2}$ & $  22.4$ \\

		\hline
	\end{tabular}
\end{table*}
\begin{center}
	\begin{table*}	
		\scriptsize
		\caption{ $L_{\infty, relative}$,$L_{2,relative}$ and $H^1_{relative}$ errors of ST-GFDM for Example 5}
		\begin{tabular}{ccccccccccc}
			\hline
			\multirow{1}{*}{$t$} & \multicolumn{3}{c}{u} & \multicolumn{3}{c}{v}& \multicolumn{3}{c}{p} & \multicolumn{1}{c}{$Time(s)$}\\  
			\hline		
			& $L_{\infty, relative}$  &$L_{2,relative}$   &   $H^1_{relative}$
			
			& $L_{\infty, relative}$  &$L_{2,relative}$   &   $H^1_{relative}$
			
			& $L_{\infty, relative}$  &$L_{2,relative}$   &   $H^1_{relative}$ & $ $\\
			\hline

			$t=0.1$ & $1.33\times10^{-1}$ & $6.92\times10^{-2}$ & $3.20\times10^{-1}$ &  $1.66\times10^{-1}$ & $ 7.09\times10^{-2}$ & $3.27\times10^{-1}$&
			$3.15\times10^{-2}$ & $2.34
			\times10^{-2}$ & $7.20\times10^{-2}$& $16.8$\\
			
			$t=0.5$ & $ 4.18\times10^{-2}$ & $ 1.33\times10^{-2}$ & $1.03\times10^{-1}$ &  $  2.64\times10^{-2}$ & $1.18\times10^{-2}$ & $9.00\times10^{-2}$&
			$3.57\times10^{-2}$ & $2.77\times10^{-2}$ & $ 7.76\times10^{-2}$ &$18.9$\\
			
			$t=1.0$ & $8.53\times10^{-3}$ & $5.36\times10^{-3}$ & $2.97\times10^{-2}$ &  $8.53\times10^{-3}$ & $5.36\times10^{-3}$ & $2.97\times10^{-2}$&
			$3.95\times10^{-2}$ & $3.15\times10^{-2}$ & $8.94\times10^{-2}$ &$13.1$\\
			\hline
		\end{tabular}
	\end{table*}
\end{center} 
\begin{figure}
	\centering
	\includegraphics[scale=.4]{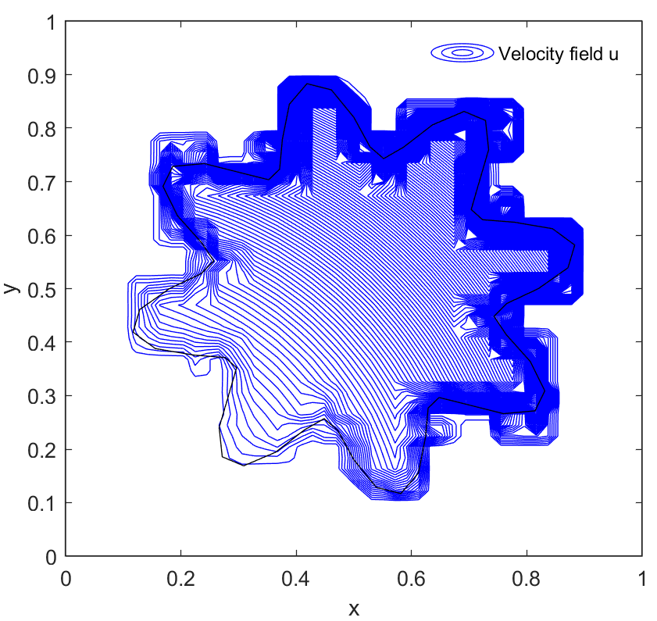}
	\includegraphics[scale=.4]{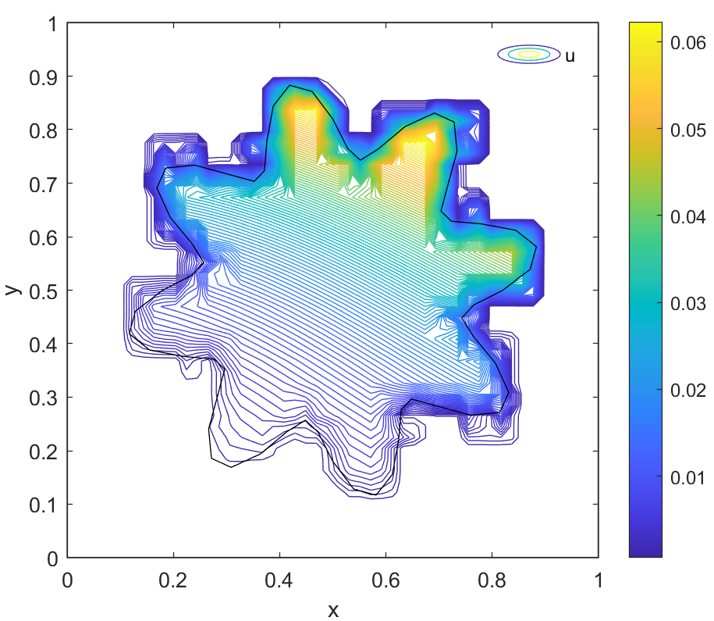}       
	\includegraphics[scale=.4]{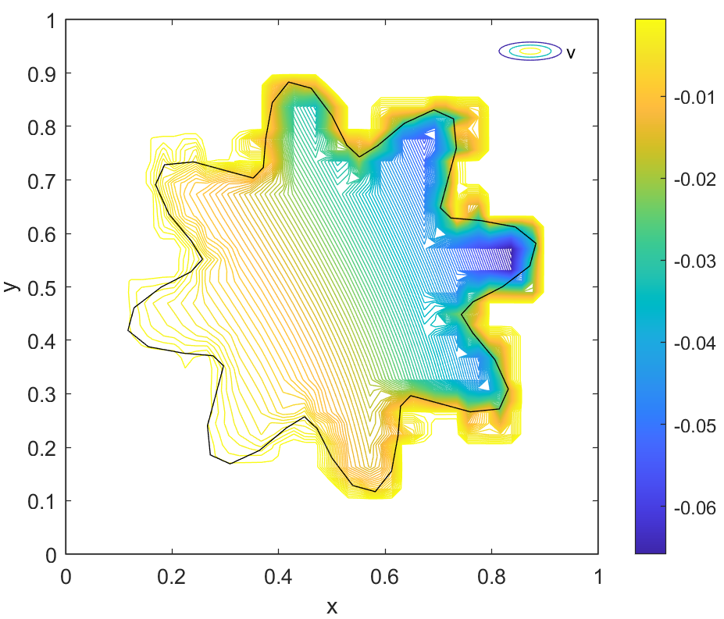}              
	\caption{ The contour of the numerical solution when $t=0.1, N_{T}=10752$ for Example 5.}
\end{figure}
\begin{figure}
	\centering
	\includegraphics[scale=.4]{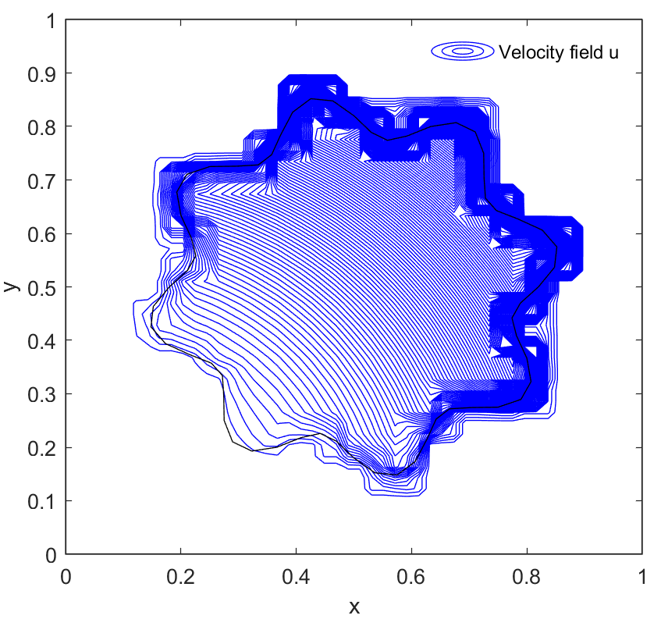}
	\includegraphics[scale=.4]{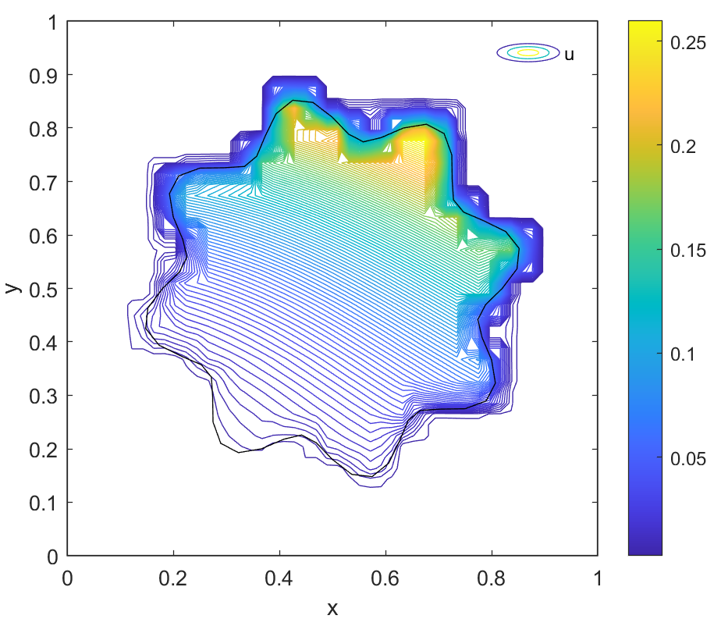}       
	\includegraphics[scale=.4]{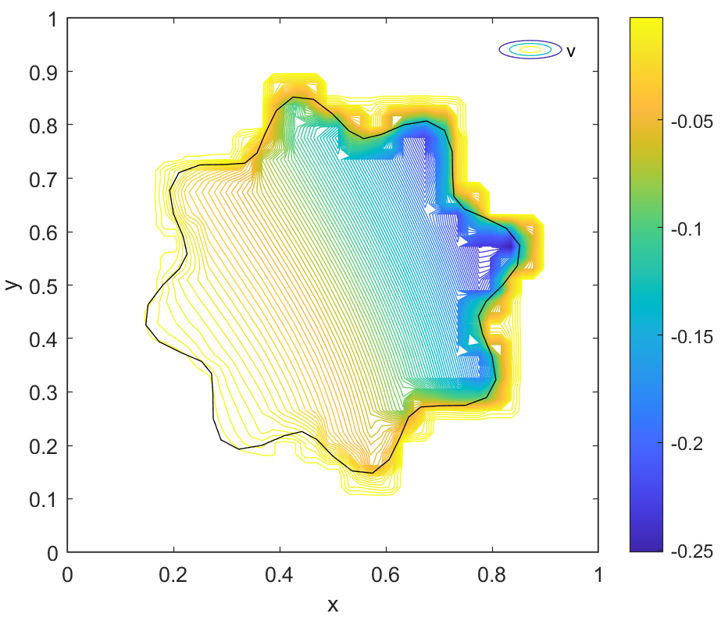}              
	\caption{ The contour of the numerical solution when $t=0.5, N_{T}=10584$ for Example 5.}
\end{figure}
\begin{figure}
	\centering
	\includegraphics[scale=.4]{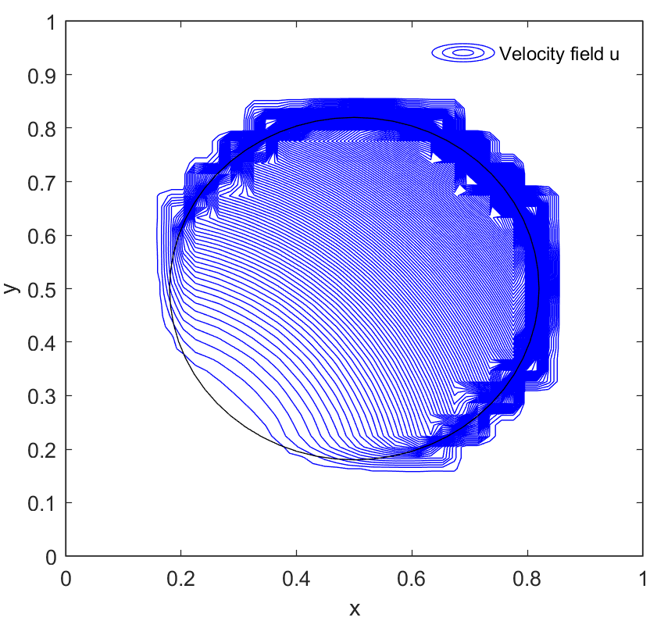}
	\includegraphics[scale=.4]{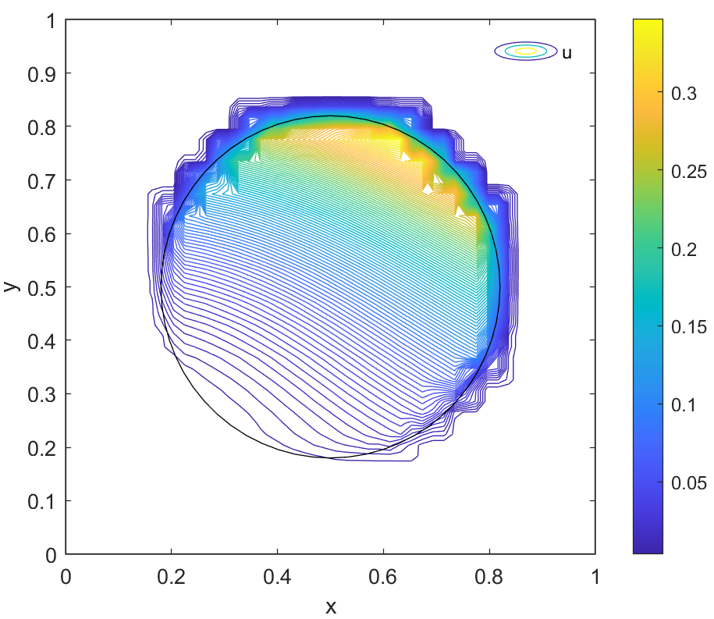}       
	\includegraphics[scale=.4]{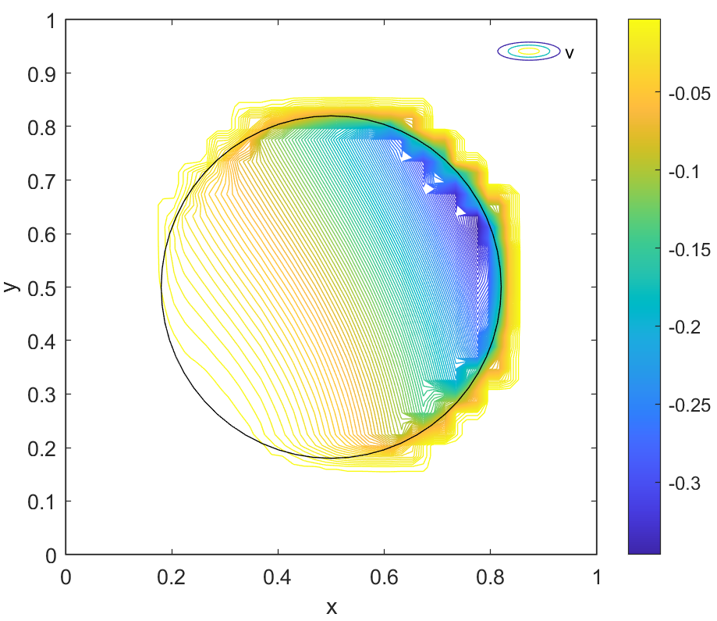}              
	\caption{ The contour of the numerical solution when $t=1, N_{T}=10668$ for Example 5.}
\end{figure}
\begin{figure}
	\centering
	\includegraphics[scale=.35]{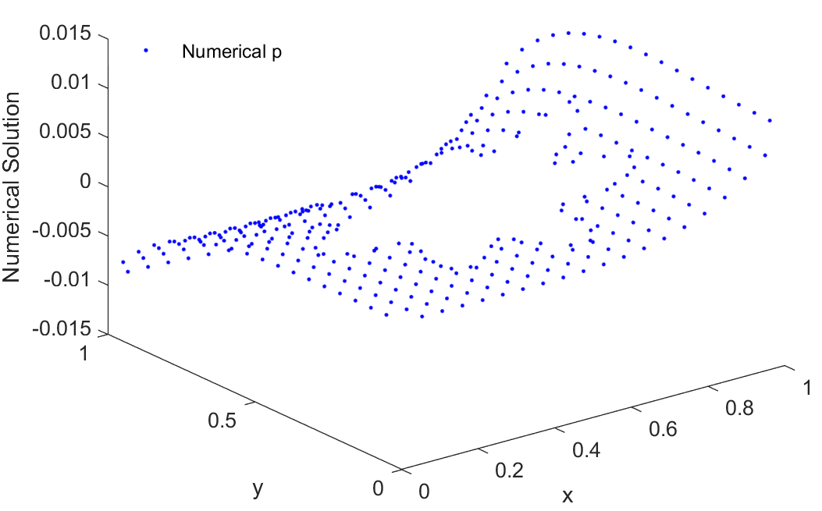}
	\includegraphics[scale=.35]{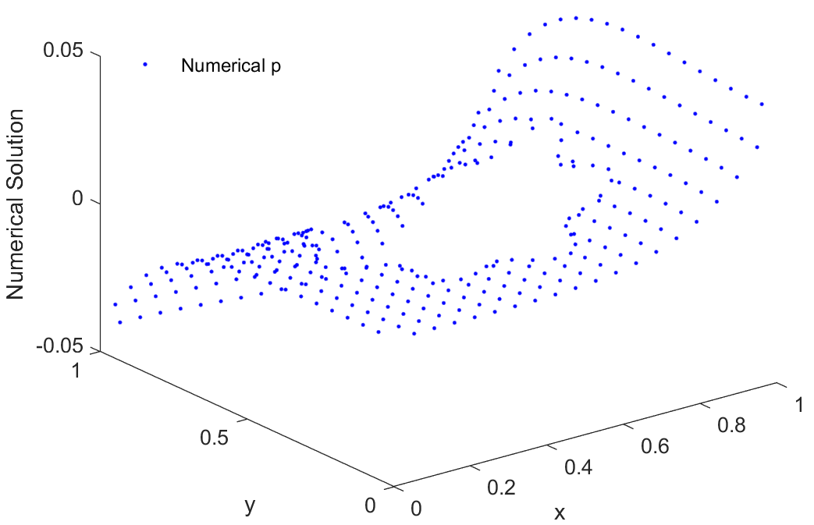}       
	\includegraphics[scale=.35]{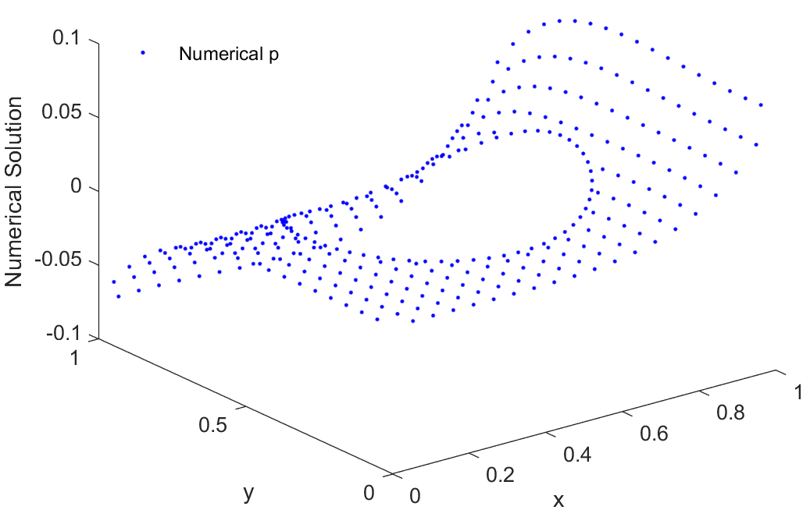}              
	\caption{ The numerical solution when $ t=0.1$(left), $t=0.5$(middle) and $t=1$(right), $N_{T}=10836$ for Example 5.}
\end{figure}
Fig. 22 (right) shows the initial shape of the interface. In Table.12, the $L_{\infty}$, $L_2$ and $H^1$ errors of ST-GFDM at $t=0.1, N_{T}=10752$, $t=0.5, N_{T}=10584$ and $t=1.0, N_{T}=10668$  are provided to show the accuracy and the stability of the ST-GFDM for solving the Stokes/Parabolic interface problem with the interface deformation. We can see that the $H^1$errors are all accurate and stable at different times. $L_{\infty, relative}$, $L_{2,relative}$ and $H^1_{relative}$ errors of the ST-GFDM are provided in Table.13.  The contours of the numerical solution when $t=0.1$, $t=0.5$ and $t=1$ are shown in Fig. 27, Fig.28 and Fig.29, respectively. From these figures, we can see that the numerical solution changes with the interface deformation and we can capture the interface shape at any time. The numerical solution of the pressure is presented in Fig.30. We can see the distribution of the numerical solution. Therefore, the efficiency of the ST GFDM for solving the Stokes/Parabolic interface problem with a deformation of the eight-petalled flower interface is verified. Furthermore, the last two examples show the advantage in dealing with the deformation of the complex interface. 

\section{Conclusion}
In this paper, the space-time generalized finite difference method is proposed to solve the transient Stokes/Parabolic moving interface problem.  For the ST-GFDM, the proposed problem which is a kind of linearize fluid-structure interaction problem, show that the ST-GFDM has advantages in dealing with the interface conditions, the translation and deformation of the interface, the complex interface shape and the terms of time derivatives. It avoids the problems of time discretization. Five examples verified the accuracy, high efficiency and stability of the ST-GFDM for solving the Stokes/Parabolic moving interface problems, and the proposed method can be extended to deal with a general fluid-structure interaction problem.

\section*{Acknowledgments}

This work is partially supported by Science and Technology Commission of Shanghai Municipality (Grant Nos.  22JC1400900, 22DZ2229014).\\

% To print the credit authorship contribution details
\printcredits
%% Loading bibliography style file
%\bibliographystyle{model1-num-names}
%		\bibliographystyle{cas-model2-names}

% Loading bibliography database
%		\bibliography{}

% Biography
%	\bio{}
% Here goes the biography details.
%	\endbio

%\bio{pic1}
% Here goes the biography details.
%	\endbio

\end{document}